\documentclass[11pt,a4paper]{article}

\usepackage[left=2.54cm,right=2.54cm, top=2.54cm,bottom=2.54cm,bindingoffset=0cm]{geometry}
\usepackage[english]{babel}
\usepackage{graphics}
\usepackage{graphicx}
\usepackage{amsfonts}
\usepackage{amsmath}
\usepackage{amssymb}
\usepackage{amsthm}
\usepackage{float}
\usepackage{url}
\usepackage{bbm}
\usepackage{setspace}
\usepackage{titlesec}
\usepackage{authblk}
\usepackage{hyperref}
\usepackage{subcaption}
\usepackage{xcolor}
\usepackage{fix-cm}
\usepackage[justification=centering]{caption}
\usepackage{algorithm}
\usepackage[noend]{algpseudocode}
\usepackage{mathtools}
\usepackage{orcidlink}

\usepackage{natbib}
 \bibpunct[, ]{(}{)}{,}{a}{}{,}

\DeclareMathOperator\arctanh{arctanh}
 
\titleformat{\section}
  {\normalfont\fontsize{13}{13}\bfseries}{\thesection}{0.5em}{}
  \titleformat{\subsection}
  {\normalfont\fontsize{11}{11}\bfseries}{\thesubsection}{0.5em}{}

  \titlespacing*{\section}
{0pt}{0.9\baselineskip}{0.3\baselineskip}
\titlespacing*{\subsection}
{0pt}{0.9\baselineskip}{0.3\baselineskip}

\usepackage{paralist}

  \pltopsep=\medskipamount
\usepackage{xr} 
\makeatletter
\newcommand*{\addFileDependency}[1]{% argument=file name and extension
  \typeout{(#1)}
  \@addtofilelist{#1}
  \IfFileExists{#1}{}{\typeout{No file #1.}}
}
\makeatother

\newtheorem{Definition}{Definition}
\newtheorem{Theorem}{Theorem}
\newtheorem*{Theorem*}{Theorem}
\newtheorem{Remark}{Remark}
\newtheorem{Example}{Example}
\newtheorem{Lemma}{Lemma}
\newtheorem{Proposition}{Proposition}
\newtheorem{Corollary}{Corollary}

\title{Assortment Optimization Given Basket Shopping Behavior \\ Using the Ising Model}

\date{\vspace{-6ex}}
\author{Andrey Vasilyev\,\orcidlink{0000-0001-5014-537X}\,\(^{a,b,}\)\textsuperscript{1}, Sebastian Maier\,\orcidlink{0000-0001-6162-0906}\,\(^{c,}\)\textsuperscript{1,}\footnote{Corresponding author at: Department of Statistical Science, University College London,  1-19 Torrington Place, London WC1E 7HB, United Kingdom. Email address: s.maier@ucl.ac.uk}\;, Ralf W. Seifert\,\orcidlink{0000-0002-2090-992X}\,\(^{a,d}\)}
\affil{ 
\footnotesize
\(^{a}\)College of Management of Technology, \'Ecole Polytechnique F\'ed\'erale de Lausanne (EPFL), Switzerland \\
\(^{b}\)Imperial College Business School,
Imperial College London, United Kingdom \\
\(^{c}\)Department of Statistical Science, University College London (UCL), United Kingdom \\
\(^{d}\)International Institute for Management Development (IMD), Lausanne, Switzerland}

\allowdisplaybreaks
\begin{document}

\maketitle

\footnotetext[1]{These authors contributed equally to this work.}

\begin{abstract}
In markets where customers tend to purchase baskets of products rather than single products, assortment optimization is a major challenge for retailers. Removing a product from a retailer's assortment can result in a severe drop in aggregate demand if this product is a complement to other products. Therefore, accounting for the complementarity effect is essential when making assortment decisions. In this paper, we develop an analytical framework designed to address this problem. We model customers' choices using the binary Ising model, which is a special case of a Markov random field. This model captures pairwise demand dependencies as well as the individual attractiveness of each product, with binary variables indicating whether products are included in the purchased basket. For estimation, we utilize existing methods developed for the spin formulation of the Ising model and then apply a parameter transformation to express the fitted model in the binary form used for assortment optimization. We show that the assortment optimization problem under the Ising model is NP-hard. We also provide multiple theoretical insights into the structure of the optimal assortments based on the graphical representation of the Ising model, and propose several heuristic algorithms that can be used to obtain high-quality solutions to the assortment optimization problem. Our numerical analysis demonstrates that the developed simulated annealing procedure leads to an expected profit gain of 15\% compared to offering an unoptimized assortment (where all products are included) and around 5\% compared to using a revenue-ordered heuristic algorithm.

\smallskip
\noindent \textit{Keywords:} assortment optimization; retailing; Ising model; multi-purchase choice
\end{abstract}
                                                                                  
\section{Introduction}
\label{2sec:intro}

Optimizing product assortments plays a crucial role in increasing a retailer's profitability by both reducing operational costs and redirecting the demand to more profitable products. A critical aspect of this process is deciding which products to eliminate from an existing assortment, as removing a single product can have a substantial impact on demand allocation for the remaining products. For example, if a retailer decides to discard a low-demand or low-margin product, it may accidentally eliminate an item that attracts customers to shop at this retailer and also buy other products \citep{Timonina1}. The importance of this phenomenon cannot be understated. In 2009, Walmart cut 15\% of its inventory, with dire consequences: ``Sales declined for seven consecutive quarters as shoppers took their entire shopping lists elsewhere. By April 2011, Walmart added 8,500 SKUs back to its mix, an average of 11\% of its products'' \citep{forbes_walm}. This example clearly shows how ignoring basket shopping behavior can result in suboptimal assortment decisions and lead to adverse economic outcomes.

The majority of assortment optimization models assume that each customer purchases at most one
product. Such an approach generally allows researchers to account only for the product substitution effect, but not for the complementarity effect. In other words,
excluding a product from the assortment can only increase the probability of customers purchasing other products instead. Discrete choice models satisfying this property are often referred to as regular choice models (see, e.g., \citealp{Berbeglia1}). However, in many industries where customers tend to buy baskets of products, the complementarity effect cannot be ignored: Eliminating a product from the assortment can reduce the probability of customers purchasing other
products. There is an evident need for suitable methodologies designed to support assortment decision-making in such complex environments. This is clearly reflected in recent interviews with practitioners, which indicate that most retailers are willing to adopt more advanced analytics in the area of assortment planning \citep{Rooderkerk2}.

In this paper, we utilize Markov random fields (MRFs) to model basket shopping behavior. In particular, we provide a new perspective on the Ising model -- which is a special case of an MRF -- by showing that it can be viewed as a multi-purchase choice model, i.e., a choice model under which each customer can select a set of products instead of just one product.  
We consider the Ising model such that each node represents a random variable indicating whether a product belongs to the basket, and an edge between two nodes represents the dependency between the two random variables. The joint probability distribution of such random variables determines the choice probabilities for all possible baskets from a given assortment. Note that we assume that each basket consists of unique products, which is a standard assumption for multi-purchase choice models (see Subsection~\ref{2subsec:ising_choice} for more details). Under this assumption, the term ``basket'' is sometimes replaced with other terms such as ``bundle'' (see, e.g., \citealp{Tulabandhula1}) or ``subset'' (see, e.g., \citealp{Benson1}). However, in this paper we refer to customer choices as baskets since this is intuitive terminology that improves the clarity of exposition.

One of the contributions of this paper is that it highlights the connection between the Ising model and the multivariate logit model (MVL), which is one of the most prominent techniques used for modeling basket shopping behavior. In the classic formulation of the Ising model, the random variables take spin values, i.e., $-1$ or~$1$. However, in the customer-choice model used for assortment optimization in this paper, these random variables are binary, with $x_i=1$ indicating that product $i$ belongs to the purchased basket and $x_i=0$ otherwise. As shown in Subsection~\ref{2subsec:ising_choice}, for a given offered set of products, these two assumptions on the domain of the random variables result in the same choice probabilities up to a parameter transformation. At the same time, the binary values assumption turns the Ising model formulation into the MVL formulation (see Subsection~\ref{2subsec:ising_choice} for more details). The purpose of this paper is thus not to introduce a novel multi-purchase choice model, but rather to establish the link between MRFs and multi-purchase choice models. This link ultimately allows us to leverage methodology developed for the spin formulation of the Ising model while formulating the assortment optimization problem under the binary Ising model.

By building upon the general Ising model, we start with a model that is ``inherently multi-purchase'' and are able to make contributions in several areas of importance to multi-purchase settings. \cite{Bai1} concluded their work with a discussion of six challenging directions for future work: model estimation, hardness results, choice probabilities, complementary multi-purchase events, revenue-ordered assortments, and extensions to additional choice models. The theoretical groundwork discussed above provides us with the means to directly tackle the first four of these and touch on the last two challenges. For example, an effective estimation procedure is needed to obtain a choice probability distribution over all possible shopping baskets. Computing the normalization coefficient (also known as the partition function), which ensures that all probabilities sum up to one, is computationally demanding even for relatively small problems. This makes the task of estimating the parameters of such models extremely challenging. In the case of the spin Ising model, however, we can utilize existing methods -- such as approximate maximum likelihood and Density Consistency (DC) estimation -- that alleviate the computational complexity issue associated with the partition function and make it possible to obtain parameter estimates from historical sales data. Moreover, the DC approach allows for the estimation of the spin Ising model’s parameters using solely closed-form expressions, which makes it particularly attractive for large-scale environments.

Further developing this modeling framework, we address the assortment optimization problem under the binary Ising model. Using the NP-hardness of computing the partition function of the Ising model \citep{Istrail1}, we show that the assortment optimization problem under the binary Ising model is NP-hard by proving NP-hardness of its decision version. We also establish several theoretical results on the structure of the optimal assortments based on the graphical representation of the Ising model which can be used to reduce the dimensionality of the assortment optimization problem. In assortment optimization, being able to estimate the expected revenue for a given assortment is essential, but directly computing such revenues might not be computationally feasible if the customer choice probabilities are defined over all possible baskets. However, these values can be approximated by simulating a large number of customer choices, and the methodology developed for the Ising model proves beneficial for this task as well. In this paper, we employ Gibbs sampling to efficiently generate such samples. Lastly, we develop a customized simulated annealing algorithm to solve the assortment optimization problem, where each candidate solution is evaluated using the aforementioned simulation procedure. We compare this method against several benchmark heuristics, showing that it not only achieves higher expected profits but also results in smaller assortment sizes.

The remainder of this paper is organized as follows: In Section~\ref{2sec:lit}, we review the relevant literature with a primary focus on modeling techniques related to multi-purchase shopping behavior. In Section~\ref{2sec:ising_intro}, we provide the background behind our modeling framework and highlight the connection between the Ising model and the MVL. In Section~\ref{2sec:estim}, we describe two methods for estimating the model parameters, followed by an illustrative example. Section~\ref{2sec:ao_single} is devoted to assortment optimization under the Ising model: We show that this optimization problem is NP-hard, provide theoretical insights into the structure of optimal assortments, and investigate the effect of removing edges on the number of isolated nodes. We then develop several heuristic algorithms to solve the assortment optimization problem and carry out an extensive numerical study. Lastly, in Section~\ref{2sec:concl}, we summarize our findings and discuss future research directions.

\section{Literature Review}
\label{2sec:lit}

Discrete choice modeling is a vast research area that aims to predict customer choices given different sets of alternatives. A common approach to discrete choice modeling suggests that each customer selects at most one alternative. This applies to the multinomial logit model (MNL) introduced by \cite{Luce1959}, which is arguably the most prominent discrete choice model. Despite being widely adopted in practice, choice models that assume a single-item customer shopping behavior are less suitable for industries where customers purchase products primarily in baskets. As a result, a number of models of customers' basket shopping behavior have been developed.
These models can be divided into three categories: Multi-purchase choice models (sometimes referred to as menu or subset selection choice models), multiple-discrete choice models, and multiple discrete-continuous choice models. In multi-purchase choice models, several alternatives can be chosen at the same time, but at most one unit of each alternative can be selected.
Multiple-discrete choice models allow an integer number of units of each of the selected alternatives to be chosen. Finally, in multiple discrete-continuous choice models, noninteger amounts of several alternatives can be selected simultaneously. 
Our focus is on the setting in which customers purchase products in sets, i.e., where each shopping basket consists of unique products. This assumption is fairly justified for many retail industries such as consumer electronics, fashion, toys, etc. It also applies if the goal is to study product categories rather than products themselves. Therefore, theory on multiple-discrete and multiple discrete-continuous choice models is less relevant to our research. We refer the reader to the notable works of \cite{Hendel1}, \cite{Dube1}, \cite{Kim1}, \cite{bhat1} and \cite{bhat2} in those fields.

With regard to multi-purchase choice modeling, one of the most prevalent models is the multivariate logit model (MVL). \cite{Hruschka1} used the MVL to analyze how cross-category sales promotion effects influence purchase probabilities.
The theoretical justification for the MVL was developed by \cite{Russell1}, who derived basket choice probabilities from conditional probabilities of purchasing each product given the purchase decisions related to all
other products. The authors showed that the only joint distribution that is consistent with the specified conditional choice probabilities is the multivariate logistic distribution \citep{Cox1}. \cite{Russell1} applied the MVL to analyze basket purchases in a fairly small setting with four product categories. Subsequently, \cite{Boztug1} proposed a way of applying the MVL in settings of a larger scale. In their paper, the authors first determine prototypical baskets that are used for segmentation of the customer base and then estimate a separate MVL for each customer segment. This approach also makes it possible to account for customer heterogeneity by considering a mixture of MVL models. 

\cite{Song1} presented the idea of leveraging the MVL framework to model customer choices across different product categories assuming that customers choose at most one product within each category (i.e., products in one category are considered to be strong substitutes). Such an extension of the MVL is often referred to as the multivariate MNL, or MVMNL (see, e.g., \citealp{Chen1} and \citealp{Jasin1}). The MVL can thus be viewed as a special case of the MVMNL where each product category contains a single product. It can, however, model an arbitrarily strong substitution effect between any two products by assigning a sufficiently large negative value to the parameter representing the pairwise demand dependency between these products (see Section~\ref{2sec:ising_intro} for parameter definition). Two other noteworthy extensions of the classic MNL model are the multi-purchase MNL (MP-MNL) recently proposed by \cite{Bai1}, and the closely related model of \cite{Abdallah1}.

An interesting model similar to the MVL was proposed by \cite{Benson1}. In their work, the utility of each basket equals the sum of utilities of individual products in that basket plus an optional correction term.  The authors proved that the problem of determining the optimal set of baskets receiving corrective utilities is NP-hard and developed several heuristic algorithms for finding such sets. Overall, the MVL and MVMNL models have been applied in various contexts such as recommendation systems \citep{Moon1} and pricing under competition \citep{Richards1} (see \cite{Jasin1} for a brief overview of different application areas of these models). 

Another classic approach to multi-purchase choice modeling is based on the multivariate probit model (MVP) introduced by \cite{Manchanda1}. In the MVP framework, vectors of unobserved parts of product utilities are assumed to follow the multivariate normal distribution. 
The MVP makes it possible to capture the correlation of customer product preferences across different baskets. However, this modeling approach does not consider the complementarity and substitution effects occurring within one purchased basket. In other words, unlike the MVL, the MVP does not account for the fact that purchasing one product may change the marginal value of adding other products to the same basket \citep{Kamakura1}. 

To the best of our knowledge, there are only three papers that address assortment optimization under the MVL model and its variations. \cite{Tulabandhula1} studied the assortment optimization problem under the BundleMVL-K model -- a version of the MVL in which the size of each basket is restricted by an exogenously given constant $K$. The authors focused on a setting in which each basket consists of at most two products, i.e., $K=2$. They showed that the decision version of the assortment problem under this model is NP-complete. This powerful result might seem somewhat counterintuitive. Indeed, the BundleMVL-2 specified for $n$ products can be viewed as the MNL defined over the set of $n^2$ basket alternatives, meaning that the optimal set of baskets (which is generally not translated into the optimal set of products) can be found in polynomial time. The authors proposed a binary search-based heuristic algorithm to solve the assortment problem and compared it against a mixed-integer programming benchmark, as well as two other heuristic algorithms: a greedy approach and the revenue-ordered heuristic (see, e.g., \citealp{Rusmevichientong1}). 
 
\cite{Chen1} studied the assortment optimization problem under the MVMNL with two product categories, i.e., where the size of each basket is at most two. However, in contrast to the work of \cite{Tulabandhula1}, \cite{Chen1} explicitly separated the product portfolio into two disjoint product categories. The authors proved that the assortment optimization problem in this setting is strongly NP-hard. They proposed the concept of adjusted-revenue-ordered assortments and showed that the assortment with the highest revenue provides a $0.5$-approximation. They also developed a $0.74$-approximation algorithm based on a linear programming relaxation of the assortment optimization problem. Lastly, they considered three extensions of the original setting: with capacity constraints, with generally defined basket prices, and with three product categories. They proved that the assortment optimization problems in all these settings do not admit constant-factor approximation algorithms assuming the Exponential Time Hypothesis.

More recently, \cite{Jasin1} addressed the assortment optimization problem under the MVMNL assuming that the parameters representing cross-category demand dependencies are the same for all products in one category. The authors showed that 
the decision version of the assortment optimization problem is NP-hard even if each category comprises no more than two products and there are no cross-category interaction terms (i.e., the utility of each basket is the sum of utilities of the corresponding product categories). They developed a fully polynomial-time approximation scheme (FPTAS) that can be used to solve this optimization problem. The authors also proposed the generalization of the FPTAS for the case with nonzero cross-category interactions assuming that the maximum number of interacting categories is two. Finally, they considered a setting with capacity constraints and developed an FPTAS for the capacitated assortment problem by building on the ideas presented for the uncapacitated case.

Unlike existing works in the multi-purchase literature, which have made important contributions under specific modeling assumptions, our Ising-based formulation allows shopping baskets to contain any possible combination of the $n$ available products (i.e., any of the $2^n$ possible baskets) without imposing artificial restrictions on customer choices. We further illustrate the value of this generality and the potential cost of model misspecification by comparing Ising-based assortments with assortments derived from MNL and BundleMVL-2 benchmarks.
As detailed in Section~\ref{2sec:intro}, leveraging the powerful framework of the Ising model has also enabled us to advance existing research on assortment optimization under basket shopping behavior by adopting methodologies and theoretical results developed for the inherently multi-purchase Ising model. The necessary theoretical background related to the Ising model is provided in the following section.

\vspace{-0.04cm}
\section{Model Formulation}
\label{2sec:ising_intro}

In this section, we recall the definition of a Markov random field (MRF) and show that the Ising model can be viewed as a multi-purchase choice model.

\subsection{Background on Markov Random Fields}
\label{2subsec:mrf_background}

The Ising model was introduced by \cite{Lenz1} and studied in a one-dimensional case in the work of \cite{ising}. It marked the beginning of the development of graphical models, including MRFs. For clarity, 
we start by providing a formal definition of an~MRF.
\begin{Definition}
\label{MRF}
Let $\xi = \{\xi_i\}_{i \in \{1, \dots, n\}}$ be a random field represented by an undirected graph  $G = (V, E)$, where $V = \{v_i\}_{i \in \{1, \dots, n\}}$ is the set of nodes such that each node $v_i$ is associated with a random variable $\xi_i$, and $E$ is the set of edges. Then, the random field $\xi$ is called a Markov random field if the joint probability distribution $p(\xi)$ is positive and satisfies the conditional independence property, i.e., each random variable $\xi_i$ is independent of all other variables given its~neighbors.
\end{Definition}

The MRF model received major attention in the 1970s when the equivalence relation between this model and the Gibbs distribution was established. To define the Gibbs distribution, recall that a subset of nodes of a graph is called a clique if the corresponding subgraph is complete.
Let~$\mathfrak{L}$ be the set of all cliques of graph $G$. Suppose that for each clique $l \in \mathfrak{L}$, there exists a strictly positive function $\psi_l$ referred to as the potential function of this clique. Then, the Gibbs distribution can be defined in the following way~(see, e.g., \citealp{Murphy1}):

\begin{Definition}
A joint probability distribution $p$ defined over the nodes of graph~$G$ is called a Gibbs distribution if it can be expressed in terms of potential functions of cliques of graph~$G$ in the following way:
\begin{align*}
&p(x) = \frac{1}{Z} \prod_{l \in \mathfrak{L}} \psi_{l}(x),
\end{align*}
where $Z = \sum_{x} \prod_{l \in \mathfrak{L}} \psi_l(x)$ is the normalization coefficient (also known as the partition function).
\end{Definition}

The connection between MRFs and the Gibbs distribution was first studied in the works of \cite{dobruschin} and \cite{spitzer}. Subsequently, the equivalence relation between these two concepts was proved in an unpublished paper by \cite{hammersley_clifford}. This important result was later revisited in the work of \cite{besag_1974}, where the author restated the Hammersley-Clifford theorem and gave an alternative proof of it. This theorem can be formulated as follows \citep{Hristopulos1}:
\begin{Theorem*}[Hammersley-Clifford]
A random field $\xi$ is a Markov random field if and only if the joint probability distribution $p(\xi)$ is a Gibbs distribution.
\end{Theorem*}
\noindent Note that the Hammersley-Clifford theorem establishes the equivalence between local properties of a random field (the conditional independence properties) and its global property (being defined by the Gibbs distribution).

The Ising model is a prototypical example of an MRF. It originates from statistical mechanics and is often defined on a lattice rather than a general graph. In our case, however, we do not assume any lattice structure in the graph representation to allow for modeling pairwise demand dependencies for any pair of products. The Ising model can be defined by the joint distribution of random variables $\xi = \{\xi_i\}_{i \in \{1, \dots, n\}}$. Let us now specify the probability mass function of this distribution.

Let $\mathcal{N}$ denote the set of indices $\{1, \dots, n\}$. Suppose that each $\xi_i$, $i \in \mathcal{N}$ is a binary random variable. For the Ising model, the probability
mass function is defined in the following way:
\begin{equation}
\label{pmf_full}
p_{\theta}(x | \mathcal{N}) = \exp\biggl(\sum_{i \in \mathcal{N}}\theta_{ii}x_i + \sum_{i \ne j, \ i, j \in \mathcal{N}}x_i\theta_{ij}x_j - A_{\theta}(\mathcal{N})\biggr),
\end{equation}
where $x \in \mathcal{X}(\mathcal{N}) = \{0, 1\}^{|\mathcal{N}|}$ is a realization of the random vector $\xi$, $\theta \in \mathbb{S}^n$ is a symmetric matrix of parameters,
and $A_{\theta}(\mathcal{N})$ is the logarithm of the partition function given by
\begin{equation}
\label{norm_coef_full}
A_{\theta}(\mathcal{N}) = \log\Biggl(\sum_{x \in \mathcal{X}(\mathcal{N})}\exp\biggl(\sum_{i \in \mathcal{N}}\theta_{ii}x_i + \sum_{i \ne j, \ i, j \in \mathcal{N}}x_i\theta_{ij}x_j \biggr)\Biggr).
\end{equation}
Clearly, distribution~(\ref{pmf_full}) is a special case of the Gibbs distribution. Therefore, from the Hammersley-Clifford theorem, it follows that the Ising model is an MRF. The conditional independence property can also be verified directly, thereby confirming that the Ising model is an MRF. 
Importantly, the Ising model can define very complex probability distributions. In fact, any pairwise MRF -- i.e., a graphical model defined by expressions similar to~(\ref{pmf_full}) and~(\ref{norm_coef_full}) but with the sums $\sum_{i}\theta_{ii}x_i$ and $\sum_{i \ne j}x_i\theta_{ij}x_j$ replaced by  $\sum_{i}E_i(x_i)$ and $\sum_{i \ne j} E_{ij}(x_i, x_j)$, respectively, where $E_i$ and $E_{ij}$ are arbitrary real-valued functions -- is equivalent to the Ising model with one extra node (see \citealp{Globerson1} and \citealp{Schraudolph1}). We leave the question of utilizing more general MRFs or the Ising model with an additional node in the context of revenue management for future research.

Now, having provided all the necessary theoretical background, we can offer a novel perspective on the Ising model as a multi-purchase choice model.

\subsection{Ising Model as a Multi-Purchase Choice Model}
\label{2subsec:ising_choice}

As mentioned in Section~\ref{2sec:lit}, we consider the setting in which each basket consists of unique products, i.e., it is assumed that customers do not purchase more than one unit of each product. Then, basket $\{i_1, \dots, i_k\} \subseteq \mathcal{N}$ can be represented by vector $x = \{x_1, \dots, x_n\}$ such that
\begin{equation*}
\begin{aligned}
x_i = &\begin{cases}
1&\text{if $i \in \{i_1, \dots, i_k\}$},\\
0&\text{otherwise.}
\end{cases}
\end{aligned}
\end{equation*}
We can thus consider the Ising model such that each node represents a binary random variable corresponding to a certain product, and an edge between two nodes
represents the dependency between the corresponding random variables. The binary value of each random variable indicates whether the corresponding product belongs to the basket.
Then, each parameter $\theta_{ij}$, $i \ne j$ represents the pairwise demand dependency between products $i$ and $j$, and parameter $\theta_{ii}$ determines the individual attractiveness of product~$i$. Consequently, the joint probability distribution of these random variables defines the probability of a random customer selecting any basket, i.e., any subset of products from $\mathcal{N}$. Finally, formulas~(\ref{pmf_full}) and~(\ref{norm_coef_full}) can be straightforwardly adjusted for the case when only a subset $S \subseteq \mathcal{N}$ of products is offered, thereby defining a multi-purchase choice model:
\begin{subequations}
\label{ising_choice}
\begin{align}
p_{\theta}(x | S) &= \exp\biggl(\sum_{i \in S}\theta_{ii}x_i + \sum_{i \ne j, \ i, j \in S}x_i\theta_{ij}x_j - A_{\theta}(S)\biggr), \label{pmf_ising}\\
A_{\theta}(S) &= \log\Biggl(\sum_{x \in \mathcal{X}(S)}\exp\biggl(\sum_{i \in S}\theta_{ii}x_i + \sum_{i \ne j, \ i, j \in S}x_i\theta_{ij}x_j \biggr)\Biggr), \label{norm_coef_ising}
\end{align}
\end{subequations}
where $x$ represents a basket, and $\mathcal{X}(S)$ is the set of all possible baskets of products from~$S$.

Importantly, in the classic formulation of the Ising model, random variables are assumed to be spin variables rather than binary variables, i.e., they are assumed to take values $-1$ or~$1$. However, for a given offered set, these two assumptions on the domain of the random variables result in models that are equivalent to each other up to a parameter transformation. In other words, one can construct a parameter transformation that maps the choice probabilities given by one model onto the choice probabilities given by another model, and vice versa. 
\begin{Proposition}
\label{prop_transf_equiv}
Suppose that $\mathcal{M}_1$ is the Ising model with parameters $\theta$ defined over binary variables, and $\mathcal{M}_2$ is the Ising model with parameters $\tilde{\theta}$ defined over spin variables. Let $\mathcal{M}_1$ and $\mathcal{M}_2$ define the choice probabilities of customers selecting baskets from a fixed assortment~$\mathcal{N}$. If $\theta$ and $\tilde{\theta}$ satisfy the following conditions:
\begin{equation}
\label{par_transform}
\begin{aligned}
\theta_{ii} &= 2\tilde{\theta}_{ii} - 4\sum\limits_{j \in \mathcal{N}}: \, j \ne i  \tilde{\theta}_{ij} &\forall i \in \mathcal{N}, \\
\theta_{ij} &= 4\tilde{\theta}_{ij} &\forall i, j \in \mathcal{N}: \ i \ne j,
\end{aligned}
\end{equation}
then the choice probabilities yielded by models $\mathcal{M}_1$ and $\mathcal{M}_2$ are identical.
\end{Proposition}
\noindent The formal proof of this proposition can be found in \ref{prop_transf_equiv_proof}. If the variables in the Ising model take binary values instead of spin values, then for a given choice set the Ising model takes exactly the same functional form as the MVL. Proposition~\ref{prop_transf_equiv} allows us to extend this equivalence relation to the Ising model defined over spin variables. Formally, we can state the following corollary of Proposition~\ref{prop_transf_equiv}:
\begin{Corollary}
\label{cor_ising_mvl}
The following relation between the Ising model and the MVL holds:
\begin{itemize}
    \item For a given offered set, the binary Ising model is equivalent to the MVL;
    \item For a given offered set, the spin Ising model is equivalent to the MVL up to a parameter transformation.
\end{itemize}
\end{Corollary}
It should be emphasized that Proposition~\ref{prop_transf_equiv} establishes the equivalence between the spin and binary formulations for a \textit{fixed} offered set. In the framework developed below, this result is used only for the full product set $\mathcal{N}$ observed in the historical sales data: after estimating spin parameters under $\mathcal{N}$, we transform them into the corresponding binary parameters. The fitted choice model used for assortment optimization is then the binary Ising model in~(\ref{ising_choice}). Once this transformation has been applied, choice probabilities for candidate assortments $S \subseteq \mathcal{N}$ are defined directly by this binary model.

The binary formulation is natural for the choice model considered in this paper: it uses the standard $0/1$ representation of basket membership in choice modeling and assortment optimization and it makes the connection with the MVL model transparent. At the same time, the spin formulation is standard in much of the Ising-model literature, and the estimation methods described in Section~\ref{2sec:estim} are stated in terms of spin variables. Proposition~\ref{prop_transf_equiv} allows us to use these existing estimation methods directly, rather than rederiving them in binary notation, which would add unnecessary technical complexity without changing the substance of the framework.

As previously stated, the purpose of this paper is not to introduce a novel multi-purchase choice model, but to underline the connection (in a certain sense, the equivalence relation) between the Ising model and the MVL. In Section~\ref{2sec:lit}, we emphasized that the MVL model is one of the most prominent multi-purchase choice models.
This model can be theoretically justified in two different ways based on the random utility theory.
First, one can specify the conditional utility of purchasing each product given the purchase decisions related to all other products in the assortment and derive the corresponding conditional choice probabilities. Then, the factorization theorem of~\cite{besag_1974} can be used to obtain the unique joint distribution consistent with
this conditional probability distribution. 
We refer the reader to \cite{Russell1}, who developed this theoretical justification for the MVL model, and to \cite{Tulabandhula1}, who adapted this approach to the BundleMVL-K model and provided a detailed description of it. Alternatively, one can assume that the utility of a basket represented by vector $x$ is given by the following expression:
\begin{equation*}
    U_\theta(x|S) = \sum_{i \in S}\theta_{ii}x_i + \sum_{i \ne j, \ i, j \in S}x_i\theta_{ij}x_j + \epsilon_{x},
\end{equation*}
where $\epsilon_x$ is a random variable representing the unobserved part of the utility. If variables~$\epsilon_x$ are independent and identically distributed (i.i.d.) for all possible baskets and normalized so that their mean is zero and the variance is $\pi^2/6$, then the MVL can be viewed as an extended version of the MNL, where the set of alternatives is given by the set of all possible baskets and the resulting choice probabilities are given by expressions~(\ref{ising_choice}). This theoretical justification approach was pursued by, for example, \cite{Jasin1} for the MVMNL model. The independence assumption might be criticized as being too restrictive since it must hold even for baskets that share some of the products (see, e.g., \citealp{Tulabandhula1}). However, this assumption is milder than it might appear. Note that the independence assumption applies to the MNL model, despite the fact that some alternatives (i.e., individual products) can share similar attributes. In fact, as \cite{Train1} pointed out with respect to the MNL model, this assumption can be viewed as a natural consequence of a well-specified model. It implies that for each basket, the constant part of utility is specified sufficiently so that the unobserved part of the utility does not provide any information about the unobserved parts of utilities of other alternatives. Therefore, this is an assumption on, above all, the quality of the model specification. To summarize, each of the two described theoretical justifications is well grounded and can be viewed as the utility theory foundation underlying the considered model.

Since the Ising model was introduced long before the MVL, we will primarily use the former term when referring to customer choice probabilities. Using such terminology also highlights the fact that most of the methodology utilized in this paper was originally developed for the Ising model, which includes parameter estimation.

\section{Parameter Estimation and Illustrative Example}
\label{2sec:estim}

In this section, we describe two methods for estimating the Ising model parameters. We begin with a more classical approach that employs approximate sparse maximum \text{likelihood\,estimation.}

\subsection{Approximate Sparse Maximum Likelihood Estimation}

Since we utilize the methodology developed for the classic version of the Ising model -- with variables taking spin values -- we first describe how to obtain estimates of parameters~$\tilde{\theta}$. The estimates of parameters~$\theta$ for the binary Ising model can then be obtained by applying transformation~(\ref{par_transform}) to~$\tilde{\theta}$. The use of spin variables in this section is therefore limited to the estimation step; after the transformation, the fitted model is expressed in the binary form used for assortment optimization. We assume that the only sales history information that is available is the list of purchased baskets assuming that all products are offered, i.e., $S = \mathcal{N}$. Such an assumption makes the parameter estimation method described in this section readily applicable in a range of practical situations and is fairly standard in the literature (see, e.g., \cite{Tulabandhula1,Jasin1}). Note that the assortment optimization problem studied later in this paper can thus be viewed as one of finding a subset of the current product portfolio (i.e. $S \subseteq\mathcal{N}$) that maximizes the expected profit. Finally, since the assortment is fixed throughout this section, we will slightly abuse notation by omitting the argument $\mathcal{N}$ of the log-partition function $\tilde{A}_{\tilde{\theta}}\coloneq\tilde{A}_{\tilde{\theta}}(\mathcal{N})$, which is defined as:
\begin{equation*}
    \tilde{A}_{\tilde{\theta}}(\mathcal{N}) = \log\biggl(\sum_{\tilde{x} \in \{-1,1\}^{\vert\mathcal{N}\vert}}\exp\Bigl(\sum_{i \in \mathcal{N}}\tilde{\theta}_{ii}\tilde{x}_i + \sum_{i \ne j, \ i, j \in \mathcal{N}}\tilde{x}_i\tilde{\theta}_{ij}\tilde{x}_j \Bigr)\biggr)\,.
\end{equation*}
Importantly, unlike the log-partition function in (\ref{norm_coef_full}), the expression above for $\tilde{A}_{\tilde{\theta}}$ uses parameters~$\tilde{\theta}$ for the spin Ising model but also a vector $\tilde{x}$ of spin random variables, i.e. $\tilde{x}\in\{-1,1\}^{\vert\mathcal{N}\vert}$.

Suppose that the available historical data sample comprises a list of purchased baskets  $B = \{b^k_i\}^{i = 1, \dots, n}_{k = 1, \dots, m}$, where each basket $b^k$ is represented as a spin vector with $b^k_i = 1$ if product $i\in \mathcal{N}$ belongs to basket $b_k$, otherwise $b^k_i = -1$.
Then, under the Ising model with parameters~$\tilde{\theta}$, the negative mean log-likelihood of selecting baskets from the given sample is as follows:
\begin{equation}
\label{neg_mean_ll}
-LL_{mean}(\tilde{\theta}) = \tilde{A}_{\tilde{\theta}} - \biggl(\sum_{i \in \mathcal{N}}\tilde{\theta}_{ii}\mu_i + \sum_{i \ne j, \ i, j \in \mathcal{N}}\tilde{\theta}_{ij}s_{ij}\biggr),
\end{equation}
where  $\mu = \frac{1}{m} \sum^m_{k=1}{b^k}$ is the sample mean
vector, and $s =  \frac{1}{m} \sum^m_{k=1}{b^k(b^k)^T} =  \frac{1}{m} B^T B$
is the sample second-order moment matrix. The model parameters can be learned by minimizing the right-hand side of~(\ref{neg_mean_ll}) over $\tilde{\theta}  \in \mathbb S^n$. Moreover, the matrix of off-diagonal entries of $\tilde{\theta}$ can be sparsified by using $\ell_1$-regularization, i.e., by adding a $\ell_1$-norm penalty term to the corresponding optimization problem. This allows us to focus on the most important pairwise dependencies while removing nonmeaningful interactions from consideration. Formally, sparse estimates of parameters $\tilde{\theta}$ can be obtained by solving the following problem:
\begin{equation}
\label{opt_est}
\min_{\tilde{\theta} \in \mathbb S^n} \tilde{A}_{\tilde{\theta}} - \biggl(\sum_{i \in \mathcal{N}}\tilde{\theta}_{ii}\mu_i + \sum_{i \ne j, \ i, j \in \mathcal{N}}\tilde{\theta}_{ij}s_{ij}\biggr) + \rho ||\tilde{\theta} - \text{diag}(\tilde{\theta})||_1,
\end{equation}
where 
$\rho\geq 0$ is the weight of the penalty term. The greater the value of $\rho$, the fewer nonzero off-diagonal entries in matrix $\tilde{\theta}$ and the fewer edges in the corresponding network. Note, we penalize only the off-diagonal entries of $\tilde{\theta}$ because including diagonal ones would lead to a further reduction in the quality of fit but not provide any additional benefits in terms of model sparsity.

The log-partition function $\tilde{A}_{\tilde{\theta}}$ contains an exponential number of terms, which makes optimization problem~(\ref{opt_est}) extremely computationally challenging. In fact, even the problem of computing the partition function for a given $\tilde{\theta}$ is NP-hard in general, 
as proven by \cite{Barahona1}. \cite{Istrail1} extended the work of \cite{Barahona1} and showed that the computational complexity of this problem arises from the topology of the underlying graphical model: If the corresponding graph is nonplanar, then the problem is NP-hard. Since we do not make any restrictive assumptions on the topological structure of the Ising model, the presence of the \text{(log-)} partition function poses a challenging problem. It can be overcome by a notable result obtained for the Ising model, as illustrated below.  \cite{Wainwright1} showed that the log-partition function can be upper bounded by the solution to a certain convex optimization problem, and \cite{Banerjee1} proved that this upper bound can be rewritten as follows: 
\vspace{-0.2cm}
\begin{equation}
\label{ising_up_bound}
\widehat{\tilde{A}}_{\tilde{\theta}} = \dfrac{n}{2} \log\Bigl(\dfrac{e\pi}{2}\Bigr) - \dfrac{1}{2}(n+1) - \dfrac{1}{2}\biggl(\max_{v \in \mathbb{R}^{n+1}} v^T q + \log\det{\bigl(-Q(\tilde{\theta}) - \text{diag}(v) \bigr)}\biggr),
\end{equation}
where $q = (1, 4/3, \dots, 4/3)^T \in \mathbb{R}^{n+1}$ and 
\vspace{-0.2cm}
\begin{equation*}
Q(\tilde{\theta}) =
\begin{pmatrix}
    0       & \tilde{\theta}_{11} & \tilde{\theta}_{22} & \dots & \tilde{\theta}_{nn} \\
    \tilde{\theta}_{11} & 0 & 2\tilde{\theta}_{12} & \dots & 2\tilde{\theta}_{1n} \\
    \tilde{\theta}_{22} & 2\tilde{\theta}_{21}  & 0 & \dots & 2\tilde{\theta}_{2n} \\
    \vdots &  \vdots  &  \vdots &  \ddots &  \vdots \\
    \tilde{\theta}_{nn} & 2\tilde{\theta}_{n1}  & 2\tilde{\theta}_{n2} & \dots & 0 \\
\end{pmatrix}.
\end{equation*}
Thus, the model parameters can be estimated by solving the following optimization problem:
\begin{equation}
\label{opt_est_bar}
\min_{\tilde{\theta} \in \mathbb S^n} \widehat{\tilde{A}}_{\tilde{\theta}} - \biggl(\sum_{i \in \mathcal{N}}\tilde{\theta}_{ii}\mu_i + \sum_{i \ne j, \ i, j \in \mathcal{N}}\tilde{\theta}_{ij}s_{ij}\biggr) + \rho ||\tilde{\theta} - \text{diag}(\tilde{\theta})||_1.
\end{equation}
Optimization problem~(\ref{opt_est_bar}) is convex and can be solved using one of the standard commercial or open-source solvers. We used SDPT3 v4.0 -- a MATLAB package for semidefinite-quadratic-linear programming (see \citealp{Toh1} and \citealp{Tutuncu1}). Once optimal values of parameters~$\tilde{\theta}$ are identified, we can obtain the values of parameters $\theta$ for the binary Ising model by applying transformation~(\ref{par_transform}) to $\tilde{\theta}$.

To illustrate the practical relevance of this estimation method, we apply it to an open-source dataset. In particular, we choose the {\scshape Bakery} dataset -- utilized by \cite{Benson1} to evaluate the performance of their model -- on account of its moderate size, which makes it very well suited to visualizing results.
The dataset comprises a list of selected baskets obtained from the receipts of purchases by customers of a bakery. It was preprocessed so that it contains baskets with less than 6 products and only those products that were selected at least 25 times, resulting in 67,488 baskets and $n=50$ products in total. 

We randomly split this dataset into training and test samples at a ratio of 80:20. Using the training sample, we estimate parameters $\tilde{\theta}$ by solving problem~(\ref{opt_est_bar}) for several values of the penalty weight, $\rho \in \{0, 0.0025, 0.005, 0.01, 0.015, 0.02, 0.03\}$. The value used for likelihood-based evaluation is selected by five-fold cross-validation on the training sample, with mean validation log-likelihood per basket as the selection criterion. This procedure selects $\rho = 0$.

The unregularized estimate, however, yields an almost complete graph, which is difficult to interpret visually. To visualize the main dependency structure in the product network, we therefore also estimate a sparse model with $\rho = 0.015$, chosen from the cross-validation grid as a representative value that substantially simplifies the graph structure while retaining a reasonable validation likelihood. After refitting on the full training sample, this value reduces the number of nonzero off-diagonal entries of $\theta$ from 1225 to 137, where an edge is counted as nonzero if $|\theta_{ij}| > 10^{-5}$. Thus, the sparse estimate keeps only a small subset of the pairwise interactions and yields a compact graph of the strongest dependencies, making the fitted model more attractive for interpretation, storage, implementation, and downstream optimization. The corresponding loss in test likelihood is reported below. Lastly, we convert $\tilde{\theta}$ to $\theta$ by applying transformation~(\ref{par_transform}) and construct a graph based on the values of~$\theta$. We visualize this sparse graph in the following way (see Figure~\ref{network_vis_sparse}):
\begin{figure}[H]
\centering
    \includegraphics[width=0.9\textwidth]{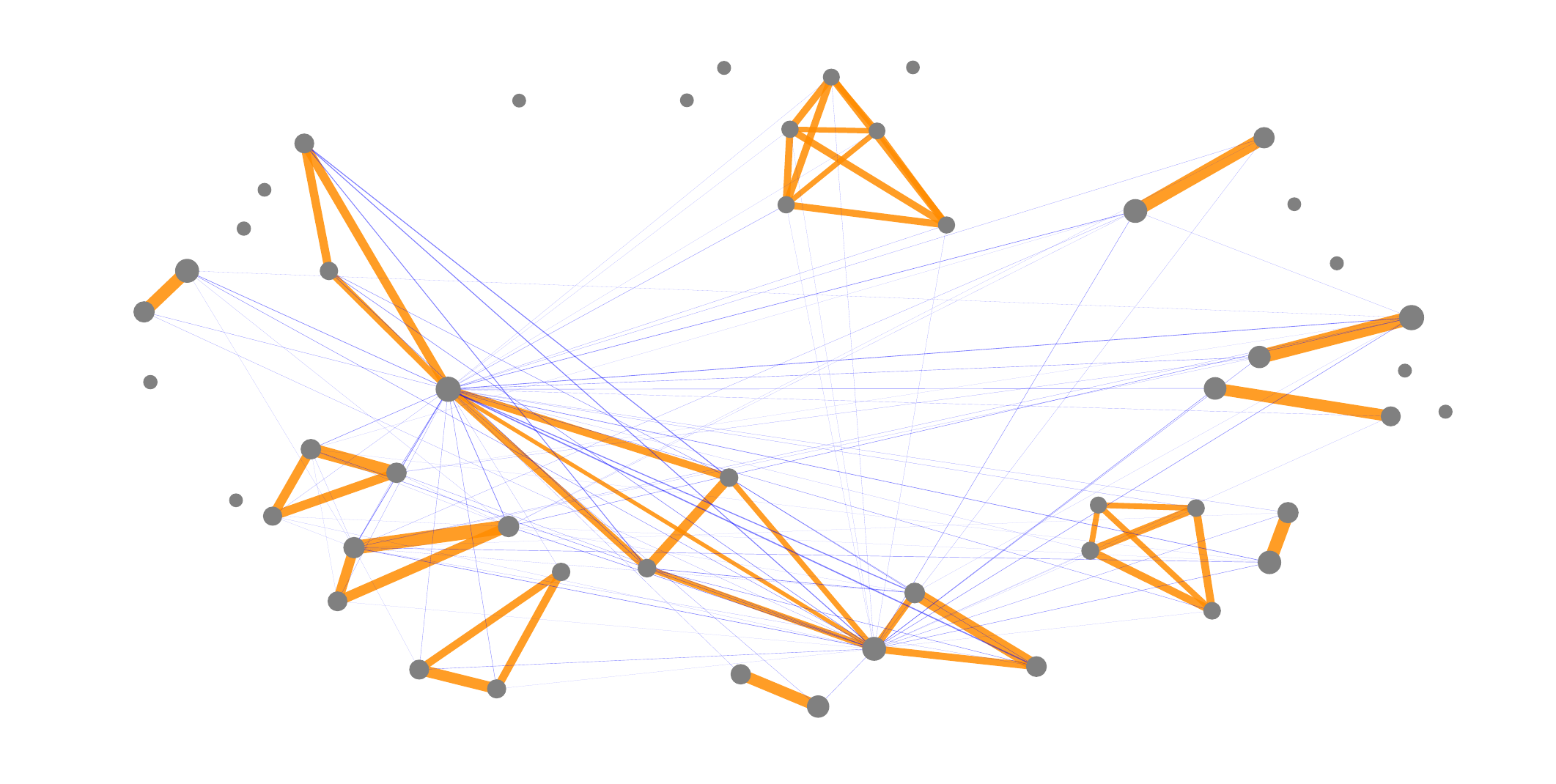}
    \caption{Sparse graphical representation of the Ising model for the {\scshape Bakery} dataset.}
    \vspace{-0.1cm}
    \label{network_vis_sparse}
\end{figure}
\begin{itemize}
\item
  The greater the value of $\theta_{ii}$, the larger the size of node
  $i$;
\item
  The greater the absolute value of $\theta_{ij}$, the thicker the edge
  $(i, j)$;
\item
  If $\theta_{ij} < 0$, then the corresponding edge is blue
  (meaning a negative direct dependency between $i$ and $j$), and if
  $\theta_{ij} > 0$, then the corresponding edge is orange
  (meaning a positive direct dependency between $i$ and $j$).
\end{itemize}
\noindent Such an illustration provides valuable insights into the structure of the product portfolio. For example, one can easily identify ``established'' baskets of products (represented by groups of nodes mutually connected by edges with positive weights), or independent-demand products (represented by isolated nodes).

\subsection{Density Consistency Estimation}

Another powerful method for estimating the parameters of the spin Ising model is the Density Consistency (DC) approach proposed by \cite{Braunstein1}. It allows for the estimation of the Ising model's parameters using solely closed-form expressions that are functions of the first and second empirical moments, which makes it especially attractive in large-scale environments. This approach is closely related to the expectation propagation algorithm and employs a refined Gaussian approximation with a modified consistency condition \citep{Braunstein2}.
The DC estimates of parameters $\tilde{\theta}$ are as follows:
\begin{equation*}
\begin{aligned}
\tilde{\theta}^{DC}_{ii} &= \tilde{\theta}^{IP}_{ii} + (n - 2) \arctanh{\mu_i} - \sum\limits_{j \ne i} \frac{\Sigma_{jj} \mu_i - \Sigma_{ij} \mu_j}{\Sigma_{ii} \Sigma_{jj} - \Sigma^2_{ij}} + (\Sigma^{-1} \mu)_i, \\
\tilde{\theta}^{DC}_{ij} &= \frac{1}{2} \biggl(\tilde{\theta}^{IP}_{ij} - \Sigma^{-1}_{ij} - \frac{\Sigma_{ij}}{\Sigma_{ii} \Sigma_{jj} - \Sigma^2_{ij}}\biggr),
\end{aligned}
\end{equation*}
where
\begin{equation*}
\begin{aligned}
\tilde{\theta}^{IP}_{ii} &= -(n - 2) \arctanh{\mu_i} + \dfrac{1}{4} \sum\limits_{j \ne i} \log{\dfrac{((1 + \mu_i) (1 + \mu_j) + C_{ij})  ((1 + \mu_i) (1 - \mu_j) - C_{ij})}{((1 - \mu_i) (1 + \mu_j) - C_{ij})  ((1 - \mu_i) (1 - \mu_j) + C_{ij})}}, \\
\tilde{\theta}^{IP}_{ij} &= \dfrac{1}{4} \log{\dfrac{((1 + \mu_i) (1 + \mu_j) + C_{ij})  ((1 - \mu_i) (1 - \mu_j) + C_{ij})}{((1 + \mu_i) (1 - \mu_j) - C_{ij})  ((1 - \mu_i) (1 + \mu_j) - C_{ij})}}, \\
\mu &= \frac{1}{m} \sum^m_{k=1}{h^k}, \quad C = s - \mu \mu^T, \quad s = \frac{1}{m} B^{T}  B, \\
\Sigma_{ii} &= \dfrac{\mu_i}{\arctanh{\mu_i}}, \quad \Sigma_{ij} = C_{ij}  \sqrt{\dfrac{\Sigma_{ii}\Sigma_{jj}}{C_{ii}  C_{jj}}} \qquad \forall i,j = 1, \dots, n.
\end{aligned}
\end{equation*}

In the above expressions, $\tilde{\theta}^{IP}_{ii}$ and $\tilde{\theta}^{IP}_{ij}$ are the so-called independent-pair estimates of the Ising model parameters, which rely on a simple approximation technique that assumes that each pair of nodes is independent of the others (see \citealp{Roudi1} for more detail). As before, $B$ is the matrix containing samples, $\mu$ is the sample mean vector (empirical first moments), and $s$ is the sample second-order moment matrix (empirical second moments). Finally, $C$ is the empirical covariance matrix.
\begin{figure}[t]
\centering
\begin{subfigure}{.49\textwidth}
  \centering
  \includegraphics[width=\linewidth]{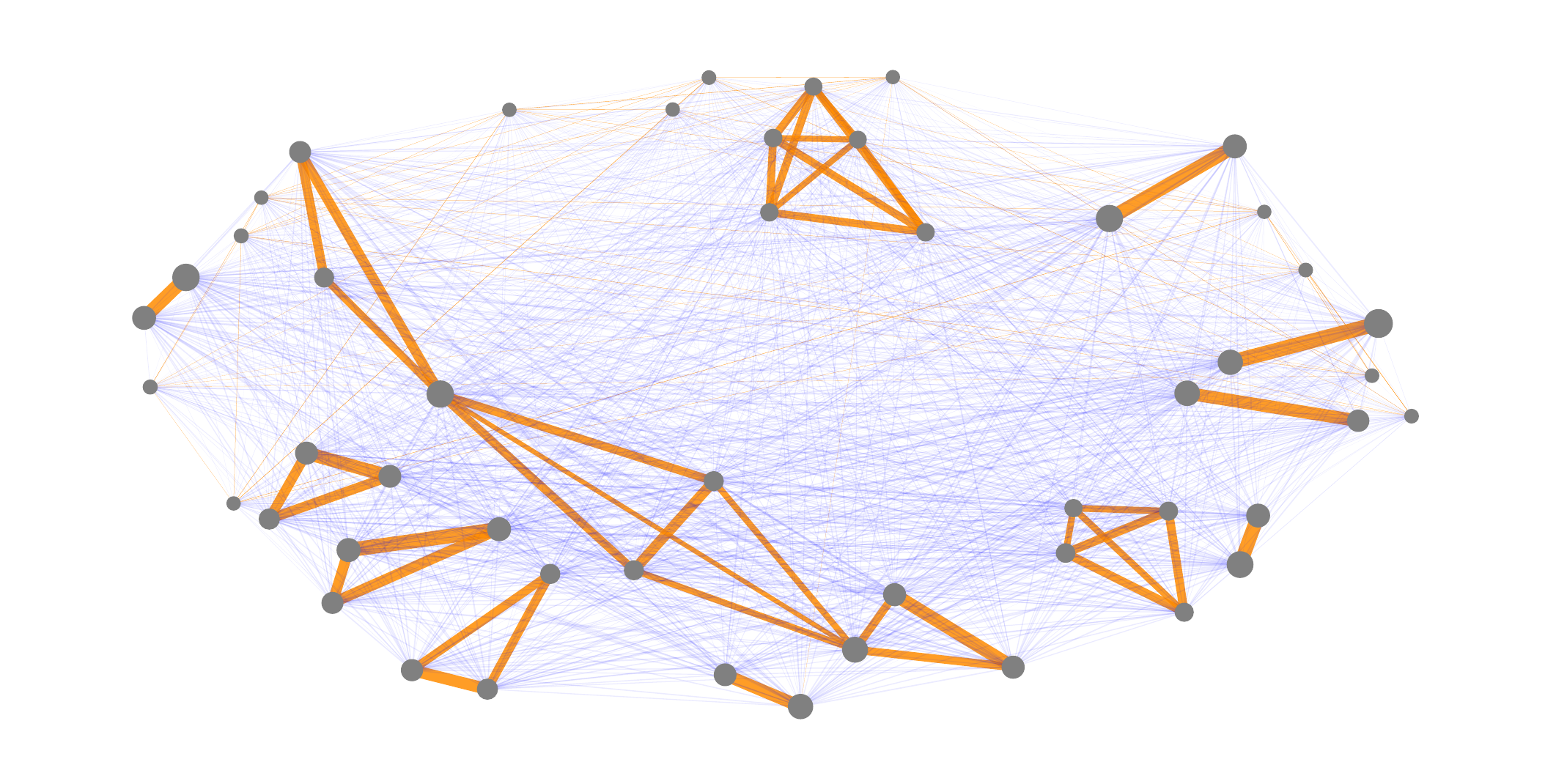}
   \caption{\footnotesize{Approximate maximum likelihood estimation.}}
\end{subfigure}
\begin{subfigure}{.49\textwidth}
  \centering
  \includegraphics[width=\linewidth]{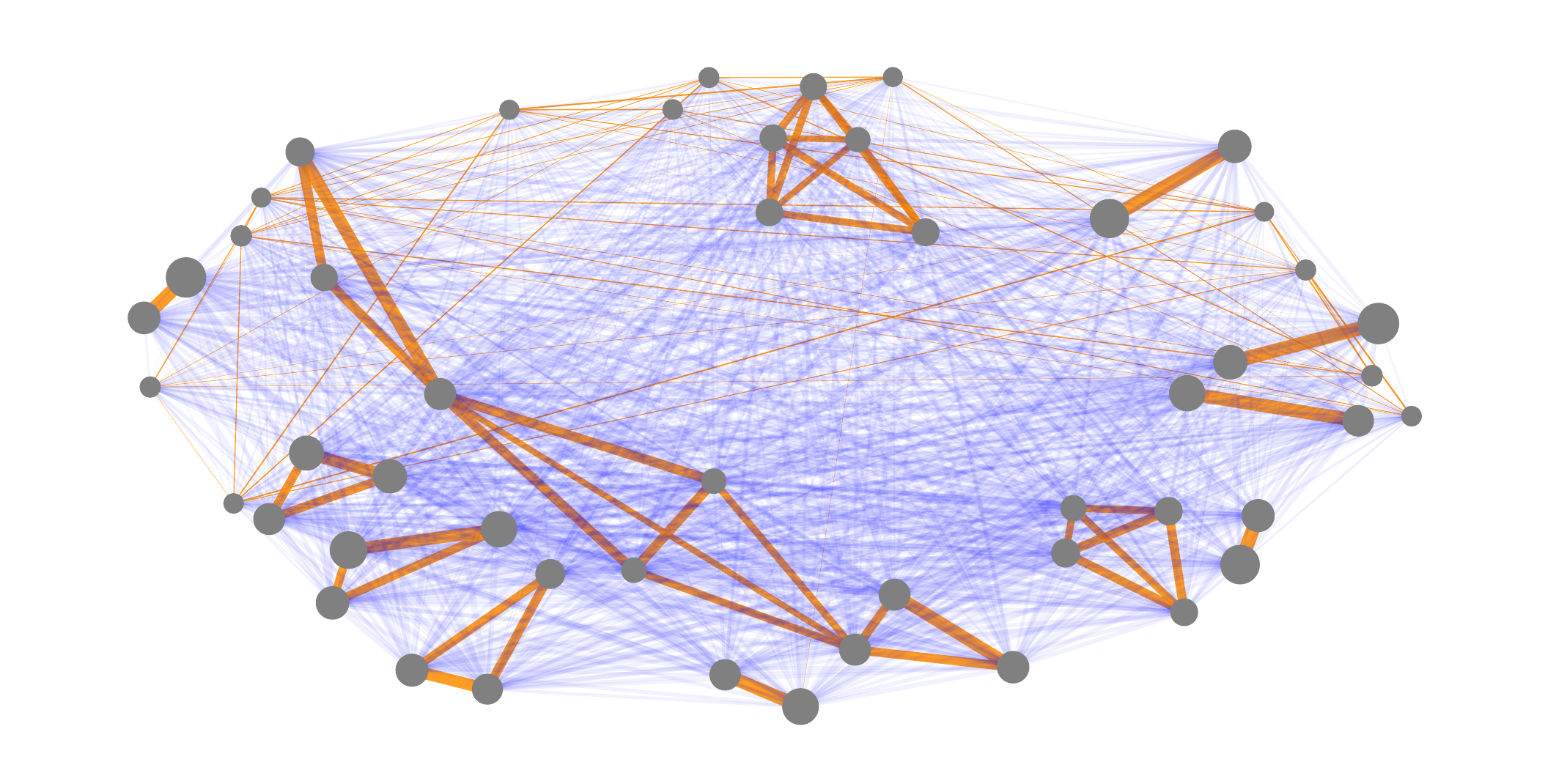}
   \caption{\footnotesize{Density Consistency estimation.}}
\end{subfigure}
\caption{Ising model parameters estimated using (a) approximate maximum likelihood estimation and (b) Density Consistency estimation for the {\scshape Bakery} dataset.}
\vspace{-0.1cm}
\label{network_vis_two}
\end{figure}

The DC approach to estimating the Ising model parameters performs surprisingly well. \cite{Braunstein1} compared this method with several state-of-the-art inference techniques across a wide range of generated instances. The comparison showed that while none of the considered methods is universally better than the others, DC consistently ranks among the most accurate and reliable methods. When applied to the {\scshape Bakery} dataset, DC infers a graph structure that is perfectly aligned with the one yielded by the approximate maximum likelihood estimation described in the previous subsection. Figure~\ref{network_vis_two}(a) visualizes the Ising model parameters estimated using the approximate maximum likelihood estimation with the sparsity-controlling parameter $\rho$ set to zero, while Figure~\ref{network_vis_two}(b) visualizes the parameters estimated using the DC approach.
It can be observed that DC assigns slightly higher weights to negative dependencies and individual product parameters $\tilde{\theta}_{ii}$, yet the overall graph structures are very similar.

\subsection{Empirical Performance of Ising Model}

To evaluate the performance of the Ising model as a multi-purchase choice model on the {\scshape Bakery} dataset, we consider different parameter estimation methods and use a benchmark for comparison. Following \cite{Benson1}, we adopt the separable model for this purpose. In the separable model, each product has a certain utility, and a basket's utility is equal to the sum of the utilities of its elements assuming a fixed basket size. The key feature of this model is that demands for different products are independent of each other. The probability of choosing basket
$\{i_1, i_2, \dots, i_k\}$ from the given choice set~$\mathcal{N}$ under the separable model is given by:
\begin{equation*}
\hat{p}_{sep}(\{i_1, i_2, \dots, i_k\}) = \hat{p}_{size}(k) \dfrac{\prod_{i = i_1}^{i_k} \hat{p}_i}{\sum_{\{j_1, j_2, \dots, j_k\} \subseteq \mathcal{N}} {\prod_{j = j_1}^{j_k} \hat{p}_j}},
\end{equation*}
where $\hat{p}_{size}(k)$ denotes the empirical probability of a customer selecting a
basket of size $k$, and $\hat{p}_i$ is the empirical probability of a customer selecting item $i \in \mathcal{N}$. Note that while procedures based on exact MLE have been successfully used for parameter estimation in rather restrictive multi-purchase settings,  utilizing standard MLE is intractable (even for small problem instances) when studying such general multi-purchase settings as we do in our work.

Let $LL_{BM}$ be the mean log-likelihood value for the benchmark model on the test sample. Furthermore, let $LL_{DC}$, $LL_{SML}$, and $LL_{ML}$ represent the mean log-likelihood values for the Ising model with parameters estimated using (1) the DC approach, (2) approximate sparse maximum likelihood estimation with $\rho = 0.015$, and (3) approximate maximum likelihood estimation (where $\rho = 0$), respectively. The value $\rho = 0$ is selected by the cross-validation procedure described above, while the sparse estimate with $\rho = 0.015$ is included to quantify the loss in fit associated with the sparse network visualization. We calculate the relative improvements of the Ising model with different parameter estimates over the benchmark model as follows: $\Delta_{DC} = e^{LL_{DC} - LL_{BM}} = 1.93$, $\Delta_{SML} = e^{LL_{SML} - LL_{BM}} = 2.09$, and $\Delta_{ML} = e^{LL_{ML} - LL_{BM}} = 2.64$. As expected, creating sparse parameter estimates through $l_1$-regularization leads to less accurate estimates and consequently a reduction in relative improvement. The DC approach results in a slightly worse quality of fit compared to both types of approximate likelihood estimation, yet it is still comparably effective. Importantly, all parameter estimation approaches demonstrate a noticeably better performance of the Ising model as a multi-purchase choice model compared to the benchmark model.

The example studied above serves as a clear illustration of the potential and practical relevance of our modeling framework. However, evaluating the performance of the Ising model as a multi-purchase choice model on large datasets in terms of likelihood is challenging due to the computational intractability of the partition function. For smaller datasets, the performance of the MVL model, which is essentially equivalent to the binary Ising model, has been explored in multiple studies (see Section~\ref{2sec:lit}). Overall, the MVL is a well-established multi-purchase choice model that has proven useful in various practical situations. As for the parameter estimation problem for the Ising model (also known as the inverse Ising problem), it represents a significant research area with numerous publications (see, e.g., \citealp{Nguyen1}, \citealp{Braunstein1}, \citealp{Sano1}). One of the contributions of our paper is to highlight that this vast methodology developed for the Ising model can be effectively utilized for assortment optimization purposes. Although there is no material disadvantage to deploying estimation methods based on spin variables given the straightforward linear transformation in Proposition \ref{prop_transf_equiv}, developing extensions of the two methods described here to the binary case could be explored in future research.

\section{Assortment Optimization}
\label{2sec:ao_single}

\subsection{Problem Definition and Theoretical Insights}

Suppose that customers' choices follow the binary Ising model in~(\ref{ising_choice}). This binary model is the choice model used throughout the assortment optimization problem. Its parameters can be obtained by applying one of the two estimation methods described in Section~\ref{2sec:estim} -- Approximate Sparse MLE and DC -- to fit the spin Ising model to sales data observed under the full assortment $\mathcal{N}$, and then transforming the obtained parameter estimates using Proposition~\ref{prop_transf_equiv}. Let $r_j$ denote the gross profit (alternatively, the revenue) per unit of product $j\in \mathcal{N}$. Then, the expected profit generated by a random customer given assortment $S \subseteq \mathcal{N}$ is as follows:
\begin{equation*}
    R(S) = \sum\limits_{x \in \mathcal{X}(S)}\Bigl(p_{\theta}(x | S)\sum\limits_{j \in S} r_{j} x_j\Bigr).
\end{equation*}
The assortment optimization problem under the Ising model can be formulated as the problem of maximizing $R(S)$ over all possible assortments $S \subseteq \mathcal{N}$, i.e.,
\begin{equation}
\label{ising_assort_single}
    \max_{S\subseteq \mathcal{N}}  R(S).
\end{equation}
Note that this formulation does not depend on the expected number of customers since the optimal solution to an optimization problem does not change if the objective function is multiplied by a positive constant.

Solving optimization problem~(\ref{ising_assort_single}) is computationally difficult. In fact, the problem is NP-hard, as demonstrated by the following theorem:
\begin{Theorem}
\label{prop_assortm_np_hard}
The assortment optimization problem under the Ising model is NP-hard.
\end{Theorem}
\begin{Remark}
While we study an unconstrained setting in our work, if a capacity constraint is introduced in (\ref{ising_assort_single}) (e.g. as considered by \cite{Tulabandhula1,Jasin1}), then the capacitated assortment optimization problem remains NP-hard since the uncapacitated problem is a special case of the capacitated one (in the case of infinite capacity).
\end{Remark}
\noindent The proof (see \ref{prop_assortm_np_hard_proof}) proceeds by showing that the decision version of problem~(\ref{ising_assort_single}) is NP-hard. Consequently, computing the exact optimal value of problem~(\ref{ising_assort_single}) is NP-hard as well. Hence, unless \(P=NP\), there is no polynomial-time algorithm for solving problem~(\ref{ising_assort_single}) exactly.

Despite the computational complexity of the considered assortment optimization problem, we can derive several theoretical insights into the structure of its optimal solution by leveraging the graphical representation of the Ising model.
\begin{Proposition}
\label{ising_opt_proper}
If a node in the graphical representation of the Ising model is isolated, then the corresponding product belongs to the optimal assortment.
\end{Proposition}
\noindent The proof of this proposition is provided in \ref{ising_opt_proper_proof}. This is an intuitively clear result -- the fact that $\theta_{kj} = 0$ $\forall j \in \mathcal{N}\backslash\{k\}$ means that the probability of any product $j \in \mathcal{N}\backslash\{k\}$ belonging to a random basket is not affected by the presence of product~$k$ in the assortment. However, if a node in the graphical representation of the Ising model is such that there are no edges with negative weights incident on it, then the corresponding product does not necessarily belong to the optimal assortment. Consider the following example:
\begin{Example}
\label{ising_opt_proper2}
Let $\mathcal{N} = \{1,2,3\}$, 
 $r_1 = 10$, $r_2 = 10$, and $r_3 = 100$. Suppose that the customer choices follow the Ising model with the following parameters (see Figure~\ref{counterex_viz} for the graphical representation of this model):
\begin{equation*}
    \theta =
\begin{pmatrix}
    1 & 5 & 2 \\
    5 & 5 & -5  \\
    2 & -5  & 5 
\end{pmatrix}.
\end{equation*}
Then, it can be checked that $P(\{1,2,3\}) \approx 47$ and $P(\{2,3\}) \approx 55$, meaning that it is more profitable to offer products $\{2, 3\}$ than products $\{1, 2, 3\}$ despite the fact that all edges connected to node $1$ have positive weights.
\begin{figure}[H]
\centering
    \includegraphics[width=0.25\textwidth]{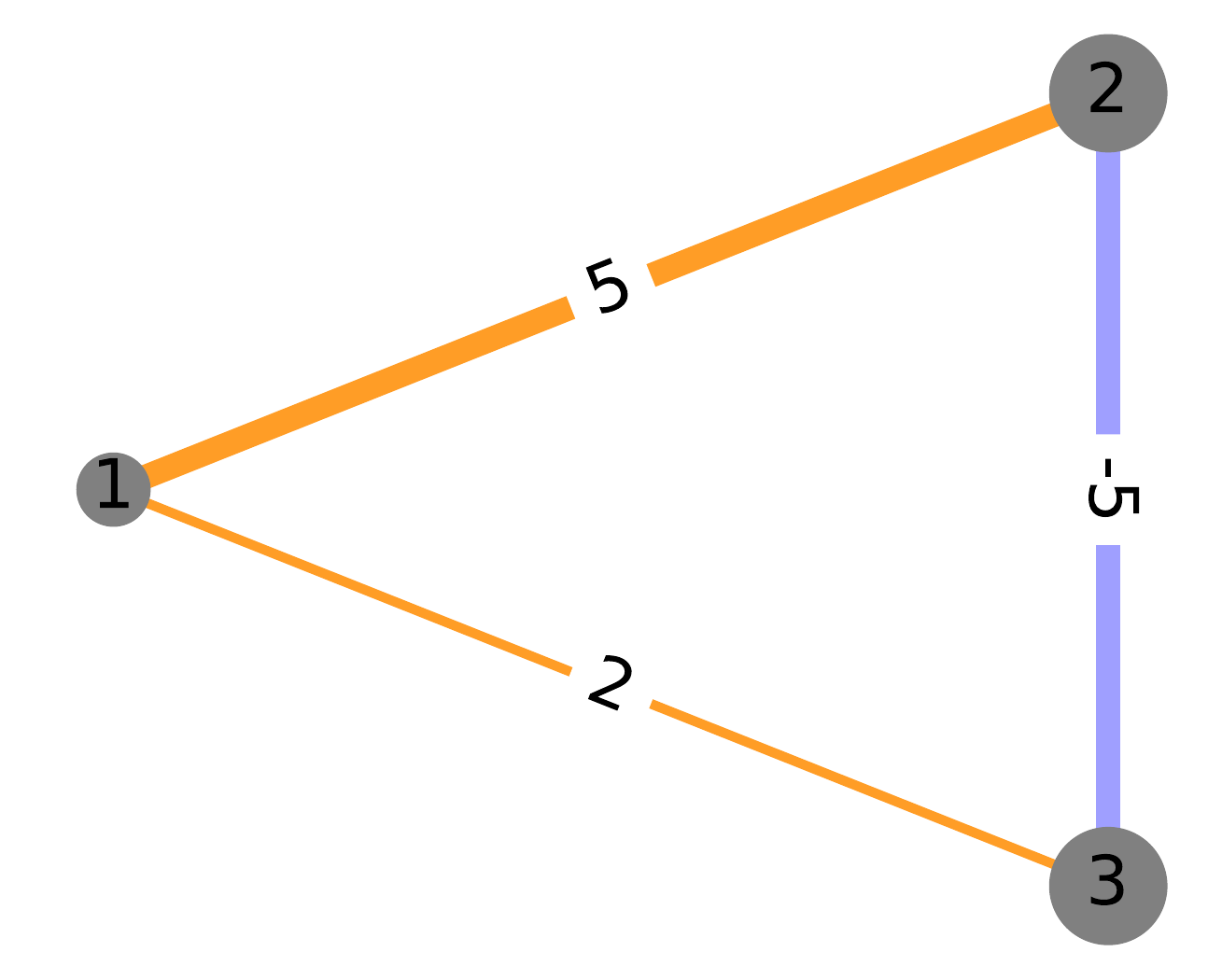}
    \caption{Example of when it is profitable to exclude a product from the assortment even though there are only positive edges connected to it.}
    \vspace{-0.1cm}
    \label{counterex_viz}
\end{figure}
\end{Example}

The above example demonstrates the complexity of the structure of the optimal assortment under the Ising model. A product added to a basket can have only positive \textit{direct} impacts on the probabilities of other products belonging to this basket. However, when all effects are taken into account jointly, it may happen that having such a product in the assortment reduces the marginal probabilities of customers buying some other products, meaning that this product does not necessarily belong to the optimal assortment after all.

Next, we can further exploit the fact that the Ising model is an MRF and formulate the following proposition:

\begin{Proposition}
\label{ising_opt_proper3}
Suppose that $\mathcal{N} = \mathcal{H} \sqcup \mathcal{K}$. If $\theta_{ij} = 0$ $\forall i \in \mathcal{H}$, $j \in \mathcal{K}$, then the assortment optimization problem~(\ref{ising_assort_single}) can be separated into assortment optimization problems for product sets $\mathcal{H}$ and $\mathcal{K}$. 
\end{Proposition}
\noindent This proposition, the proof of which is provided in \ref{ising_opt_proper3_proof}, means that the properties of optimal assortments can be formulated on the level of isolated subgraphs of the graphical representation of the Ising model rather than on the level of the whole graph.

In Theorem~\ref{ising_opt_proper4}, we provide a sufficient condition under which all products in an isolated subgraph belong to the optimal assortment. We are able to do so by showing that removing a product from such a subgraph can only reduce the marginal probabilities of customers choosing other products. When it comes to marginal probabilities, one can think of removing product~$k$ from the assortment as considering conditional marginal probabilities under condition $x_k = 0$. 
Formally,
\begin{equation*}
p_{\theta}(x_l = 1|x_k = 0, S) = p_{\theta}(x_l = 1|S \backslash \{k\}).
\end{equation*}
Indeed, one can easily see that
\begin{align*}
&p_{\theta}(x_l = 1|x_k = 0, S) = \dfrac{p_{\theta}(x_l = 1, x_k = 0| S)}{p_{\theta}(x_k = 0 | S)} \\
& \hspace{1cm} = \dfrac{\sum\limits_{x \in \mathcal{X}(S): \ x_l = 1, x_k = 0}\exp\Bigl(\sum\limits_{i \in S}\theta_{ii}x_i + \sum\limits_{i \ne j, \ i, j \in S}x_i\theta_{ij}x_j\Bigr)}{\sum\limits_{x \in \mathcal{X}(S): \ x_k = 0}\exp\Bigl(\sum\limits_{i \in S}\theta_{ii}x_i + \sum\limits_{i \ne j, \ i, j \in S}x_i\theta_{ij}x_j \Bigr)} \\ 
& \hspace{1cm}=  \dfrac{\sum\limits_{x \in \mathcal{X}(S\backslash\{k\}): \ x_l = 1}\exp\Bigl(\sum\limits_{i \in S\backslash\{k\}}\theta_{ii}x_i + \sum\limits_{i \ne j, \ i, j \in S\backslash\{k\}}x_i\theta_{ij}x_j\Bigr)}{\sum\limits_{x \in \mathcal{X}(S\backslash\{k\})}\exp\Bigl(\sum\limits_{i \in S\backslash\{k\}}\theta_{ii}x_i + \sum\limits_{i \ne j, \ i, j \in S\backslash\{k\}}x_i\theta_{ij}x_j \Bigr)} = p_{\theta}(x_l = 1|S \backslash\{k\}). 
\end{align*} 
This observation plays an important role in the proof of Theorem~\ref{ising_opt_proper4}.

\begin{Theorem}
\label{ising_opt_proper4}
If an isolated subgraph in the graphical representation of the Ising model does not contain edges with negative weights, then all products from this subgraph belong to the optimal assortment.
\end{Theorem}

\noindent
Although the theorem statement is intuitively clear, its formal proof is rather involved and requires a combination of different techniques and supporting statements. The full proof can be found in \ref{ising_opt_proper4_proof}. Besides explicit theoretical insights, Theorem~\ref{ising_opt_proper4} also provides intuition on the sizes of optimal assortments in certain cases. In particular, if there are few pairwise negative dependencies, or if the absolute values of negative parameters $\theta_{ij}$ are small compared to the rest of the Ising model parameters, then one can expect that it is optimal to offer almost all products.

We can also identify a specific condition for the parameters associated with a product under which this product has to belong to the optimal assortment. Consider the following proposition (see \ref{ising_opt_proper5_proof} for the proof): 
\begin{Proposition}
\label{ising_opt_proper5}
Let $\mathcal{H} \subseteq \mathcal{N}$ induce an isolated subgraph in the graphical representation of the Ising model. If product $k \in \mathcal{H}$ is such that $\theta_{ki}= \alpha r_i$  $ \forall i \in \mathcal{H} \backslash \{k\}$, where $\alpha$ is an arbitrary positive constant, then product $k$ belongs to the optimal assortment.
\end{Proposition}
\noindent The condition in Proposition~\ref{ising_opt_proper5} has an intuitive interpretation. It describes a product whose positive demand dependencies with the other products in the same isolated subgraph are aligned with their margins. Such an alignment may arise when the complements that are most strongly tied to this product are also the ones with higher margins. This can occur, for example, when the strongest complements are product-specific or differentiated, while weaker complements are more generic. For instance, a printer may be more strongly tied to proprietary ink cartridges than to standard paper. The exact proportionality condition is stylized, but it provides a clean sufficient condition under which a product must be retained in the optimal assortment. This is illustrated by the following example:
\begin{Example}
Let $\mathcal{N} = \{1,2,3\}$, $r_1 = 10$, $r_2 = 4$, and $r_3 = 16$, and suppose that the customer choices follow the Ising model with the following parameters:
\begin{equation*}
    \theta =
\begin{pmatrix}
    1 & 2 & 8 \\
    2 & 5 & -5  \\
    8 & -5  & 5 
\end{pmatrix}.
\end{equation*}
Then, from Proposition \ref{ising_opt_proper5} it follows that product 1 belongs to the optimal assortment. Indeed, by setting $\alpha=0.5$ one has $\theta_{12}=\alpha r_2$ and $\theta_{13}=\alpha r_3$, which satisfies the required condition.
\end{Example}

The theoretical results obtained can be used towards assortment optimization, in particular, to reduce the dimensionality of problem~(\ref{ising_assort_single}). After producing the graphical representation of the Ising model as discussed in Section~\ref{2sec:estim}, the following preprocessing procedure can be applied:
\begin{enumerate}
    \item[1.] Add isolated nodes to the optimal assortment, thereby removing them from consideration;
    \item[2.] Separate the graph into connected components (maximal connected subgraphs);
    \item[3.] If some connected components do not contain edges with negative weights, add all nodes from these components to the optimal assortment;
    \item[4.] For each remaining connected component, formulate the assortment optimization problem considering nodes from this component as a separate product portfolio.
\end{enumerate}
Given that problem~(\ref{ising_assort_single}) grows exponentially with the size of the product portfolio $n$, the preprocessing procedure presented above can lead to a significant reduction in the size of the assortment optimization problem when the graphical representation of the corresponding Ising model has a sparse structure. It is important to note that removing edges from the graph has a twofold effect on the graph structure: it may lead to the creation of both isolated nodes and isolated subgraphs. The significance of the former effect will be discussed in the next subsection.

\subsection{Impact of Edge Removal on Isolated Nodes}

Although the positive impact of edge removal on optimization performance is not necessarily directly attributable to the creation of either isolated edges or isolated subgraphs individually, in this subsection we investigate the effect of removing edges on the number of isolated nodes. The number of isolated nodes created can be seen as an indirect metric for measuring problem size reduction. To evaluate the practical relevance of creating isolated nodes, we examine the graph structures of the Ising models that correspond to six real-world datasets from \cite{Benson1} that contain basket choices: {\scshape Bakery} (comprising 50 options), {\scshape WalmartDepts} (66), {\scshape WalmartItems} (183), {\scshape Kosarak}~(2605), {\scshape LastfmGenres} (413), and {\scshape Instacart} (9544). Please refer to \cite{Benson1} for a description of these datasets.

\begin{figure}[t]
\centering
\begin{subfigure}{.32\textwidth}
  \centering
  \includegraphics[width=\linewidth]{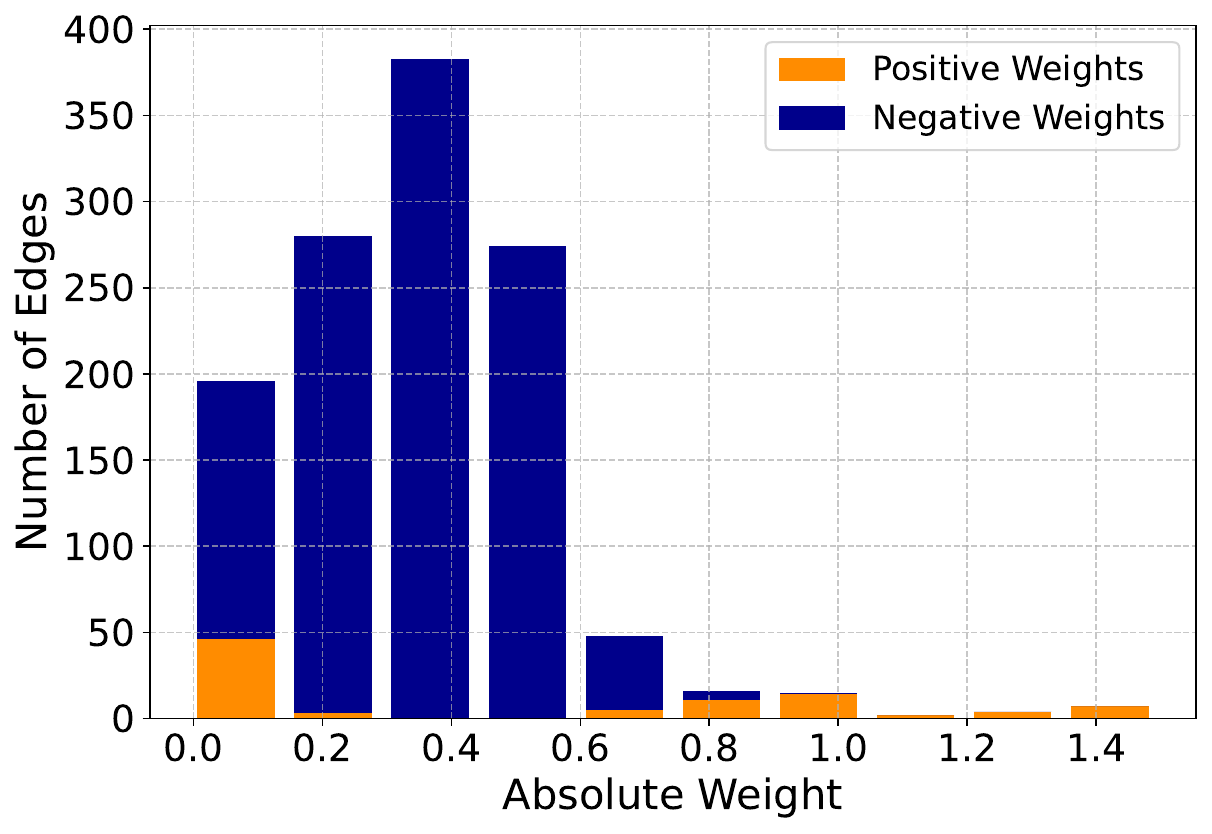}
  \caption{\footnotesize{{\scshape Bakery.}}}
\end{subfigure}
\begin{subfigure}{.32\textwidth}
  \centering
  \includegraphics[width=\linewidth]{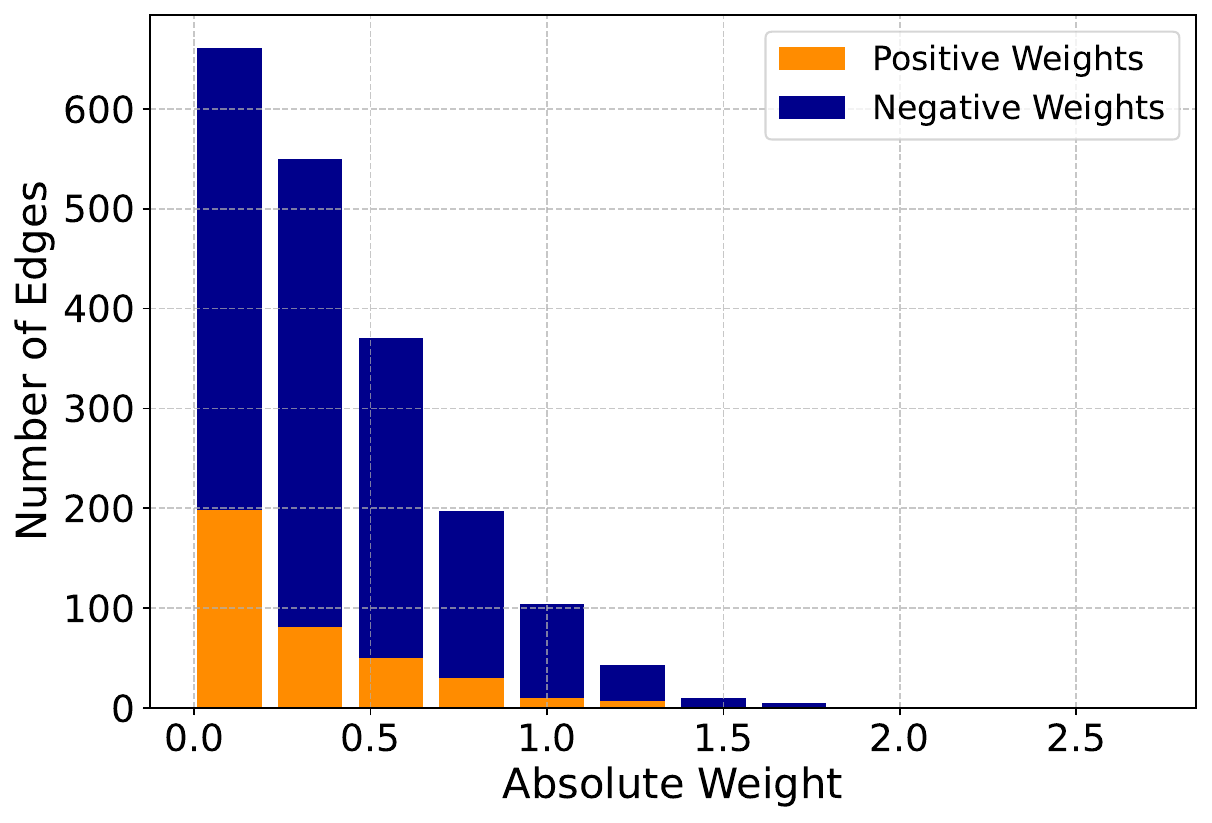}
  \caption{\footnotesize{{\scshape WalmartDepts.}}}
\end{subfigure}
\begin{subfigure}{.325\textwidth}
  \centering
  \includegraphics[width=\linewidth]{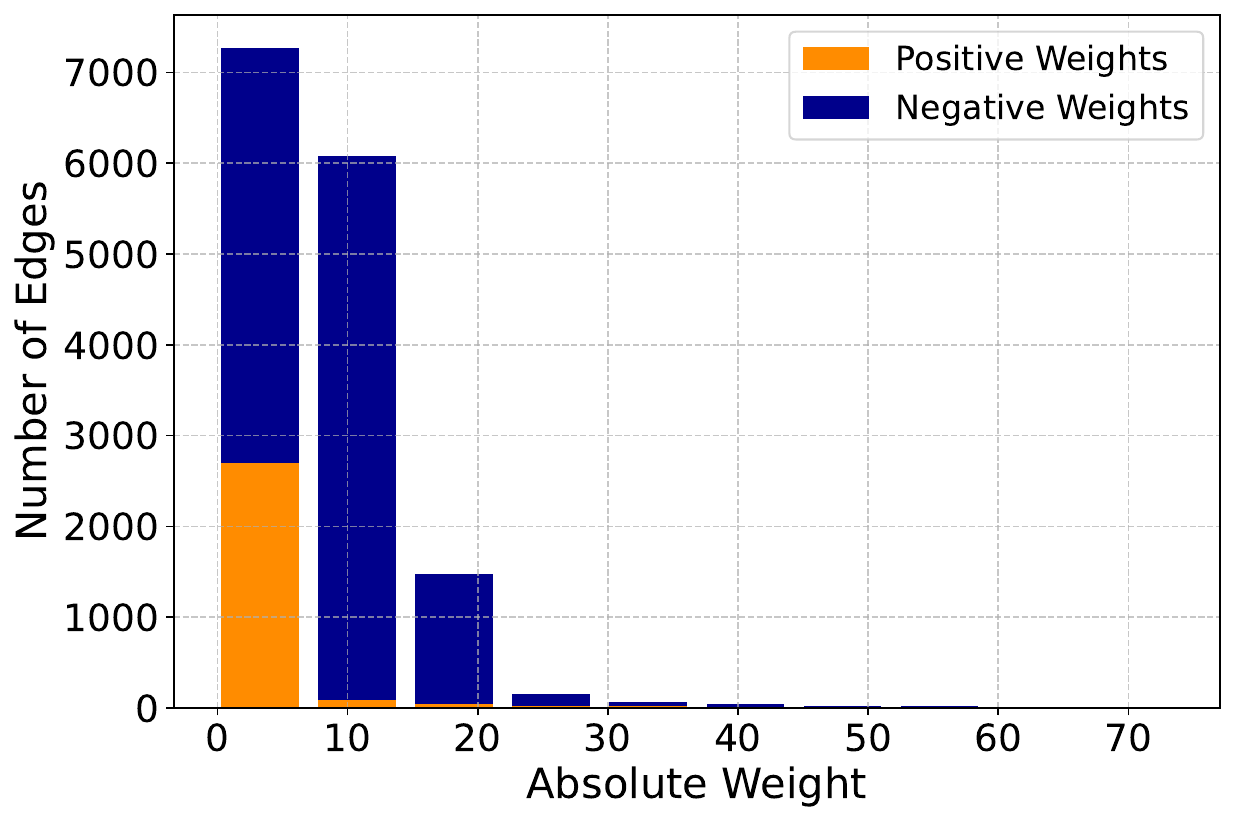}
  \caption{\footnotesize{{\scshape WalmartItems.}}}
\end{subfigure}
\begin{subfigure}{.31\textwidth}
  \centering
  \includegraphics[width=\linewidth]{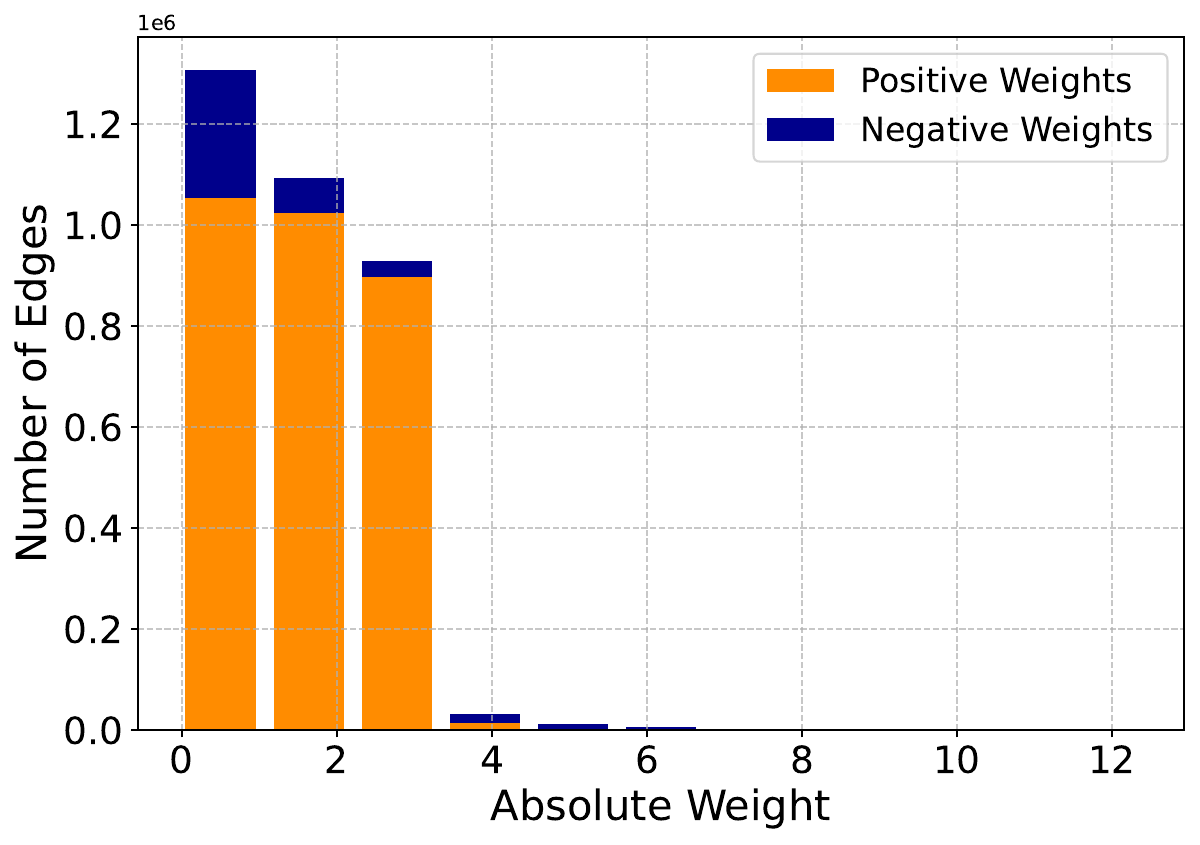}
  \caption{\footnotesize{{\scshape Kosarak.}}}
\end{subfigure}
\begin{subfigure}{.33\textwidth}
  \centering
  \includegraphics[width=\linewidth]{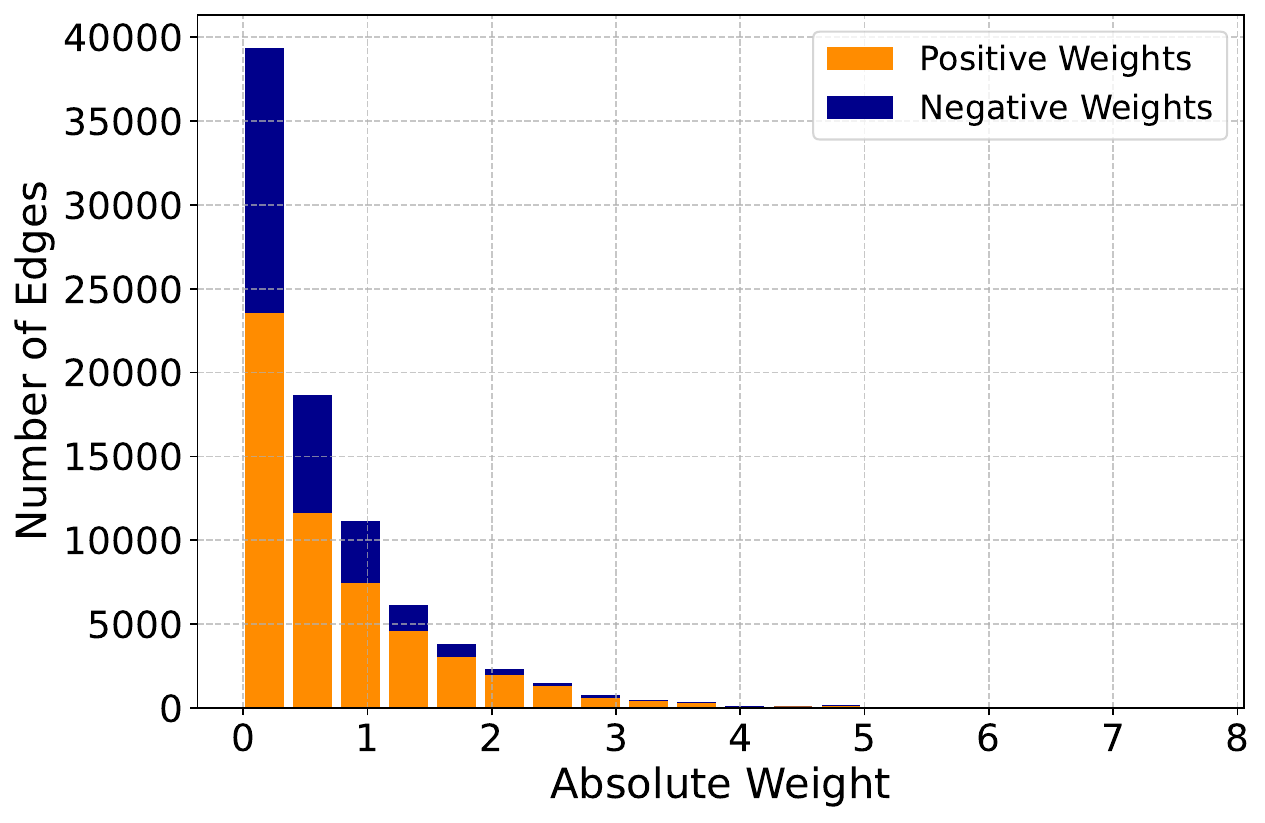}
  \caption{\footnotesize{{\scshape LastfmGenres.}}}
\end{subfigure}
\begin{subfigure}{.315\textwidth}
  \centering
  \includegraphics[width=\linewidth]{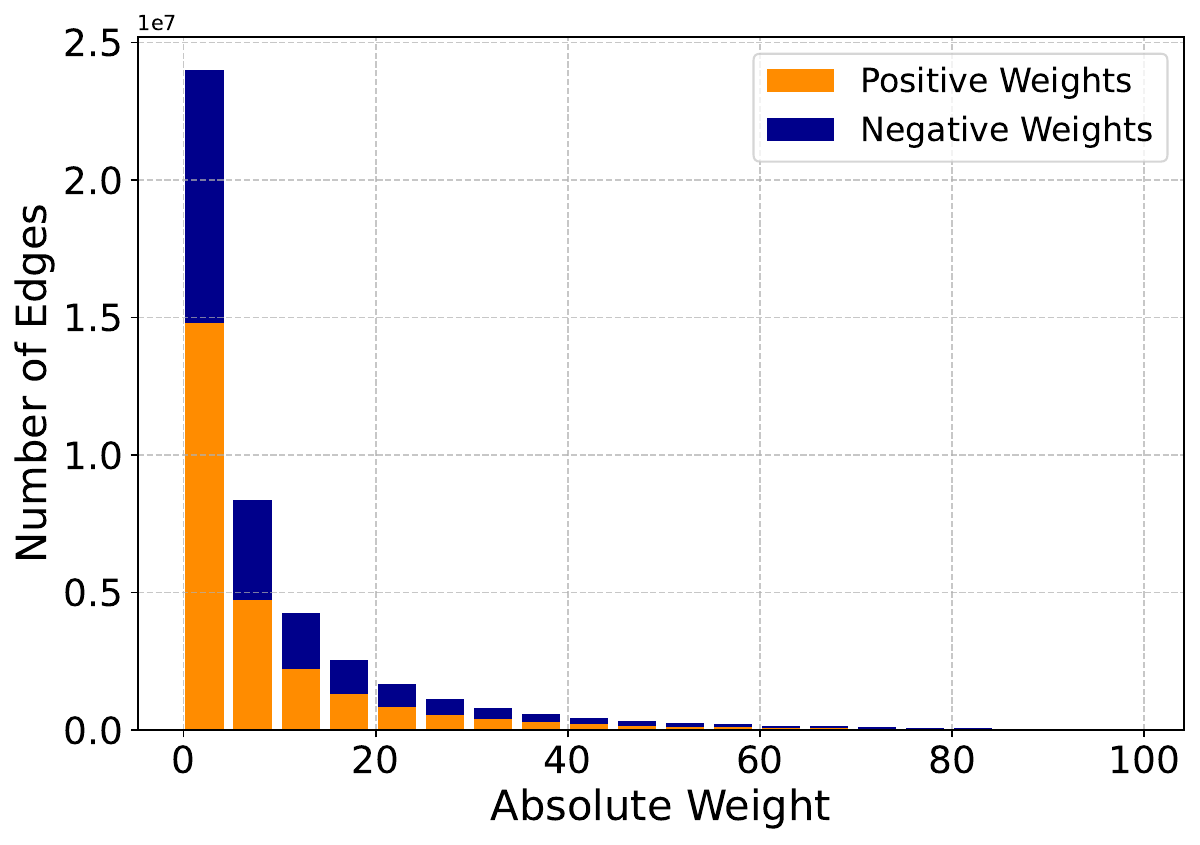}
  \caption{\footnotesize{{\scshape Instacart.}}}
\end{subfigure}
\caption{Distribution of edge weights (Ising model couplings) for real-world datasets.}
\vspace{-0.1cm}
\label{edges_distributions}
\end{figure}

We estimate the parameters of the corresponding Ising models using the DC approach (see Section~\ref{2sec:estim}). The choice of the DC approach is driven by the computational tractability of the estimation process. Although the approximate sparse MLE method may be more effective in reducing problem dimensionality by directly generating sparse graphs, its application to the larger datasets -- {\scshape Kosarak} and {\scshape Instacart} -- is computationally demanding and much less empirically tractable compared with the DC approach, which represents a good trade off between model simplicity (i.e. sparsity) and computational tractability of estimation. The distributions of edge weights (Ising model couplings) are displayed in Figure~\ref{edges_distributions}. 
Overall, the graph structures vary significantly across the datasets. For example, {\scshape Bakery, WalmartDepts}, and {\scshape WalmartItems} predominantly have edges with negative weights, whereas {\scshape Kosarak, LastfmGenres}, and {\scshape Instacart} predominantly feature edges with positive weights. The ranges of the absolute values of edge weights also vary considerably across the datasets. Nevertheless, the majority of edges in all six histograms typically fall into the bins that correspond to the lowest weights. The only exception is the {\scshape Bakery} dataset, where the third bin contains the most entries.

\begin{figure}[t]
\centering
\begin{subfigure}{.32\textwidth}
  \centering
  \includegraphics[width=\linewidth]{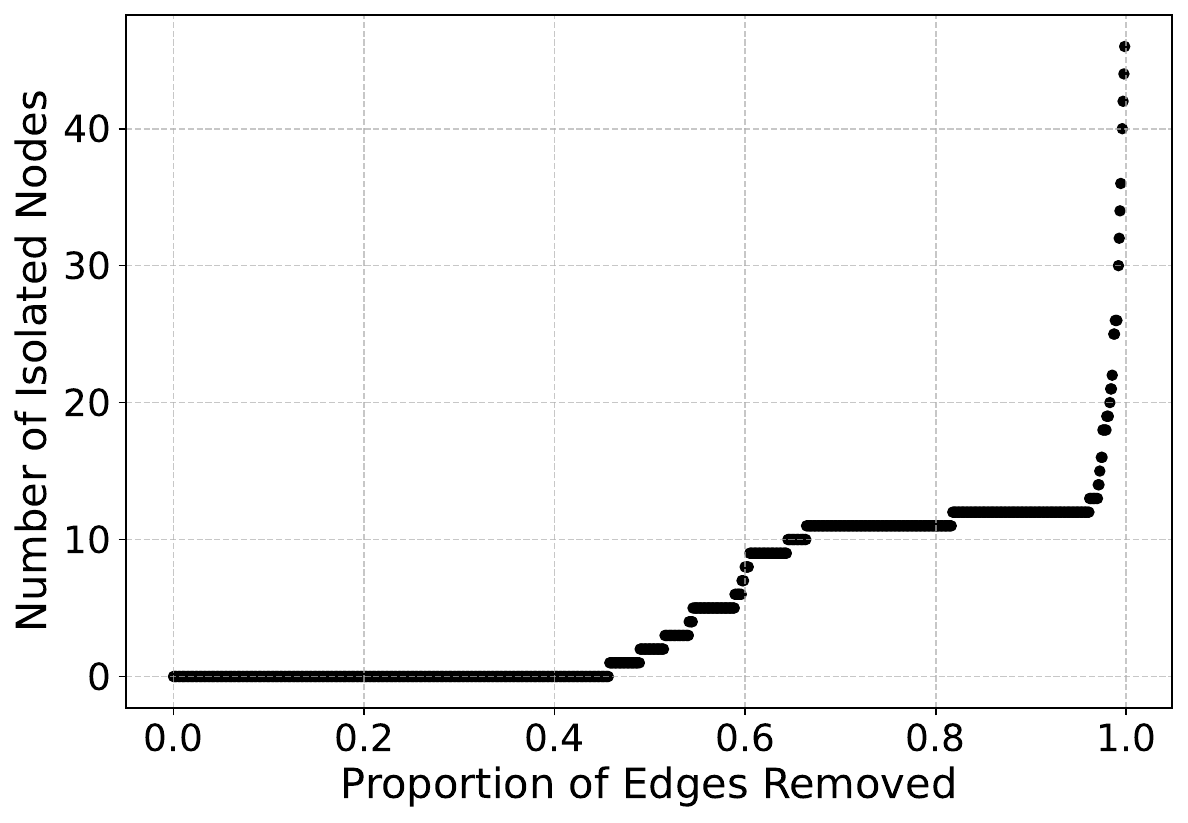}
  \caption{\footnotesize{{\scshape Bakery.}}}
\end{subfigure}
\begin{subfigure}{.322\textwidth}
  \centering
  \includegraphics[width=\linewidth]{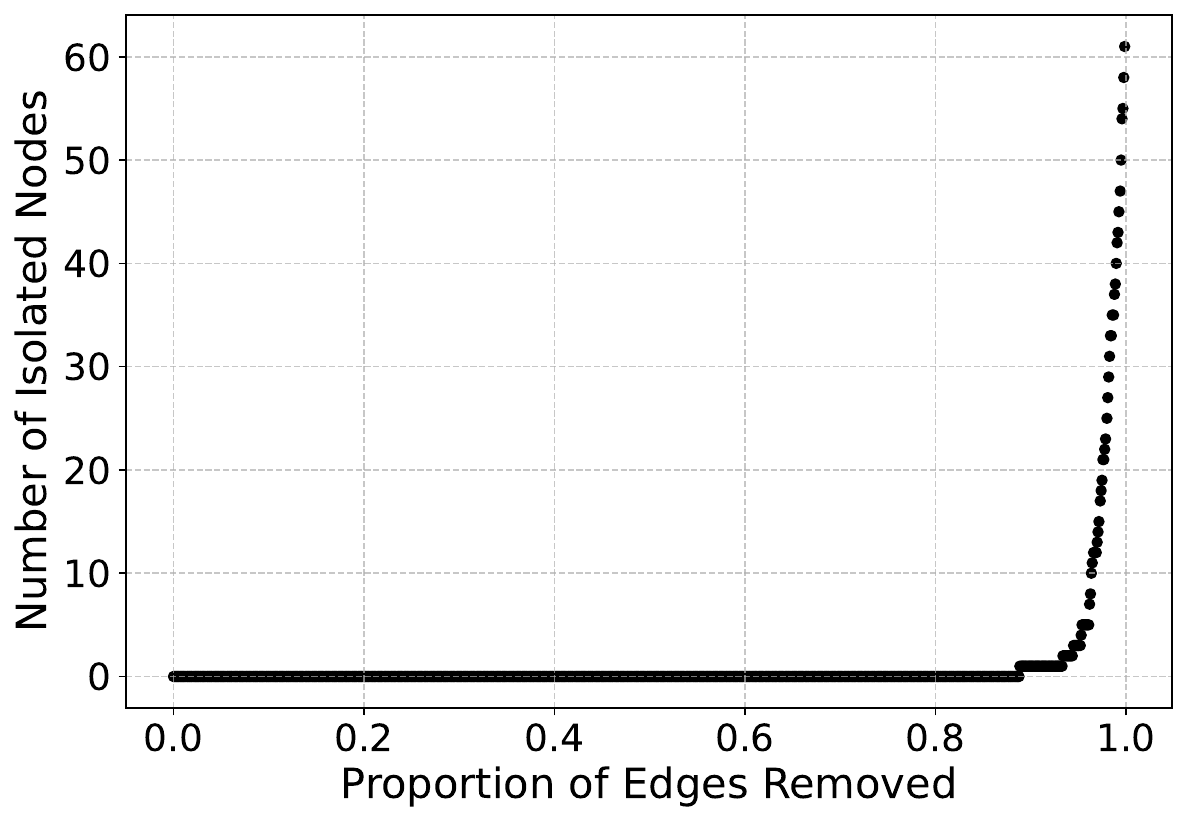}
  \caption{\footnotesize{{\scshape WalmartDepts.}}}
\end{subfigure}
\begin{subfigure}{.33\textwidth}
  \centering
  \includegraphics[width=\linewidth]{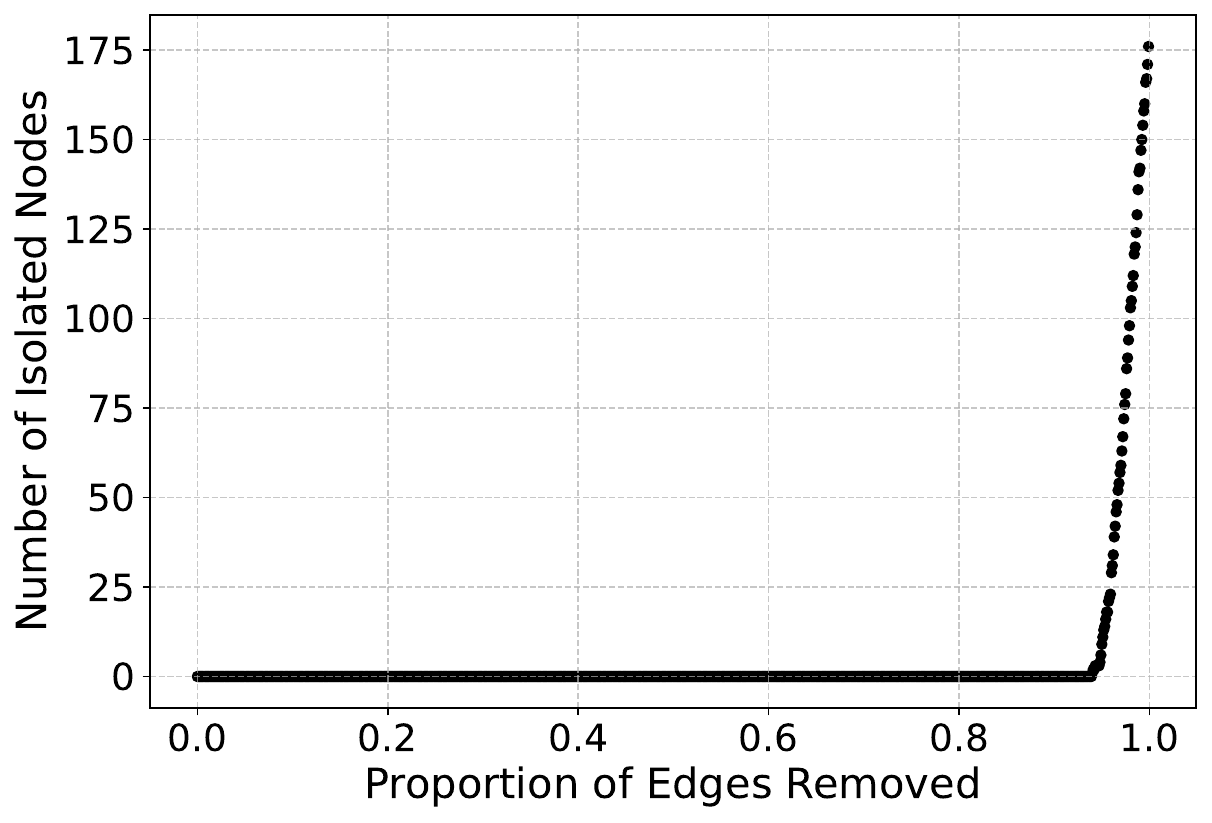}
  \caption{\footnotesize{{\scshape WalmartItems.}}}
\end{subfigure}
\begin{subfigure}{.32\textwidth}
  \centering
  \includegraphics[width=\linewidth]{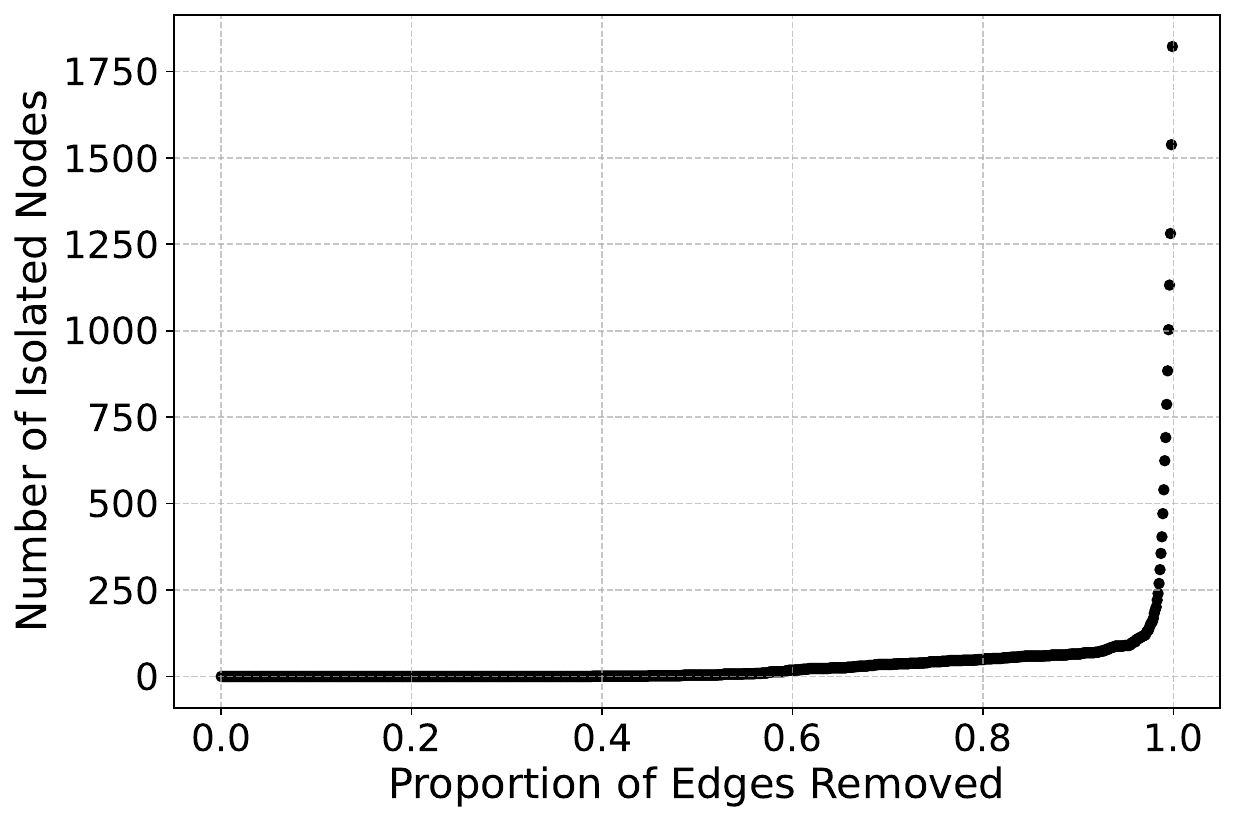}
  \caption{\footnotesize{{\scshape Kosarak.}}}
\end{subfigure}
\begin{subfigure}{.315\textwidth}
  \centering
  \includegraphics[width=\linewidth]{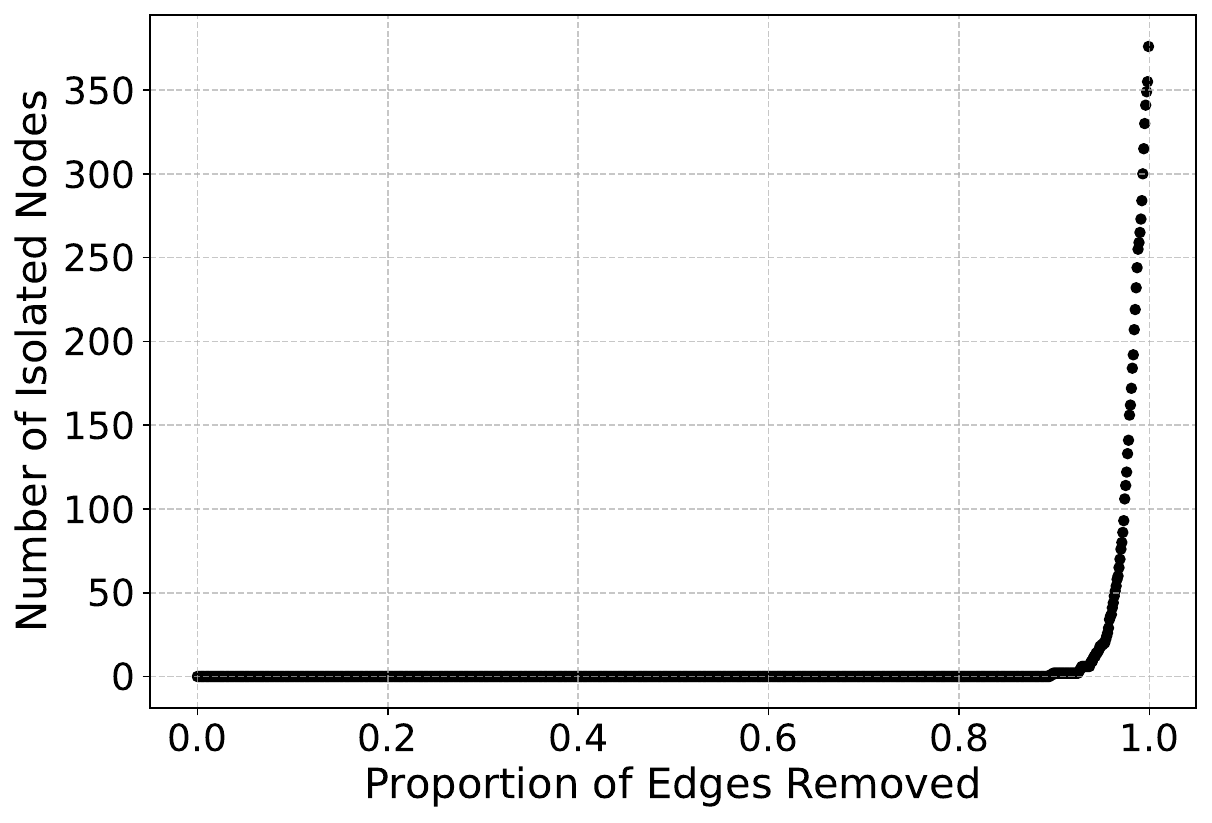}
  \caption{\footnotesize{{\scshape LastfmGenres.}}}
\end{subfigure}
\begin{subfigure}{.32\textwidth}
  \centering
  \includegraphics[width=\linewidth]{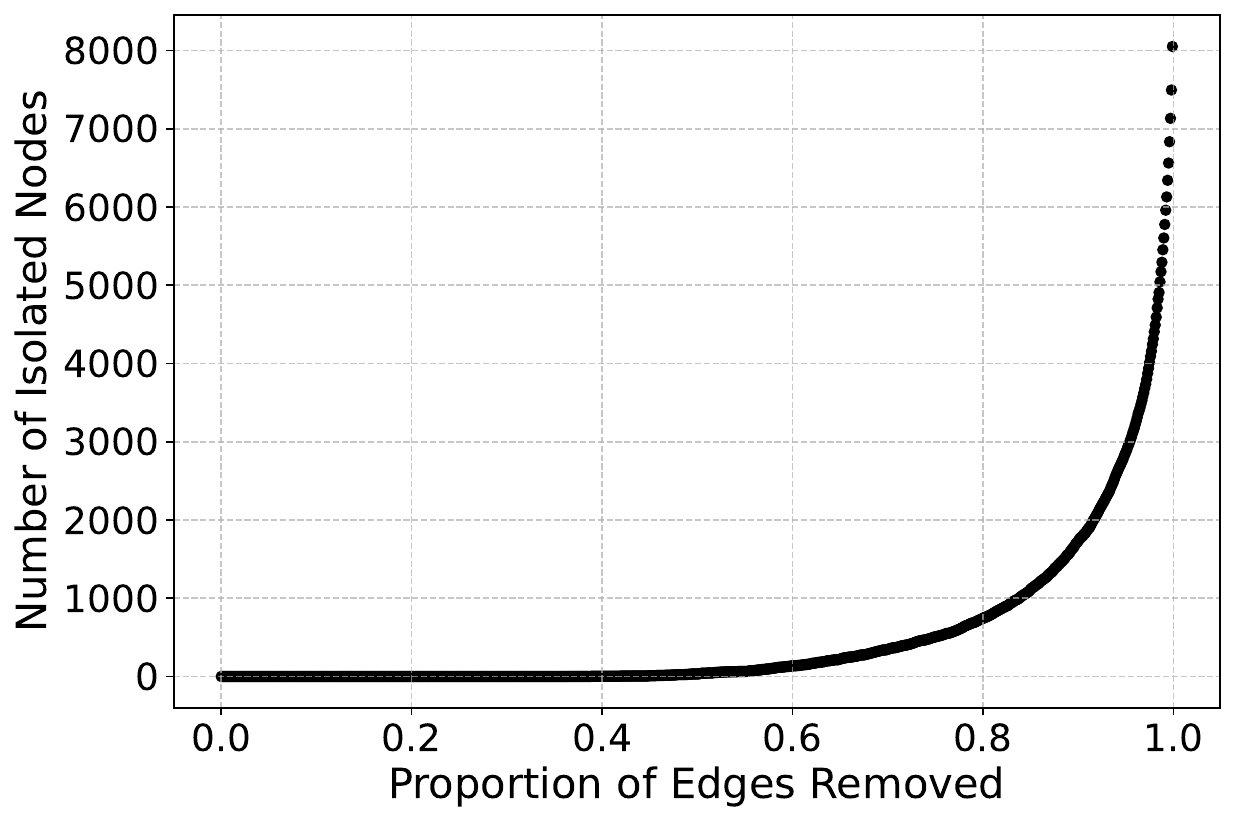}
  \caption{\footnotesize{{\scshape Instacart.}}}
\end{subfigure}
\caption{Number of isolated nodes as a function of the proportion of edges removed.}
\label{isolated_nodes_dependencies}
\end{figure}

Figure~\ref{isolated_nodes_dependencies} illustrates the relationships between the numbers of isolated nodes and the percentages of edges with the lowest absolute weights that have been removed across all six datasets. 
It shows that approximately 45\% of edges should be removed for isolated nodes to appear in {\scshape Bakery}, 90\% in {\scshape WalmartDepts}, 95\% in {\scshape WalmartItems}, 60\% in {\scshape Kosarak}, 90\% in {\scshape LastfmGenres}, and 55\% in {\scshape Instacart}. Taking the distributions of edge weights into account (Figure~\ref{edges_distributions}), this indicates that in some cases, removing only the edges that belong to the lowest-weight bins is sufficient for isolated nodes to emerge, whereas in others, removing edges with more substantial weights is required. For example, removing edges from the first bin would be sufficient for isolated nodes to appear in the case of {\scshape Instacart}. In contrast, for {\scshape LastfmGenres}, removing edges from the first five bins would be necessary, which could negatively affect the Ising model's ability to accurately represent the basket shopping behavior of customers. 

To conclude, the actual improvement in problem scale reduction due to edge removal is highly problem-specific, and in practice, an in-depth investigation of the structure of the graphical representation of the Ising model would be required for each individual dataset. While this subsection demonstrates that the creation of isolated nodes can play an important role in considerably reducing the size of the assortment optimization problem for some of the real-world datasets studied (e.g., {\scshape Bakery} or {\scshape Instacart}), it would be promising to investigate whether the creation of smaller, isolated subgraphs can reduce the problem size more effectively and more consistently. In fact, since isolated components tend to appear earlier as edges are gradually removed, it can reasonably be expected that both the number of isolated subgraphs and their maximum size are practically more important metrics for measuring problem size reduction, making the effect on isolated subgraphs an interesting avenue for future research.

\subsection{Heuristic Algorithms}

Next, we need a computationally tractable method that will allow us to solve the assortment optimization problem after completing the preprocessing stage. As mentioned earlier, problem~(\ref{ising_assort_single}) is NP-hard. Thus, we developed a metaheuristic algorithm that can be used to obtain high-quality solutions to the assortment optimization problem. Before providing a detailed description, let us first discuss the integral part of this algorithm, namely, how to evaluate potential solutions to the considered problem. 

Finding the exact value of the expected profit $R(S)$ might not be computationally feasible even for moderately sized assortments $S$. This is because of the need to compute the normalizing coefficient of $p_{\theta}(x|S)$, which contains an exponential number of terms. We thus need an efficient way of approximating the expected profit. This can be done by drawing a large number of samples from the Ising model, i.e., by generating a large number of purchased baskets. 
We employ Gibbs sampling for this purpose, which is a special variant of the Metropolis-Hastings algorithm. There are two well-known variations of Gibbs sampling: systematic scan Gibbs sampling and random scan Gibbs sampling (see, e.g., \citealp{He1}). In the former method, each new sample $x^{(l)}$ is obtained by initializing $x^{(l)} = x^{(l-1)}$ and iteratively sampling vector components $x^{(l)}_k$ from $p_{\theta}(\cdot | x^{(l)}_{-k})$ for $k = 1, \dots, n$, whereas in the latter method, each new sample is obtained by randomly selecting index $k$ and sampling $x^{(l)}_k$ from $p_{\theta}(\cdot | x^{(l)}_{-k})$. Algorithms~\ref{gibbs_gen1} and~\ref{gibbs_gen2} describe these procedures adapted to sampling from the Ising model.

\begin{minipage}[ht]{0.47\textwidth}
\begin{algorithm}[H]
    \centering
    \captionsetup{font=small}
    \caption{Systematic scan Gibbs generator}\label{gibbs_gen1}
    \scriptsize
\algloop{Input}
\algblock{Input}{EndIf}
\algloop{Output}
\algblock{Output}{EndIf}
\algloop{Initialization}
\algblock{Initialization}{EndIf}
\begin{algorithmic}[1]
\Input
    \State $x^{(l-1)}$: previous sample
\EndIf
\Output
    \State $x^{(l)}$: new sample
\EndIf
\Initialization
    \State $n$ $\gets$ $\texttt{length}(x^{(l-1)})$
    \State $x^{(l)}$ $\gets$ $x^{(l-1)}$
\EndIf
 \For{$k$ \textbf{in} $[1, \dots, n]$}:
    \State Draw a random number $r \sim U(0,1)$
    \If{$r \leq g(1, x_{-k}^{(l)})/\bigl(g(0, x_{-k}^{(l)}) + g(1, x_{-k}^{(l)})\bigr)$ where \\
    \hskip\algorithmicindent \text{\:} $g(x_k, x_{-k})=\exp\bigl(\sum\limits_{i}\theta_{ii}x_i + \sum\limits_{i \ne j}x_i\theta_{ij}x_j\bigr)$}:
    \State $a \gets 1$ 
    \Else:
    \State $a \gets 0$
    \EndIf
\smallskip       
    \State $x^{(l)}_k$ $\gets$ $a$
\EndFor
\State \Return{$x^{(l)}$}
    \end{algorithmic}
\end{algorithm}
\end{minipage}
\hfill
\begin{minipage}[ht]{0.47\textwidth}
\begin{algorithm}[H]
    \centering
        \captionsetup{font=small}
    \caption{Random scan Gibbs generator}\label{gibbs_gen2}
    \scriptsize
\algloop{Input}
\algblock{Input}{EndIf}
\algloop{Output}
\algblock{Output}{EndIf}
\algloop{Initialization}
\algblock{Initialization}{EndIf}
\begin{algorithmic}[1]
\Input
    \State $x^{(l-1)}$: previous sample
\EndIf
\Output
    \State $x^{(l)}$: new sample
\EndIf
\Initialization
    \State $n$ $\gets$ $\texttt{length}(x^{(l-1)})$
    \State $x^{(l)}$ $\gets$ $x^{(l-1)}$
\EndIf
 \State Draw a random index $k$ from $[1, \dots, n]$
    \State Draw a random number $r \sim U(0,1)$
    \If{$r \leq g(1, x_{-k}^{(l)})/\bigl(g(0, x_{-k}^{(l)}) + g(1, x_{-k}^{(l)})\bigr)$ where \\
    \text{\:} $g(x_k, x_{-k})=\exp\bigl(\sum\limits_{i}\theta_{ii}x_i + \sum\limits_{i \ne j}x_i\theta_{ij}x_j\bigr)$}:
    \State $a \gets 1$ 
    \Else:
    \State $a \gets 0$
    \EndIf
\smallskip    
    \State $x^{(l)}_k$ $\gets$ $a$
\State \Return{$x^{(l)}$}
    \end{algorithmic}
\end{algorithm}
\end{minipage}
\begin{figure}[H]
\centering
\begin{subfigure}{.49\textwidth}
  \centering
  \includegraphics[width=\linewidth]{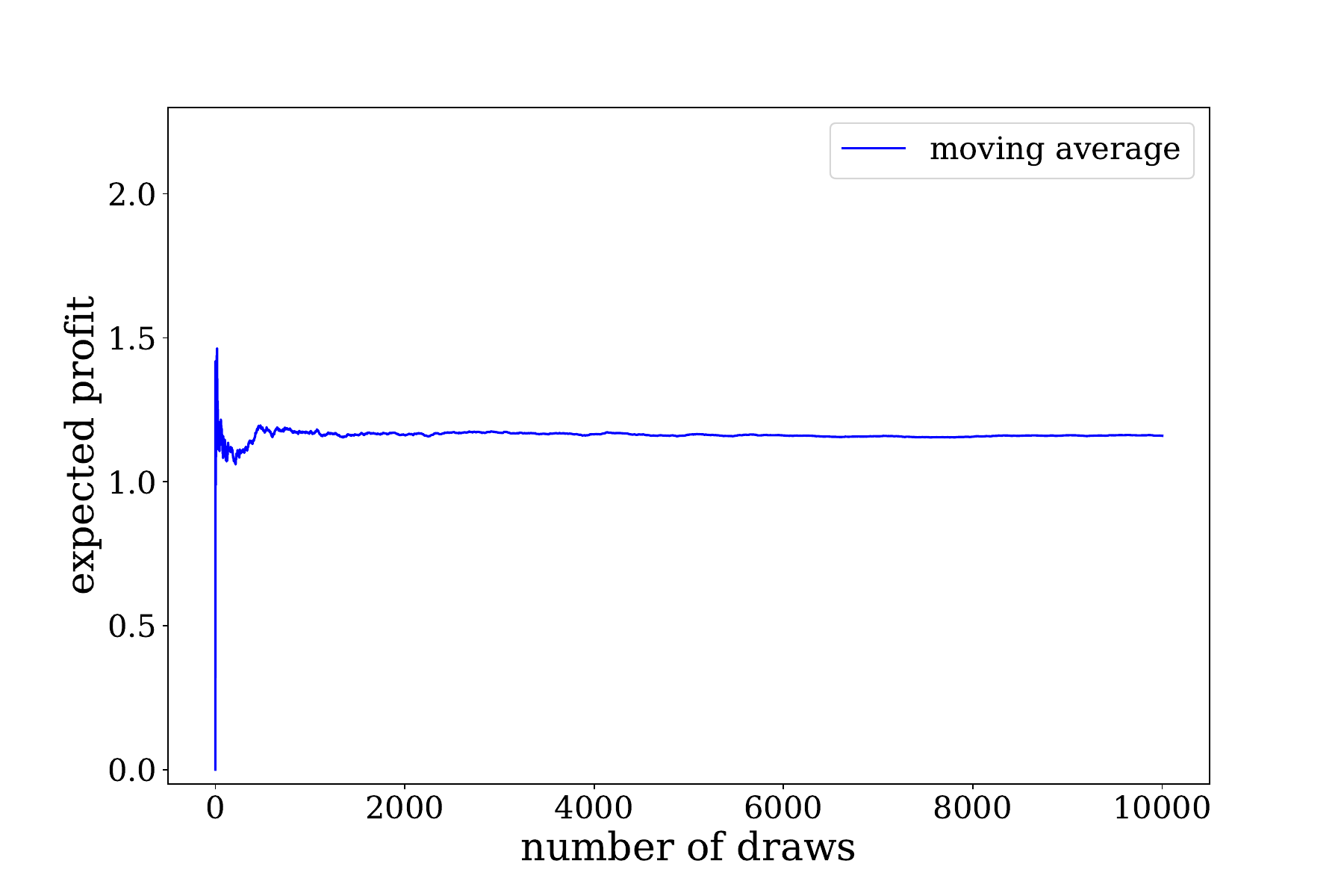}
  \caption{\footnotesize{Systematic scan Gibbs sampling.}}
\end{subfigure}
\begin{subfigure}{.49\textwidth}
  \centering
  \includegraphics[width=\linewidth]{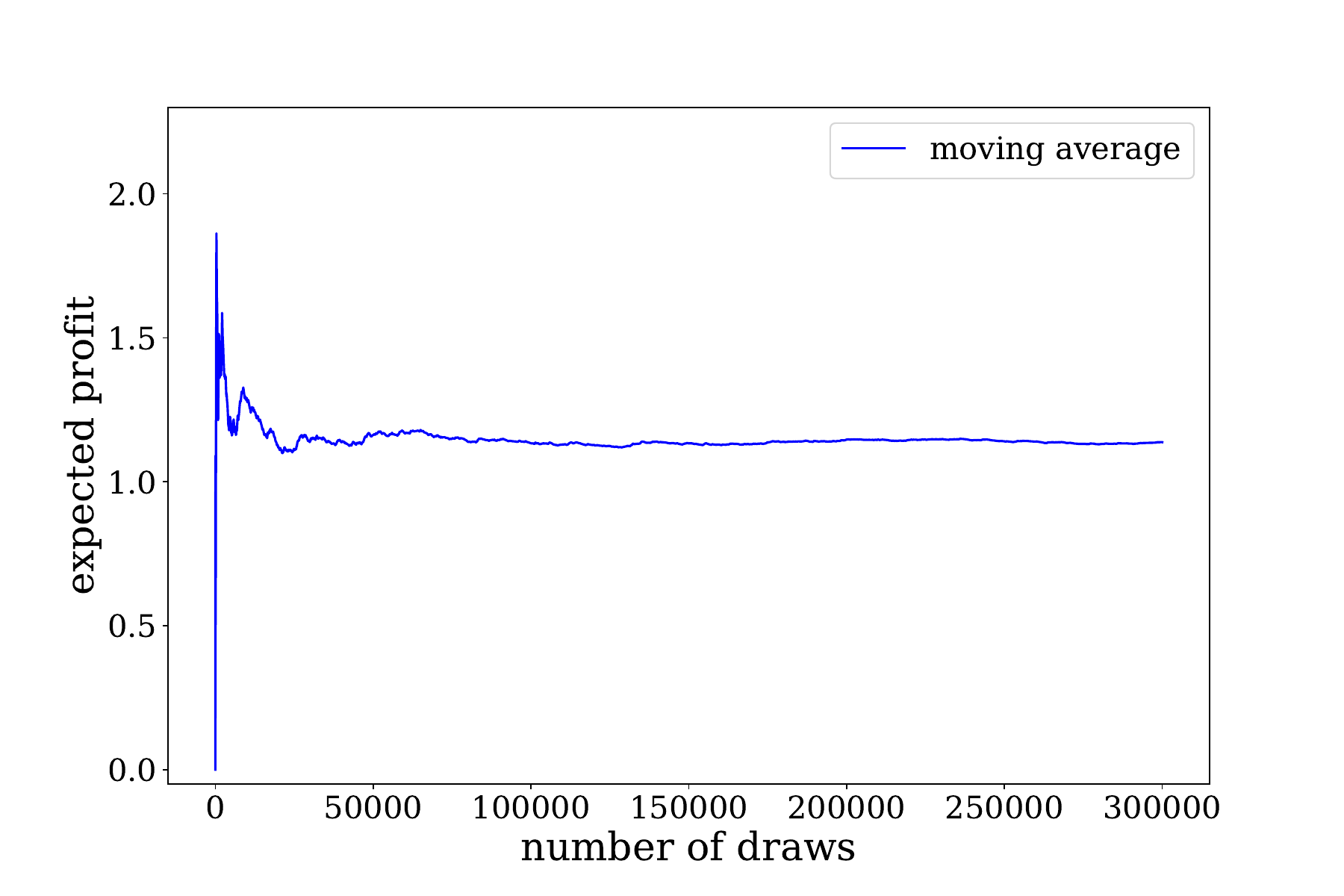}
  \caption{\footnotesize{Random scan Gibbs sampling.}}
\end{subfigure}
\caption{Example of expected profit estimation using Gibbs sampling for Ising model.}
\label{gibbs_sampl_perf}
\end{figure}

There is no universal answer as to whether systematic scan Gibbs sampling or random scan Gibbs sampling is more efficient in terms of convergence to the target distribution -- the choice is problem-specific, and either method can significantly outperform the other \citep{Gareth1}.  In practical applications, systematic scan Gibbs sampling is commonly used; however, most theoretical insights into convergence rates focus on the random update scheme \citep{Guo2}. To select the sampling procedure, we conducted a preliminary numerical analysis using the approximate maximum likelihood estimates of the Ising model parameters obtained in Section~\ref{2sec:estim} together with randomly generated profit margins $r_j \sim U(0,1)$. Based on this preliminary analysis, we decided to use systematic scan Gibbs sampling as it appeared to converge faster than random scan Gibbs sampling for a similar amount of time (see Figure~\ref{gibbs_sampl_perf} for an illustration).

Being able to approximate the expected profit for any assortment allows us to implement a metaheuristic algorithm for finding promising solutions to problem~(\ref{ising_assort_single}). 
One of the most prominent metaheuristic algorithms is Simulated Annealing. There are numerous variations of this procedure. In this work, we adopt the Simulated Annealing design described by \cite{Bierlaire1}. We set the number of algorithm iterations per temperature to one.  At each iteration, a candidate solution is generated and evaluated using simulations. Each candidate solution is selected randomly from the neighborhood of the current solution. We chose a fairly basic neighborhood structure: The neighborhood of an assortment comprises all possible assortments that can be obtained from the assortment in question by either removing or adding a product. If the candidate solution is better than the current solution, then the candidate solution is accepted and becomes the current solution. Otherwise, the candidate solution is accepted with a certain probability depending on the annealing temperature at this iteration. The temperature decreases with each new iteration. The maximum and minimum annealing temperatures are calculated based on the target maximum and minimum acceptance probabilities, respectively. It is done using parameter $d_{obj}$ representing the typical increase of the objective function in the given neighborhood structure, which can be estimated empirically. In this algorithm variation, the acceptance probability decreases linearly with respect to the iteration counter. Algorithm~\ref{alg:sim_an_single} provides a detailed description of this procedure.
\begin{algorithm}[t]
\caption{Simulated Annealing procedure to solve problem~(\ref{ising_assort_single})}
\label{alg:sim_an_single}
\footnotesize
\algloop{Input}
\algblock{Input}{EndIf}
\algloop{Output}
\algblock{Output}{EndIf}
\algloop{Initialization}
\algblock{Initialization}{EndIf}
\begin{algorithmic}[1]
\Input
    \State $\mathcal{N}$: whole product set
    \State $S_{start}$: assortment representing the starting point of the algorithm
    \State $d_{obj}$: typical increase of the objective function
    \State $k_{temps}$: number of temperatures
    \State $p_{min}$, $p_{max}$: minimum and maximum target acceptance probabilities
\EndIf
\Output
    \State $S_{heur}$: optimized assortment
\EndIf
\Initialization

    \State $S_{cur} \gets S_{start}$, \ $R_{cur}$ $\gets$ 
    \texttt{estimateProfit}($S_{cur}$)
    \Comment{Initialize the current solution and its value}
    \State $S_{heur} \gets S_{cur}$, \ $R_{heur}$ $\gets$ $R_{cur}$
    \Comment{Initialize the heuristic solution and its value}
    \State $T \gets -d_{obj} / \log(p_{max})$ 
    \Comment{Initialize the temperature}
\EndIf
\For{$i$ \textbf{in} $[1, \dots, k_{temps}]$}:
\Comment{Main loop}

    \State Draw a random product $j \in \mathcal{N}$
    
    \If{$j \in S_{cur}$}:
    \Comment{Generate the candidate solution}
    \State $S_{can} \gets S_{cur} \backslash \{j\}$
    \Else:
    \State $S_{can} \gets S_{cur} \cup \{j\}$
    \EndIf
    
    \State $R_{can}$ $\gets$ \texttt{estimateProfit}($S_{can}$)
    \Comment{Estimate the value of the candidate solution}   
    \If{$R_{can} > R_{cur}$}:
    \Comment{Update the current solution}
    \State $S_{cur} \gets S_{can}$, \ $R_{cur}$ $\gets$ $R_{can}$
    \Else:
    \State Draw a random number $r \sim U(0,1)$
    \If{$r < \exp\bigl((R_{can} - R_{cur})/T\bigr)$}:
    \State $S_{cur} \gets S_{can}$, \ $R_{cur}$ $\gets$ $R_{can}$
    \EndIf
    \EndIf
    
    \If{$R_{cur} > R_{heur}$}:
    \Comment{Update the heuristic solution}
    \State $S_{heur} \gets S_{cur}$, \ $R_{heur}$ $\gets$ $R_{cur}$
    \EndIf
    \State $T \gets -d_{obj} / \log(p_{max} + (p_{min} - p_{max}) i /k_{temps})$
\Comment{Update the temperature}
\EndFor
\State \Return{$S_{heur}$, $R_{heur}$}
\end{algorithmic}
\end{algorithm}

To evaluate the performance of the developed Simulated Annealing algorithm, we will compare it with several benchmark algorithms. These algorithms follow a similar structure: First, a weight $w_j$ is assigned to each product $j \in \mathcal{N}$. Next, the products are ranked by weight in decreasing order, resulting in the sequence $j_1, j_2, j_3, \ldots, j_n$. We then examine all possible assortments of the form $j_1, j_2, \ldots, j_k$, and select $k$ that results in an assortment with the highest expected profit. We consider three methods to compute the weights:
\begin{itemize}
    \item {\bfseries Revenue Weights:} $w_j = r_j$, where $r_j$ is the gross profit margin (alternatively, the revenue) associated with product $j$. This approach results in the well-known revenue-ordered heuristic algorithm (see, e.g., \citealp{Berbeglia1}). 
    \item {\bfseries Parameter Weights:} $w_j = r_j \exp\bigl(\theta_{jj} + \sum_{i \ne j}\theta_{ij}\bigr)$.  The idea behind this approach is to account for each product's role in the product portfolio by incorporating Ising model parameters.
    \item {\bfseries Katz Weights:} $w_j = r_j C_{Katz}(j)$, where $C_{Katz}(j)$ is the Katz centrality measure associated with node $j$ in the graphical representation of the Ising model. The idea is similar to the previous point, but here we employ a more advanced approach specifically designed to exploit the graphical representation of the Ising model and network analytics techniques.
    
\smallskip
The Katz centrality score for product $j$ satisfies:
\begin{equation*}
C_{Katz}(j) = \alpha \sum_{i \ne j} \theta_{ij} C_{Katz}(i) + \beta,
\end{equation*}
where parameter $\alpha>0$ controls the influence of distant neighbors, and parameter $\beta>0$ ensures that all nodes have a nonzero centrality  (thus contributing to centrality calculation for other nodes). The value of $\alpha$ is typically set to be less than $1/\lambda_{\text{max}}$ --  where $\lambda_{\text{max}} $ is the largest eigenvalue of the weighted graph's adjacency matrix -- for the Katz centrality calculation to converge (see \citealp{Newman1}).
\smallskip \\
Katz centrality inherently accounts for both direct and indirect relationships, amplifying the centrality of nodes (products) that are well-connected, especially to other ``central'' nodes. 
While Katz centrality typically assumes positive edge weights, it can be adapted to accommodate negative weights, provided that the Katz centrality calculation converges. Negative weights will result in lower centrality scores of products that have many negative dependencies. However, the centrality scores become challenging to interpret in this context and must be treated with caution.
\end{itemize}
Importantly, it should be emphasized that since the underlying multi-purchase models in the literature prescribe vastly different choice behavior (e.g., shoppers purchasing bundles containing two products at most), which typically leads to completely different assortment optimization problems, there is no benchmark that would allow for a meaningful comparison. We test the comparative performance of the proposed heuristic algorithms in the following subsection.

\section{Numerical Analysis}
\label{2subsec:num_analysis}

In this section we perform two sets of numerical experiments. In the first set we evaluate the performance of the heuristic algorithms for the assortment optimization problem under the Ising model in terms of their solution quality and computing times. The aim of the second set is to show the merit of the assortment optimization approach based on the Ising model compared with benchmark assortment models from the existing literature.

\subsection{Performance of Heuristic Algorithms for the Ising Model}
We conduct numerical experiments using simulated data to gain insights into the performance of the heuristic algorithms described in the previous subsection. The use of simulated data as opposed to real data enables us to obtain detailed and robust numerical evidence of the comparative performance of the heuristics across a large number of (choice behavior) instances, rather than just for a single dataset. We start by describing the parameter generation procedure. To evaluate the quality of solutions, we set the product portfolio size to $n = 50$ and generate 100 problem instances. The generation process for parameters~$\theta$ is more clearly described in terms of the corresponding graphs $G$. Consider the following graph generation procedure:
\begin{enumerate}
\item[1.] Create nodes with weights $\theta_{ii}$ $\forall i \in \mathcal{N}$ sampled randomly from the interval $[2, 4]$;
\item[2.] For all unordered pairs $i,j \in \mathcal{N}$, create an edge between nodes $i$ and $j$ with probability~$p_{edge} = 0.2$;
\item[3.] For each edge, sample its absolute weight $\theta_{ij}$ from the interval $[1, 2]$ and multiply the obtained value by $-1$ with probability $p_{neg}=0.8$.
\end{enumerate}

The chosen values of $p_{edge}$ and $p_{neg}$ imply that product substitution effects between pairs of adjacent products (i.e., products connected via an edge) are somewhat prevalent compared to complementarity effects. Nevertheless, direct product substitution effects exist with probability $p_{edge} \cdot p_{neg} = 0.16$, which is a reasonable value for illustrative purposes. Finally, for each graph~$G$, we randomly sample the gross profit values $r_j$, $j\in \mathcal{N}$ from the interval $[0.01, 1]$. 

\begin{figure}[t]
\centering
\begin{subfigure}{.485\textwidth}
  \centering
  \includegraphics[width=\linewidth]{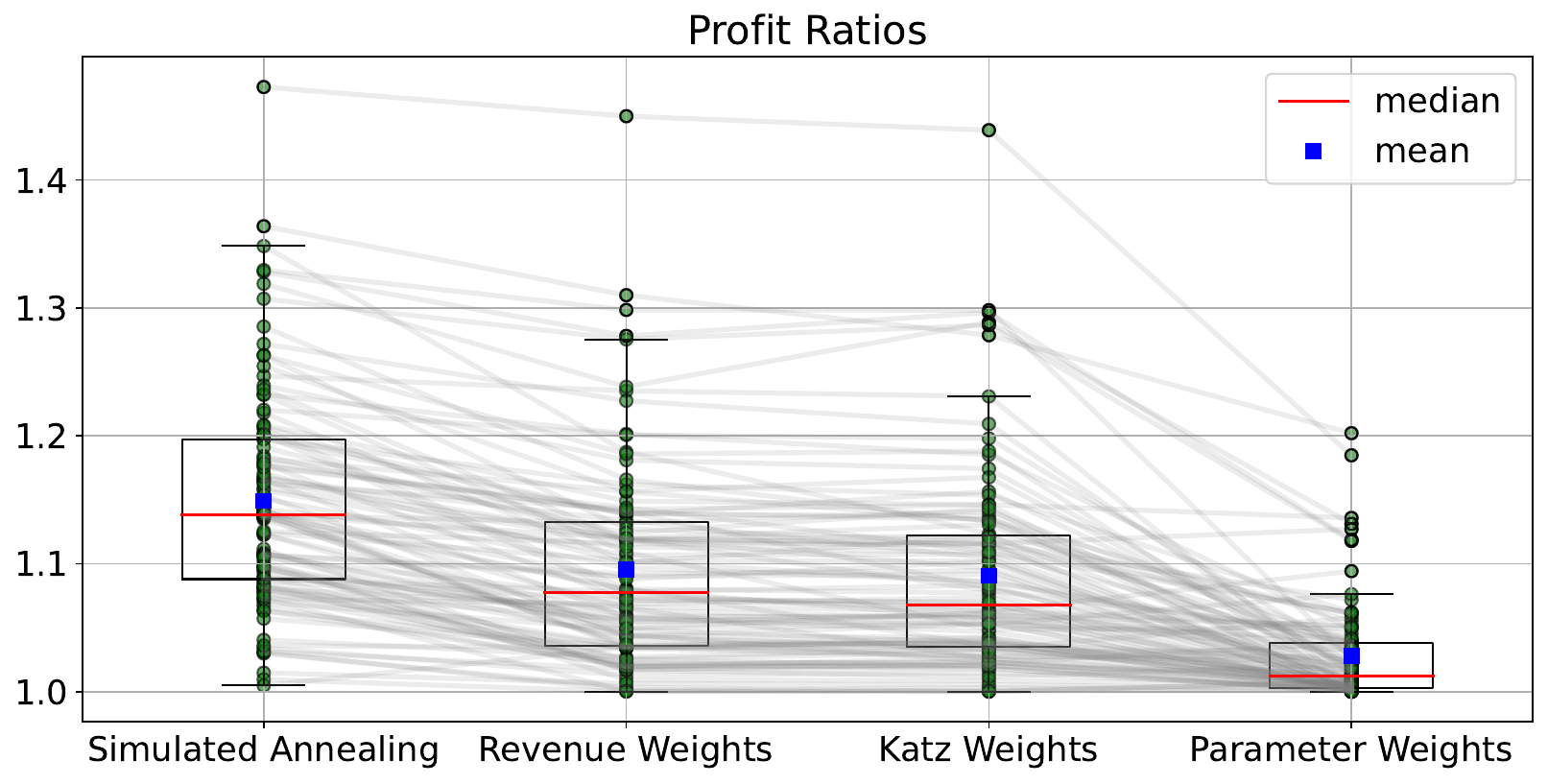} 
  \caption{\footnotesize{Profit gains from optimized assortments compared to unoptimized assortments (all products included).}}
\end{subfigure}
\quad
\begin{subfigure}{.48\textwidth}
  \centering
  \includegraphics[width=\linewidth]{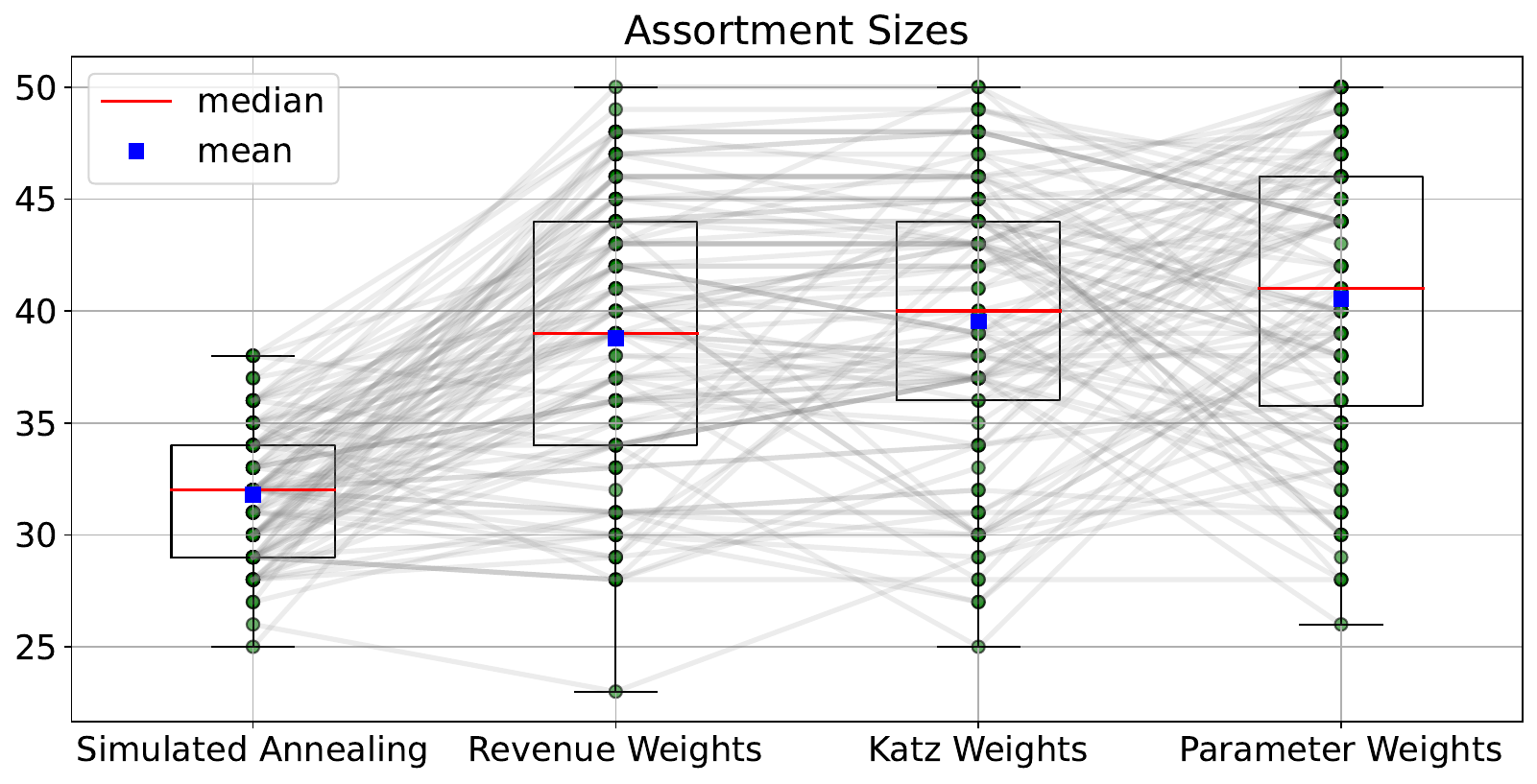}
  \caption{\footnotesize{Assortment sizes.\\ \quad}}
\end{subfigure}
\caption{Comparison of heuristic algorithms.}
\vspace{-0.1cm}
\label{ao_scomparison_all}
\end{figure}

\begin{figure}[t]
\centering
\includegraphics[width=0.55\linewidth]{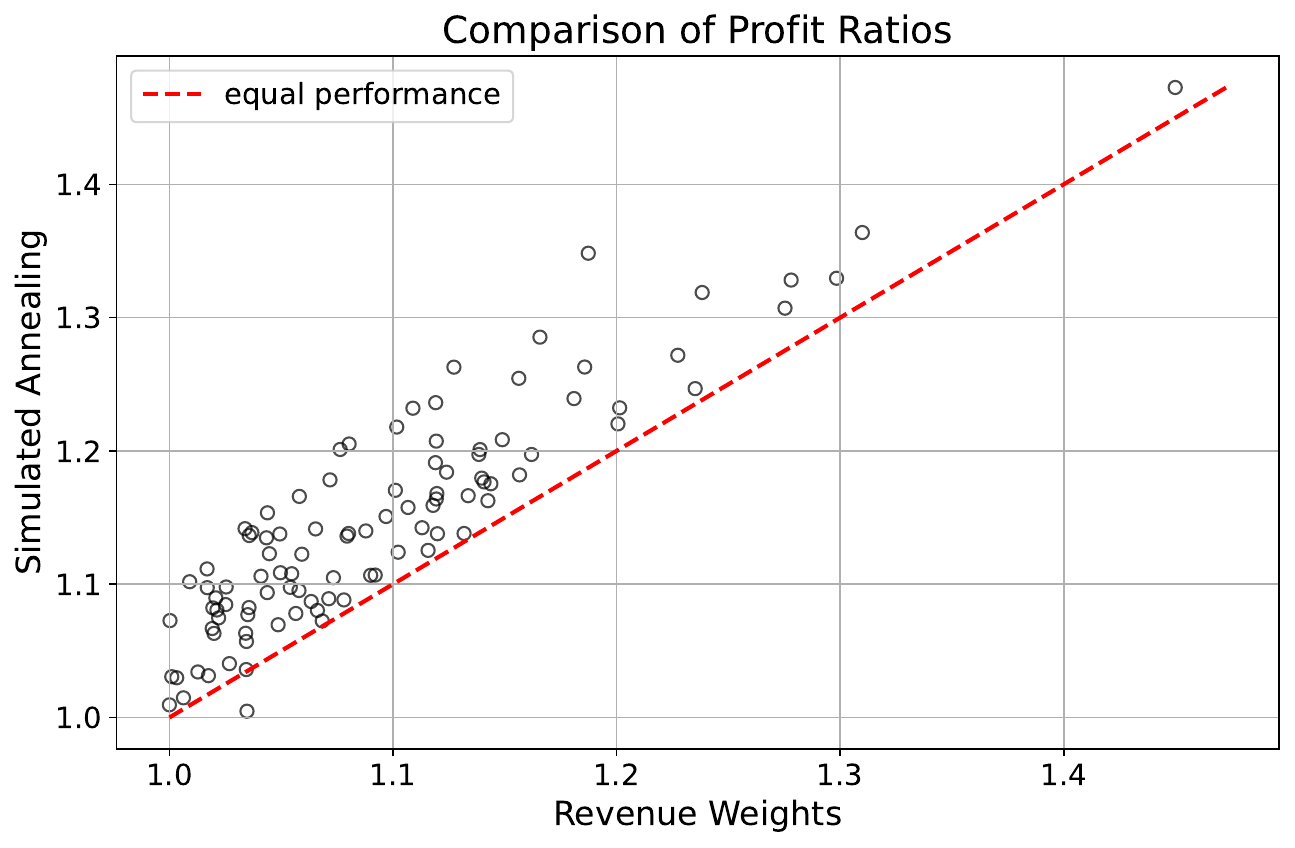}
\caption{Comparison of profit ratios for Simulated Annealing and revenue-ordered heuristics.}
    \vspace{-0.1cm}
    \label{ao_scomparison_pair}
\end{figure}

Figure~\ref{ao_scomparison_all} contains the performance statistics across the generated instances for all four heuristic algorithms. 
For the Simulated Annealing algorithm, we used the following input parameters: $p_{min} = 0.001$, $p_{max} = 0.999$, $k_{temps} = 10,000$, $d_{obj} = 0.25$, and $S_{start} = \mathcal{N}$. For the Katz Weights algorithm, we set $\beta = 1$ and $\alpha = 1/\lambda_{\text{max}} - 0.01$, where $\lambda_{\text{max}}$ is the largest eigenvalue of the corresponding weighted graph's adjacency matrix. The number of samples used to estimate the profits yielded by candidate assortments was set to $10,000$ for all four algorithms. Figure~\ref{ao_scomparison_all}(a) illustrates the profit gains yielded by optimized assortments compared to unoptimized assortments (where all products are included), and Figure~\ref{ao_scomparison_all}(b) shows the sizes of the optimized assortments. The Simulated Annealing algorithm significantly outperforms the benchmark algorithms, achieving an average profit gain of $14.9\%$, followed by the Revenue Weights algorithm at $9.5\%$, the Katz Weights algorithm at $9.1\%$, and the Parameter Weights algorithm at $2.8\%$. The revenue-ordered heuristic algorithm shows a surprisingly strong performance. It outperforms the Katz Weights and Parameter Weights algorithms, highlighting the complexity of the demand structure in the considered problem setting. Moreover, the better-performing algorithms achieve higher expected profits while leading to smaller average assortment sizes: $31.8$ units for Simulated Annealing, $38.8$ units for Revenue Weights, $39.5$ units for Katz Weights, and $40.5$ units for Parameter Weights. Furthermore, despite similar average performance in terms of expected profits, Revenue Weights and Katz Weights often result in noticeably different assortment sizes for the same generated problem instances. Finally, a case-by-case comparison of the profit gains yielded by the two leading methods -- Simulated Annealing and Revenue Weights -- is provided in Figure~\ref{ao_scomparison_pair}. 
We can observe that Simulated Annealing demonstrates superior performance in $99$ of the $100$ generated instances.

\begin{figure}[t]
\centering
\includegraphics[width=0.5\linewidth]{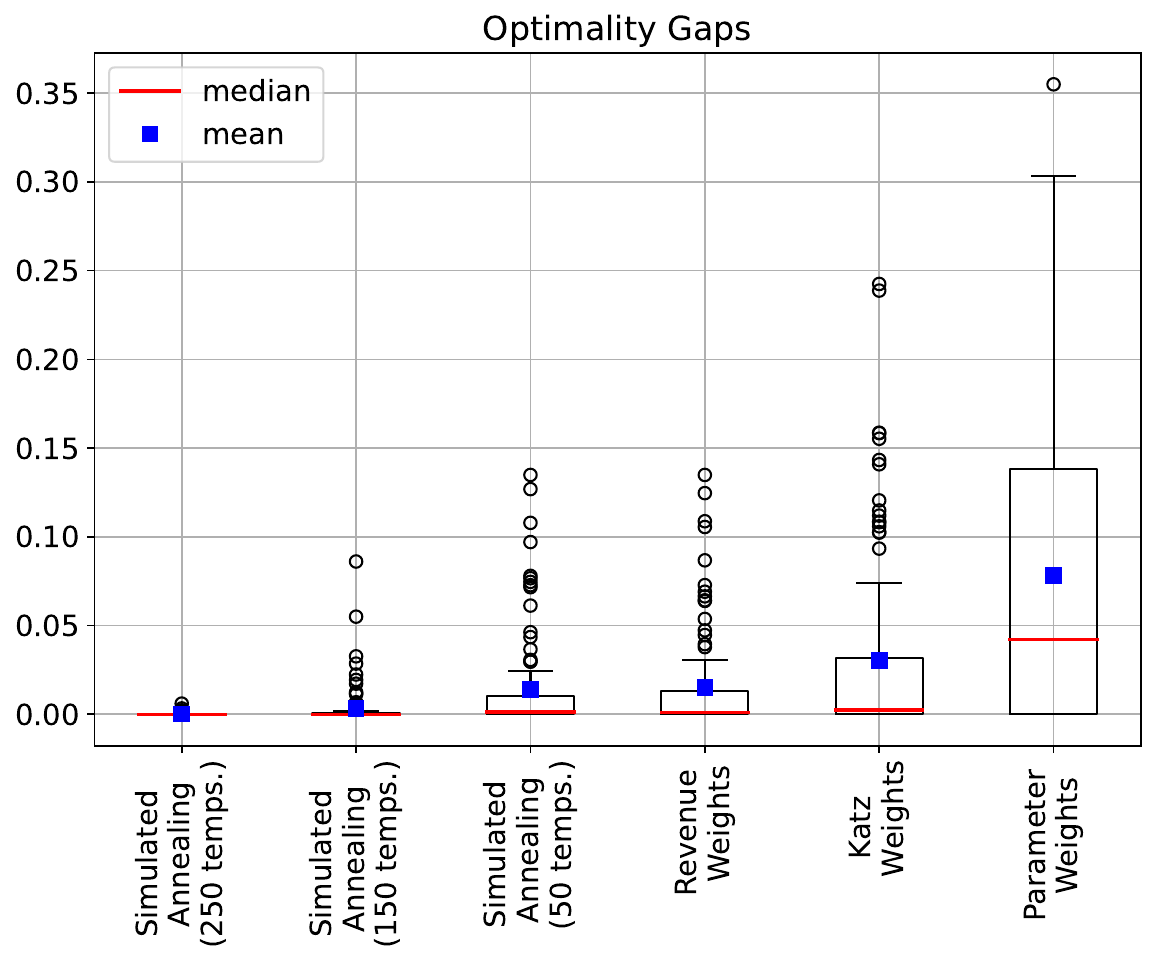}
\caption{Optimality gaps for heuristic algorithms for small-size problems.}
    \label{optim_gaps}
\end{figure}

Next, for small problem instances, we can effectively evaluate the optimality gaps between the values of our proposed heuristic algorithms and the optimal objective values obtained through brute-force search. We set the product portfolio size at $n=10$ and generate $100$ problem instances using the procedure described at the beginning of this subsection, with the exception that we use $p_{\text{edge}} = 1$ (i.e., we generate complete graphs). This adjustment ensures that the generated problem instances represent complex shopping behaviors, which might be difficult to achieve with sparse graphs in small-scale settings. We then solve the corresponding assortment optimization problems using the described heuristic algorithms. Importantly, for the considered small problem instances, the expected profit for each candidate assortment can be calculated exactly without resorting to Gibbs sampling. Furthermore, Simulated Annealing's performance hinges on the number of temperatures (which, in the case of Algorithm~\ref{alg:sim_an_single}, is equivalent to the number of iterations). If the number of temperatures is sufficiently high, the algorithm will eventually explore the entire solution space and find the optimal solution, making the performance analysis redundant. Therefore, we report results for three variants of Algorithm~\ref{alg:sim_an_single}: with $k_{\text{temps}} = 250$, $k_{\text{temps}} = 150$, and $k_{\text{temps}} = 50$. The values of the optimality gaps are provided in Figure~\ref{optim_gaps}.
The results indicate that with $k_{\text{temps}} = 250$, Simulated Annealing almost always finds the optimal solution, leading to an average optimality gap of 0.02\%. Similarly, when $k_{\text{temps}} = 150$, the average optimality gap remains very low at 0.3\%. At $k_{\text{temps}} = 50$, Simulated Annealing's performance is similar to that of the revenue-ordered heuristic, resulting in optimality gaps of 1.4\% and 1.5\%, respectively. The Katz Weights algorithm demonstrates somewhat worse performance, with an average optimality gap of 3\%. Finally, the Parameter Weights algorithm leads to a relatively high optimality gap of 7.8\%.
These results are consistent with our comparative analysis on larger problems, demonstrating that the developed Simulated Annealing procedure outperforms other considered algorithms, although the Revenue Weights algorithm is remarkably effective despite its simplicity.

\begin{figure}[t]
\centering
\begin{subfigure}{.31\textwidth}
  \centering
  \includegraphics[width=\linewidth]{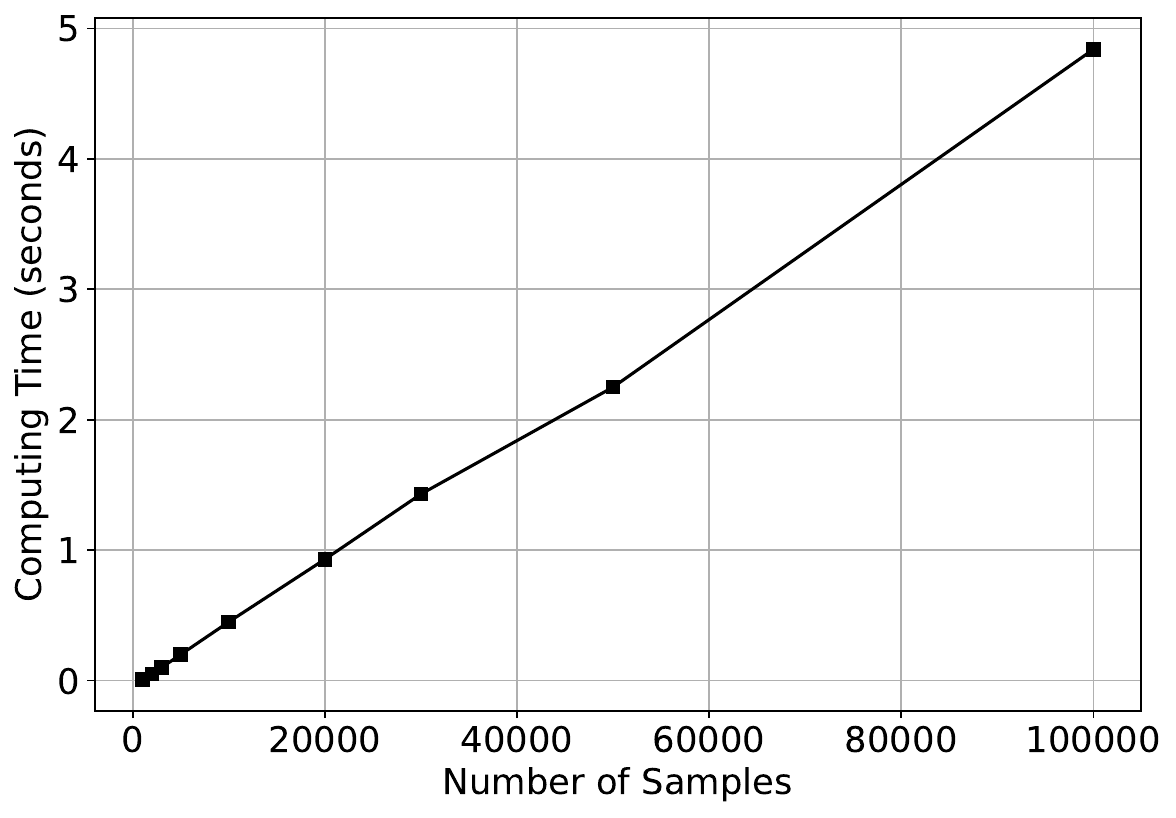}
  \caption{\footnotesize{Computing time depending on number of samples.}}
\end{subfigure}
\begin{subfigure}{.32\textwidth}
  \centering
  \includegraphics[width=\linewidth]{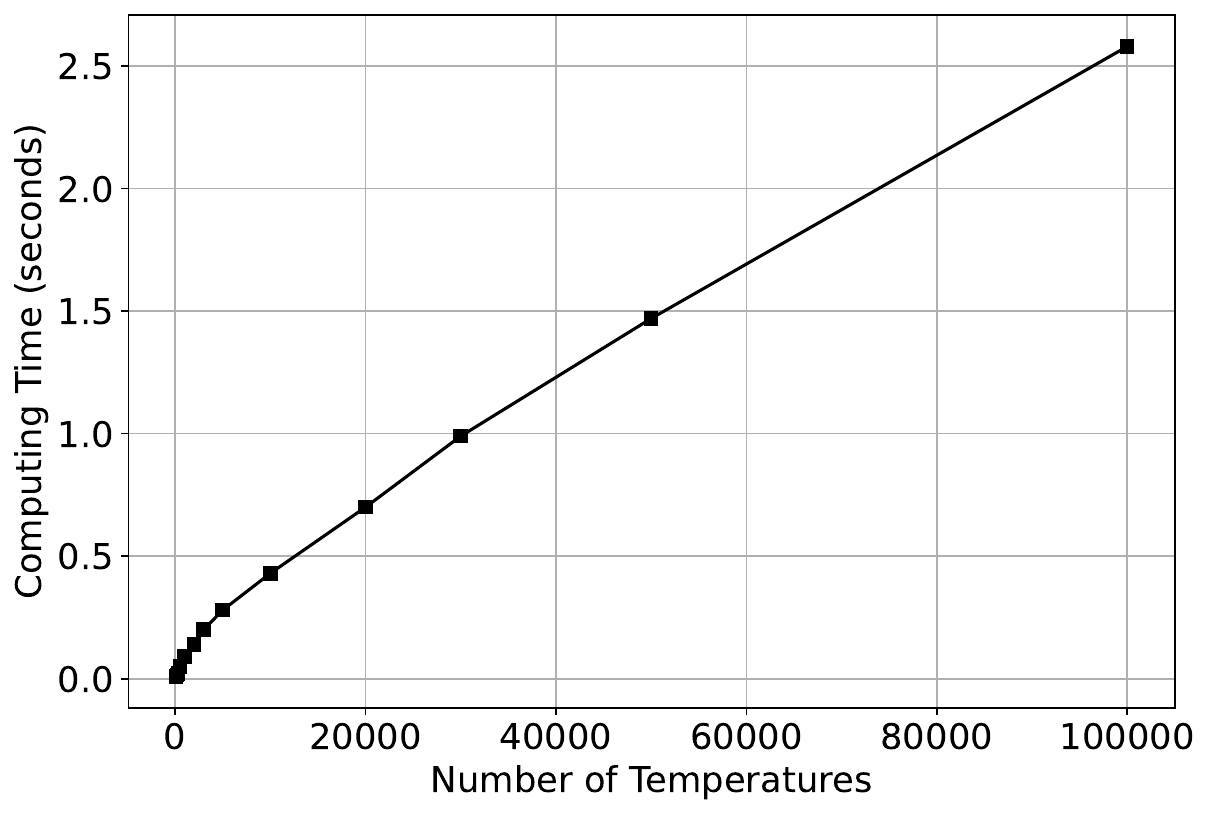}
  \caption{\footnotesize{Computing time depending on number of temperatures.}}
\end{subfigure}
\begin{subfigure}{.327\textwidth}
  \centering
    \vspace{0.3cm}
  \includegraphics[width=\linewidth]{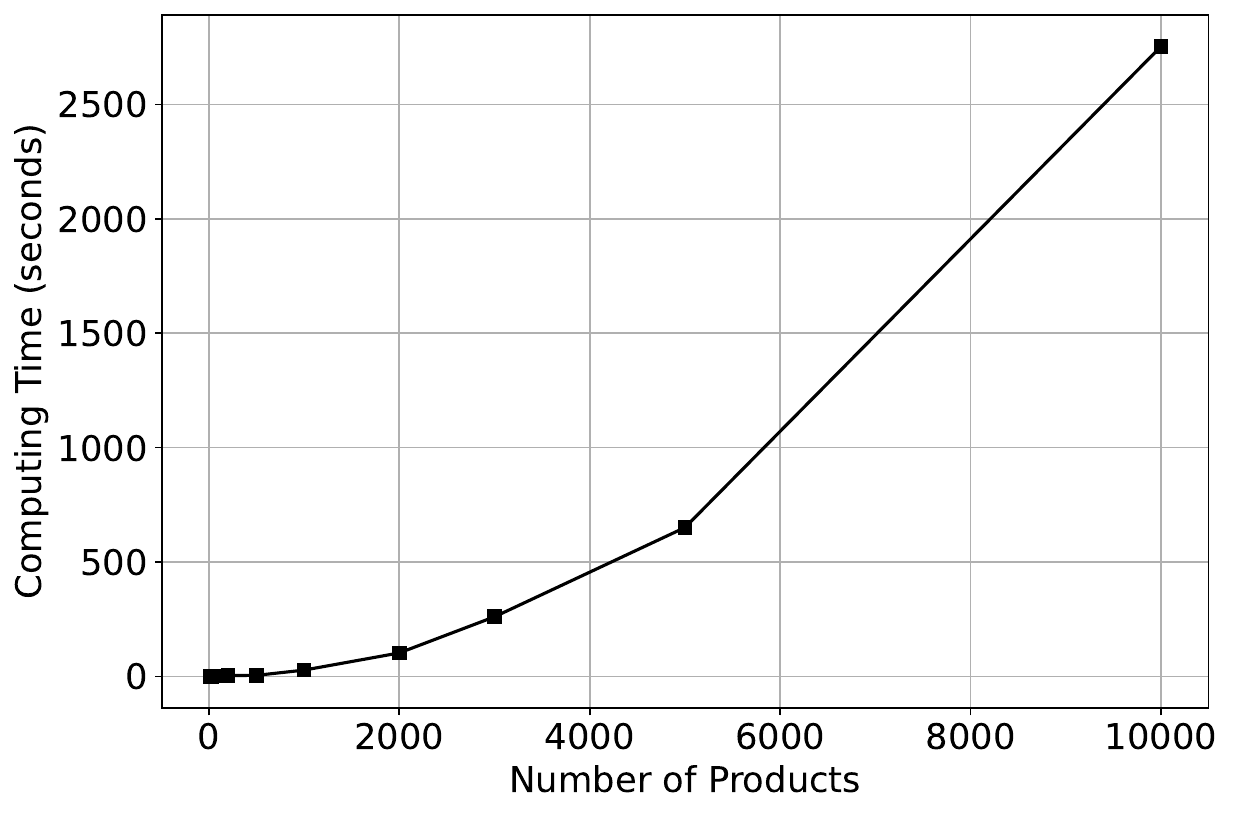}
  \caption{\footnotesize{Computing time depending on number of products.}}
\end{subfigure}
\caption{Computing time statistics for Simulated Annealing algorithm.}
\vspace{-0.1cm}
\label{ao_times}
\end{figure}

The primary challenge associated with the described Simulated Annealing algorithm is the required computing time, which depends on three factors: the size of the product portfolio, the number of temperatures, and the number of samples used to estimate the profit yielded by a candidate assortment. We analyze the impact of these parameters on computing time in Figure~\ref{ao_times}.
The parameters were generated following the procedure described previously in this subsection. We first set $n=10$, $k_{temps} = 100$, and varied the number of samples $k_{samples}$ from $100$ to $100,000$ as shown in Figure~\ref{ao_times}(a). We then set $n=10$, $k_{samples} = 100$, and varied $k_{temps}$ from $10$ to $100,000$ as shown in Figure~\ref{ao_times}(b). Finally, we set $k_{temps} = 100$, $k_{samples} = 100$, and varied $n$ from $10$ to $10,000$ as shown in Figure~\ref{ao_times}(c). All numerical experiments were carried out on a MacBook Pro with an Apple M2 Pro chip and 16 GB RAM. The algorithm was implemented using Python version 3.10.9.

We observe that the computing time increases linearly in the number of samples and temperatures, as illustrated in Figures~\ref{ao_times}(a) and (b). Additionally, it increases exponentially in the number of products, as shown in Figure~\ref{ao_times}(c). It is important to note that as the size of the product portfolio grows, the number of temperatures and samples must also increase to maintain the efficacy of the Simulated Annealing procedure. For this reason, the developed algorithm becomes less attractive for large-scale environments -- typically starting from a couple of hundred products, depending on available computing resources. Therefore, for such cases, the option of separating the product portfolio into independent segments should be explored, as it would allow for solving the assortment optimization problem for these segments individually. If this is not feasible, adopting the revenue-ordered heuristic algorithm is a reasonable alternative, given its surprisingly good performance under the Ising model.

\subsection{Ising Model versus Other Single/Multi-purchase Models}

To further demonstrate the benefit of using the Ising model, which allows for basket choices of any size, we compare our approach with two benchmark models -- namely, the classic single-choice MNL and the multi-purchase BundleMVL-2 considered by \cite{Tulabandhula1}. We use the same small-instance setting as in the previous subsection, setting $n=10$ and generating 100 problem instances using the same parameter generation procedure. For each generated matrix $\theta$, we take the full product portfolio $\mathcal{N}$ as the offered set, enumerate all possible baskets, compute the exact Ising choice probabilities, and draw a sample of baskets from this exact distribution. The resulting sample is then treated as the sales history available to the benchmark models.

The two benchmark models require different treatment of the sampled baskets. For the BundleMVL-2 model, we follow the preprocessing procedure proposed in Appendix A.3 of \cite{Tulabandhula1}. Baskets containing at most two products are used directly, whereas larger baskets are transformed into weighted observations corresponding to all possible partitions into subsets of size at most two, with equal weights assigned to the partitions. The resulting weighted sample is then used in the maximum likelihood estimation problem for BundleMVL-2. For the MNL benchmark, each sampled basket is converted into singleton purchases, so that a basket containing $k$ products contributes $k$ single-item observations, after which the model parameters are estimated by maximum likelihood. In both cases, the no-purchase probability is calibrated separately for each benchmark model using a simple grid search such that the subsequent optimization performance (in terms of expected profit) of each benchmark assortment problem is maximized individually.

\begin{figure}[t]
	\centering
	\includegraphics[clip,width=0.55\linewidth]{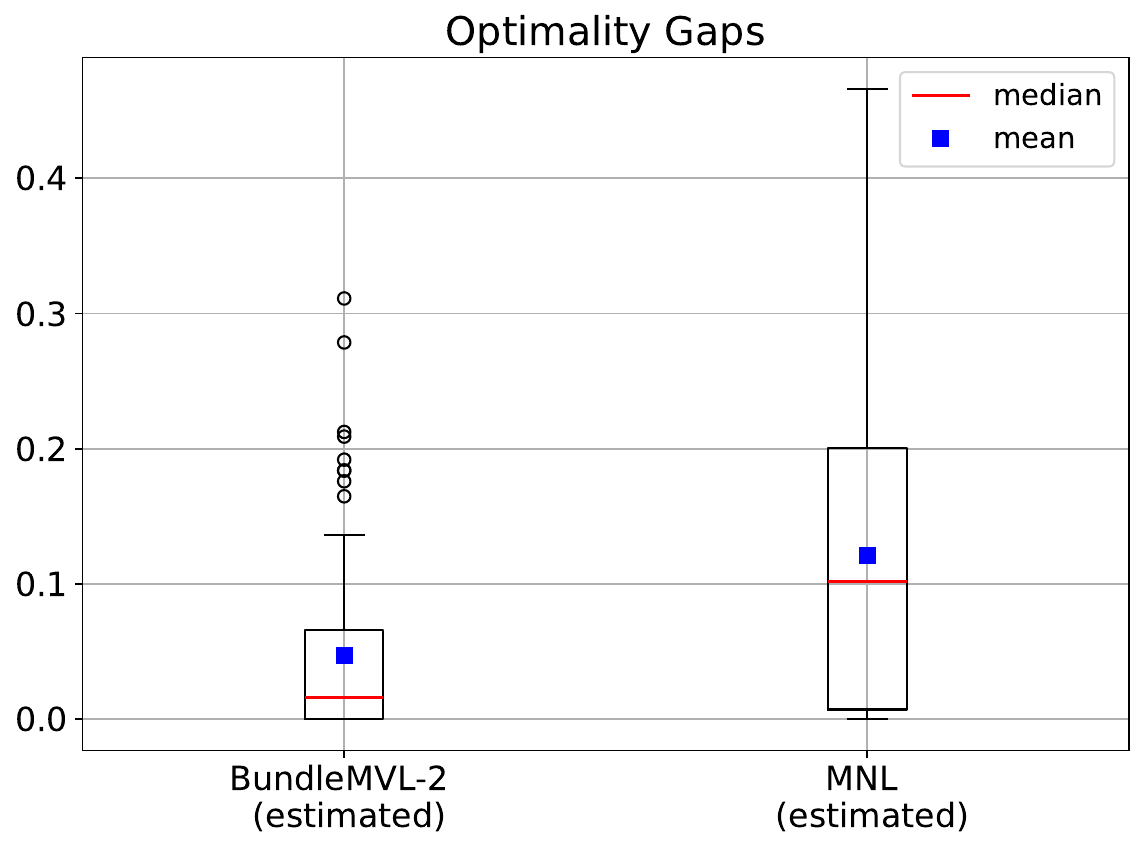}
	\caption{Optimality gap of BundleMVL-2 and MNL with Ising model as the ground truth.}	\label{fig:res1}
\end{figure}

Once the benchmark models are estimated, we solve the corresponding assortment optimization problems exactly over all subsets of $\mathcal{N}$ and then evaluate the resulting assortments under the true Ising model. Since the considered instances are small, all expected profits can be computed exactly, without resorting to Gibbs sampling and approximation techniques. We then compare the optimal benchmark assortments with the optimal assortment under the Ising model and compute the corresponding optimality gaps. The results shown in Figure \ref{fig:res1} indicate that the optimal assortments under the Ising model lead to considerably higher expected profits than those obtained under the two considered benchmark models. 
In fact, approximating multi-purchase choice behavior through the BundleMVL-2 and the MNL model leads to average optimality gaps of  4.7\% and 12.1\%, respectively. At the same time, BundleMVL-2 provides a substantially better approximation than the single-choice MNL, with a difference in average optimality gaps of almost 7.5\%.

It is important to note that the magnitude of these optimality gaps depends on the basket-size distribution. Under the baseline parameter generation procedure used in this subsection, the distribution of basket sizes across all generated problem instances and all sampled baskets has a mean of 4.07 and a standard deviation of 0.86. If the parameters are generated in such a way that most baskets contain at most two products, then the benchmark based on BundleMVL-2 becomes nearly optimal. This is not surprising, since BundleMVL-2 naturally accommodates this type of purchase behavior and has a formulation that is conceptually close to the Ising model. By contrast, when the generated baskets are larger on average, the performance of both BundleMVL-2 and MNL deteriorates substantially. Additional numerical results illustrating this sensitivity can be found in \ref{sens_bench_perfo_basket}. Overall, these findings provide further evidence that allowing for baskets of arbitrary size plays a crucial role in assortment optimization when customers tend to purchase products in large baskets.

\section{Conclusions and Future Work}
\label{2sec:concl}

In this paper, we addressed the important question of assortment optimization and its impact on firm profitability in the context of customer basket shopping behavior. To tackle this problem, we developed a comprehensive methodology at the interface of operations research and computer science. We brought into the spotlight the equivalence relation between the Ising model and the multivariate logit model. This allowed us to leverage powerful theoretical results devised for the Ising model for various purposes, including parameter estimation and deriving complexity results for the assortment optimization problem. Furthermore, we introduced a preprocessing procedure that can be used to reduce the dimensionality of the assortment optimization problem based on the graphical representation of the Ising model. 
To solve the assortment optimization problem after completing the preprocessing stage, we proposed several heuristic algorithms. We then carried out an extensive numerical analysis to compare the performance of these algorithms. In our experiments, the developed Simulated Annealing procedure significantly outperformed the other heuristic algorithms, leading to an expected profit gain of 15\%. Nevertheless, the revenue-ordered heuristics showed surprisingly good performance as well -- with a 10\% expected profit gain -- proving to be an efficient alternative if applying Simulated Annealing is not computationally feasible.

Our work motivates a number of interesting research questions. We believe that the graphical representation of the Ising model can be further exploited to boost assortment optimization capabilities through network analytics techniques. For example, one can identify ``important'' products using different centrality measures, or ``established'' baskets using community detection algorithms. This information can be incorporated into the developed Simulated Annealing procedure, for example, in such a way that at each iteration, high centrality products are added to the assortment with a higher probability and removed from the assortment with a lower probability. Another interesting avenue is to explore different weighting schemes in the weight-ordered heuristics,  potentially incorporating centrality measures other than Katz centrality. Last but not least, more general MRFs with other joint distribution functions could be used for more accurate modeling of basket shopping behavior (e.g., by enabling our model to account for product quantities in each basket). Ultimately, our work advances the development of a comprehensive framework for assortment optimization that takes into account basket shopping behavior, providing a base for future studies on assortment decision-making.

\section*{Acknowledgments}

The work of the first and the last author was supported by the Swiss National Science Foundation (SNSF) under Grant 192545. Maier's work was supported by the Engineering and Physical Sciences Research Council [grant number UKRI3639].

\appendix
\renewcommand{\thesection}{Appendix \Alph{section}}

\section{Proof of Proposition~\ref{prop_transf_equiv}}
\label{prop_transf_equiv_proof}
Let $x$ be a binary vector representing a basket of products from $\mathcal{N}$. Furthermore, let $\tilde{x}$ be the corresponding spin vector, i.e. the vector such that its entries satisfy
\begin{equation*}
\begin{aligned}
\tilde{x}_i  = &\begin{cases}
1&\text{if $x_i = 1$},\\
-1&\text{otherwise.}
\end{cases}
\end{aligned}
\end{equation*}
Note that $x = (\tilde{x} + 1)/2$. Since $\theta$ and $\tilde{\theta}$ satisfy conditions~(\ref{par_transform}), we obtain that
\begin{equation*}
\begin{aligned}
    \sum_{i \in \mathcal{N}}\theta_{ii}x_i + \sum_{i,j\in \mathcal{N}: \, i \ne j}x_i\theta_{ij}x_j &= \sum_{i\in \mathcal{N} }\Bigl(2\tilde{\theta}_{ii} - 4\sum\limits_{j\in \mathcal{N}: \, j \ne i}  \tilde{\theta}_{ij}\Bigr)\dfrac{\tilde{x_i} + 1}{2} + \sum_{i,j\in \mathcal{N}: \, i \ne j} \dfrac{\tilde{x_i} + 1}{2} 4\tilde{\theta}_{ij} \dfrac{\tilde{x_j} + 1}{2}, \\
    &= \sum_{i \in \mathcal{N}}\Bigl(\tilde{\theta}_{ii} - 2\sum\limits_{j\in \mathcal{N}: \, j \ne i}  \tilde{\theta}_{ij}\Bigr)(\tilde{x_i} + 1) + \sum_{i,j\in \mathcal{N}: \,i \ne j} \tilde{\theta}_{ij}(\tilde{x_i} \tilde{x_j} + \tilde{x_i} + \tilde{x_j} + 1) \\
    &= \sum_{i\in \mathcal{N}}\tilde{\theta}_{ii}\tilde{x_i}+ \sum_{i,j\in \mathcal{N}: \,i \ne j} \tilde{x_i} \tilde{\theta}_{ij} \tilde{x_j} + \sum_{i\in \mathcal{N} }\Bigl(\tilde{\theta}_{ii} - \sum\limits_{j\in \mathcal{N}: \, j \ne i}  \tilde{\theta}_{ij}\Bigr). 
\end{aligned}
\end{equation*}
Therefore,
\begin{equation*}
 \exp\biggl(\sum_{i\in \mathcal{N}}\theta_{ii}x_i + \sum_{i,j\in \mathcal{N}: \, i \ne j}x_i\theta_{ij}x_j\biggr) = \exp\biggl(\sum_{i\in \mathcal{N} }\tilde{\theta}_{ii}\tilde{x_i}+ \sum_{i,j\in \mathcal{N}: \,i \ne j} \tilde{x_i} \tilde{\theta}_{ij} \tilde{x_j}\biggr) \cdot const,
\end{equation*}
from which it follows that $ p_{\theta}(x \vert \mathcal{N}) = p_{\tilde{\theta}}(\tilde{x} \vert \mathcal{N})$, which concludes the proof.
\qed

\section{Proof of Theorem~\ref{prop_assortm_np_hard}}
\label{prop_assortm_np_hard_proof}
The decision version of the assortment optimization problem under the Ising model is as follows:
\begin{equation*}
\label{D-AO}
\tag{D-AO}
\textit{For any given $K$, is there assortment $S$ such that $R(S) \geq K$?}
\end{equation*}
Recall that $R(S)$ is defined by the following expression:
\begin{equation*}
    R(S) = \dfrac{\sum\limits_{x \in \mathcal{X}(S)}\exp\Bigl(\sum\limits_{i \in S}\theta_{ii}x_i + \sum\limits_{i \ne j, \ i, j \in S}x_i\theta_{ij}x_j\Bigr)\sum\limits_{j \in S} r_{j} x_j}{\sum\limits_{x \in \mathcal{X}(S)}\exp\Bigl(\sum\limits_{i \in S}\theta_{ii}x_i + \sum\limits_{i \ne j, \ i, j \in S}x_i\theta_{ij}x_j \Bigr)}.
\end{equation*}
Therefore, the inequality in the decision problem~(\ref{D-AO}) can be rewritten in the following way:
\begin{equation}
\label{ineq_rewrit}
    \begin{aligned}
    \sum\limits_{\substack{x \in \mathcal{X}(S),\\ x\ne 0}}\exp\Bigl(\sum_{i \in S}\theta_{ii}x_i + \sum_{i \ne j, \ i, j \in S}x_i\theta_{ij}x_j\Bigr)\biggl(\sum\limits_{j \in S} r_{j} x_j - K\biggr)\geq K.
    \end{aligned}
\end{equation}
Let $r_1 = K+1$, $\theta_{11} = 0$, and $r_j = 0$ $\forall j \in \mathcal{N} \backslash \{1\}$. Since product~$1$ is the only product with a nonzero profit margin, it has to belong to $S$ in order for inequality~(\ref{ineq_rewrit}) to be satisfied assuming that $K>0$. Then, inequality~(\ref{ineq_rewrit}) takes the following form:
\begin{equation*}
    \begin{aligned}
    \sum\limits_{\substack{x \in \mathcal{X}(S),\\% x\ne 0,\\
    x_1 = 1}}\exp\Bigl(\sum_{i \in S\backslash\{1\}}(\theta_{ii} + 2 \theta_{1i}) x_i + \sum_{\substack{i \ne j, \\ i, j \in S\backslash\{1\}}}x_i\theta_{ij}x_j\Bigr)\geq K.
    \end{aligned}
\end{equation*}
Lastly, let $\theta_{ii}' = \theta_{ii} + 2\theta_{1i}$ $\forall i \in \mathcal{N}\backslash\{1\}$ and $\theta_{ij}' = \theta_{ij}$ $\forall i,j \in \mathcal{N}\backslash\{1\}$. We can see that solving problem~(\ref{D-AO}) requires answering the question of whether the partition function of the Ising model with parameters $\theta'$ defined over nodes $\mathcal{N} \backslash \{1\}$ is greater than or equal to a given constant. Such a problem 
is NP-hard as shown by \cite{Istrail1} for nonplanar graphs, meaning that problem~(\ref{D-AO}) is NP-hard as well. \qed

\section{Proof of Proposition~\ref{ising_opt_proper}}
\label{ising_opt_proper_proof}
Let product~$k$ be such that $\theta_{kj} = 0$ $\forall j \in \mathcal{N} \backslash \{k\}$. Then, our goal is to show that $R(S \cup \{k\}) \geq R(S)$ for any assortment $S \subseteq \mathcal{N}$ such that $k \notin S$. This can be verified directly:
\begin{align*}
R(S\cup\{k\}) &= \dfrac{\sum\limits_{x \in \mathcal{X}(S\cup\{k\})}\exp\Bigl(\sum\limits_{i \in S\cup\{k\}}\theta_{ii}x_i + \sum\limits_{i \ne j, \ i, j \in S\cup\{k\}}x_i\theta_{ij}x_j\Bigr)\sum\limits_{j \in S\cup\{k\}} r_{j} x_j}{\sum\limits_{x \in \mathcal{X}(S\cup\{k\})}\exp\Bigl(\sum\limits_{i \in S\cup\{k\}}\theta_{ii}x_i + \sum\limits_{i \ne j, \ i, j \in S\cup\{k\}}x_i\theta_{ij}x_j \Bigr)} \\ 
&\geq  \dfrac{\sum\limits_{x \in \mathcal{X}(S\cup\{k\})}\exp\Bigl(\sum\limits_{i \in S\cup\{k\}}\theta_{ii}x_i + \sum\limits_{i \ne j, \ i, j \in S\cup\{k\}}x_i\theta_{ij}x_j\Bigr)\sum\limits_{j \in S} r_{j} x_j}{\sum\limits_{x \in \mathcal{X}(S\cup\{k\})}\exp\Bigl(\sum\limits_{i \in S\cup\{k\}}\theta_{ii}x_i + \sum\limits_{i \ne j, \ i, j \in S\cup\{k\}}x_i\theta_{ij}x_j\Bigr)}
\\ 
&=  \dfrac{\sum\limits_{x \in \mathcal{X}(S)}\bigl(1 + \exp(\theta_{kk} + 2\sum\limits_{i \in S} x_i\theta_{ki})\bigr)\exp\Bigl(\sum_{i \in S}\theta_{ii}x_i + \sum\limits_{i \ne j, \ i, j \in S}x_i\theta_{ij}x_j\Bigr)\sum\limits_{j \in S} r_{j} x_j}{\sum\limits_{x \in \mathcal{X}(S)}\bigl(1 + \exp(\theta_{kk} + 2\sum\limits_{i \in S} x_i\theta_{ki})\bigr)\exp\Bigl(\sum_{i \in S}\theta_{ii}x_i + \sum\limits_{i \ne j, \ i, j \in S}x_i\theta_{ij}x_j \Bigr)}
\\ 
&=  \dfrac{\sum\limits_{x \in \mathcal{X}(S)}\bigl(1 + \exp(\theta_{kk})\bigr)\exp\Bigl(\sum_{i \in S}\theta_{ii}x_i + \sum\limits_{i \ne j, \ i, j \in S}x_i\theta_{ij}x_j\Bigr)\sum\limits_{j \in S} r_{j} x_j}{\sum\limits_{x \in \mathcal{X}(S)}\bigl(1 + \exp(\theta_{kk})\bigr)\exp\Bigl(\sum_{i \in S}\theta_{ii}x_i + \sum\limits_{i \ne j, \ i, j \in S}x_i\theta_{ij}x_j \Bigr)}
\\ 
&=  \dfrac{\sum\limits_{x \in \mathcal{X}(S)}\exp\Bigl(\sum_{i \in S}\theta_{ii}x_i + \sum\limits_{i \ne j, \ i, j \in S}x_i\theta_{ij}x_j\Bigr)\sum\limits_{j \in S} r_{j} x_j}{\sum\limits_{x \in \mathcal{X}(S)}\exp\Bigl(\sum_{i \in S}\theta_{ii}x_i + \sum\limits_{i \ne j, \ i, j \in S}x_i\theta_{ij}x_j \Bigr)} = R(S). \qed
\end{align*}

\section{Proof of Proposition~\ref{ising_opt_proper3}}
\label{ising_opt_proper3_proof}
Let $S_{\mathcal{H}} = S \cap \mathcal{H}$ and $S_\mathcal{K} = S \cap \mathcal{K}$, and let us fix product $l \in S_{\mathcal{H}}$. To prove the proposition, it is sufficient to show that the marginal probability of a customer choosing product $l \in S_{\mathcal{H}}$ does not depend on $S_{\mathcal{K}}$. This is a direct implication of the fact that the Ising model satisfies global Markov properties. This can also be verified directly:
\begin{align*}
    &p_{\theta}(x_l = 1|S) =  \dfrac{\sum\limits_{x \in \mathcal{X}(S): \ x_l = 1}\exp\Bigl(\sum\limits_{i \in S}\theta_{ii}x_i + \sum\limits_{i \ne j, \ i, j \in S}x_i\theta_{ij}x_j\Bigr)}{\sum\limits_{x \in \mathcal{X}(S)}\exp\Bigl(\sum\limits_{i \in S}\theta_{ii}x_i + \sum\limits_{i \ne j, \ i, j \in S}x_i\theta_{ij}x_j \Bigr)} \\
    &= \dfrac{\sum\limits_{\substack{x \in \mathcal{X}(S_{\mathcal{H}}): \\ x_l = 1}}\biggl(\exp\Bigl(\sum\limits_{i \in S_{\mathcal{H}}}\theta_{ii}x_i + \sum\limits_{i \ne j, \ i, j \in S_{\mathcal{H}}}x_i\theta_{ij}x_j\Bigr)\sum\limits_{x' \in \mathcal{X}(S_{\mathcal{K}})}\exp\Bigl(\sum\limits_{i \in S_{\mathcal{K}}}\theta_{ii}x'_i + \sum\limits_{i \ne j, \ i, j \in S_{\mathcal{K}}}x'_i\theta_{ij}x'_j \Bigr)\biggr)}{\sum\limits_{x \in \mathcal{X}(S_{\mathcal{H}})}\biggl(\exp\Bigl(\sum\limits_{i \in S_{\mathcal{H}}}\theta_{ii}x_i + \sum\limits_{i \ne j, \ i, j \in S_{\mathcal{H}}}x_i\theta_{ij}x_j\Bigr)\sum\limits_{x' \in \mathcal{X}(S_{\mathcal{K}})}\exp\Bigl(\sum\limits_{i \in S_{\mathcal{K}}}\theta_{ii}x'_i + \sum\limits_{i \ne j, \ i, j \in S_{\mathcal{K}}}x'_i\theta_{ij}x'_j \Bigr)\biggr)} \\
    &=
    \dfrac{\sum\limits_{x \in \mathcal{X}(S_{\mathcal{H}}): \ x_l = 1}\exp\Bigl(\sum\limits_{i \in S_{\mathcal{H}}}\theta_{ii}x_i + \sum\limits_{i \ne j, \ i, j \in S_{\mathcal{H}}}x_i\theta_{ij}x_j\Bigr)}{\sum\limits_{x \in \mathcal{X}(S_{\mathcal{H}})}\exp\Bigl(\sum\limits_{i \in S_{\mathcal{H}}}\theta_{ii}x_i + \sum\limits_{i \ne j, \ i, j \in S_{\mathcal{H}}}x_i\theta_{ij}x_j \Bigr)} = p_{\theta}(x_l = 1|S_{\mathcal{H}}). \qed
\end{align*}

\section{Proof of Theorem~\ref{ising_opt_proper4}}
\label{ising_opt_proper4_proof}

Let us first formulate the following auxiliary lemma:
\begin{Lemma}
\label{aux_lem_1}
\begin{equation*}
    p_{\theta}(x_l = 1| S) \geq p_{\theta}(x_l = 1|x_k = 0, S) \iff p_{\theta}(x_l = 1|x_k = 1, S)  \geq p_{\theta}(x_l = 1|x_k = 0, S).
\end{equation*}
\end{Lemma}
\noindent {\itshape Proof of Lemma~\ref{aux_lem_1}.} Suppose that
\begin{equation}
\label{lem_ineq_1}
    p_{\theta}(x_l = 1| S) \geq p_{\theta}(x_l = 1|x_k = 0, S).
\end{equation}
Using the law of total probability, we obtain that:
\begin{equation*}
\begin{aligned}
    &p_{\theta}(x_l = 1|x_k = 0, S) p_{\theta}(x_k = 0| S) + p_{\theta}(x_l = 1|x_k = 1, S) p_{\theta}(x_k = 1| S) \geq p_{\theta}(x_l = 1|x_k = 0, S),
\end{aligned}
\end{equation*}
or
\begin{equation*}
     p_{\theta}(x_l = 1|x_k = 1, S) p_{\theta}(x_k = 1| S) \geq p_{\theta}(x_l = 1|x_k = 0, S)(1 -  p_{\theta}(x_k = 0| S)).
\end{equation*}
Since $p_{\theta}(x_k = 1| S) = 1 -  p_{\theta}(x_k = 0| S)$  and $p_{\theta}(x_k = 1| S) > 0$, we can divide both sides of the inequality by $p_{\theta}(x_k = 1| S)$ while preserving the inequality sign, thus obtaining the desired inequality:
\begin{equation}
\label{lem_ineq_2}
     p_{\theta}(x_l = j|x_k = 1, S)  \geq p_{\theta}(x_l = 1|x_k = 0, S).
\end{equation}
Finally, one can prove that inequality~(\ref{lem_ineq_1}) follows from inequality~(\ref{lem_ineq_2}) by repeating the above steps from the bottom up. \qed

\smallskip

\noindent {\itshape Proof of Theorem~\ref{ising_opt_proper4}.}  From Proposition~\ref{ising_opt_proper3} it follows that we can consider products in an isolated subgraph as a separate product portfolio $\mathcal{N}$. Suppose that $\theta_{ij}\geq0$ $\forall i,j \in \mathcal{N}$, $i\ne j$. To prove the theorem, it is sufficient to show that removing any product from the assortment can only reduce the marginal probabilities of customers choosing other products, i.e., for any $k, l \in \mathcal{N}$:
\begin{equation*}
p_{\theta}(x_l = 1| S) \geq p_{\theta}(x_l = 1|x_k = 0, S).
\end{equation*}
Without loss of generality, suppose that $S = \{1, \dots, m\}$, $l = 1$, and $k = 2$. As shown in Lemma~\ref{aux_lem_1}, the above inequality is equivalent to the following one:
\begin{equation}
\label{theor_start_point}
     p_{\theta}(x_1 = 1|x_2 = 1, S)  \geq p_{\theta}(x_1 = 1|x_2 = 0, S).
\end{equation}

Let $f(a,b) = \sum\limits_{x \in \mathcal{X}(S): \ x_1 = a, x_2=b}\exp\Bigl(\sum\limits_{i \in S}\theta_{ii}x_i + \sum\limits_{i \ne j, \ i, j \in S}x_i\theta_{ij}x_j\Bigr)$. Then, inequality~(\ref{theor_start_point}) can be rewritten as
\begin{equation*}
 \dfrac{f(1,1)}{f(0,1) + f(1, 1)} \geq \dfrac{f(1,0)}{f(0,0) + f(1, 0)},
\end{equation*}
or
\begin{equation}
\label{th3_main_ineq_short}
f(0,0)f(1,1) \geq f(1,0)f(0, 1).
\end{equation}

Let $g(x|S) = \exp\Bigl(\sum\limits_{i \in S}\theta_{ii}x_i + \sum\limits_{i \ne j, \ i, j \in S}x_i\theta_{ij}x_j\Bigr)$, and let $S^{-}_m = \{3, \dots, m\}$. Then, inequality~(\ref{th3_main_ineq_short}) takes the following form:
\begin{equation}
\label{th3_main_ineq}
\sum\limits_{x \in \mathcal{X}(S^{-}_m)}g\bigl[(0,0,x)|S\bigr] \sum\limits_{x \in \mathcal{X}(S^{-}_m)}g\bigl[(1,1,x)|S\bigr] \geq
\sum\limits_{x \in \mathcal{X}(S^{-}_m)}g\bigl[(1,0,x)|S\bigr] \sum\limits_{x \in \mathcal{X}(S^{-}_m)}g\bigl[(0,1,x)|S\bigr].
\end{equation}

We prove inequality~(\ref{th3_main_ineq}) by induction over $m$.

\smallskip

\noindent \textit{Induction base.} If $m=2$, then inequality~(\ref{th3_main_ineq}) can be written as:
\begin{equation*}
    1 \cdot \exp(\theta_{11} + \theta_{22} + 2 \theta_{12})  \geq \exp(\theta_{11}) \exp(\theta_{22}),
\end{equation*}
which holds true since $\theta_{12} \geq 0$. We also check that inequality~(\ref{th3_main_ineq}) holds if $m=3$. In this case, this inequality becomes:
\begin{equation*}
\begin{aligned}
    & \bigl(1 + \exp(\theta_{33})\bigr) \bigl(\exp(\theta_{11} + \theta_{22} + 2 \theta_{12}) + \exp(\theta_{11} + \theta_{22} + \theta_{33} + 2 \theta_{12} + 2 \theta_{13} + 2 \theta_{23})\bigr) \geq \\
    &\hspace{1cm} \bigl(\exp(\theta_{11}) + \exp(\theta_{11} + \theta_{33} + 2 \theta_{13})\bigr)\bigl(\exp(\theta_{22}) + \exp(\theta_{22} + \theta_{33} + 2 \theta_{23})\bigr).
\end{aligned}
\end{equation*}
It is easy to see that the above inequality holds if:
\begin{equation}
\label{induct_base_ineq}
    \exp(2 \theta_{12}) + \exp(2 \theta_{12} + 2 \theta_{13} + 2 \theta_{23}) \geq \exp( 2 \theta_{23}) + \exp( 2 \theta_{13}).
\end{equation}
Note that it is sufficient to show that the latter inequality holds for $\theta_{12} = 0$, in which case this inequality is equivalent to the following one:
\begin{equation*}
(\exp(2\theta_{13}) - 1)(\exp(2\theta_{23}) - 1) \geq 0.
\end{equation*}
The above inequality is true for any $\theta_{13}, \theta_{23} > 0$, meaning that inequality~(\ref{th3_main_ineq}) holds for $m=3$.

\smallskip
\noindent \textit{Induction step.} Assuming that inequality~(\ref{th3_main_ineq}) holds for dimension $m-1$, our goal is to show that it also holds for dimension $m$. This inequality can be rewritten in the following way:
\begin{equation*}
\sum\limits_{x, x' \in \mathcal{X}(S^{-}_m)}g\bigl[(0,0,x)|S\bigr] \cdot g\bigl[(1,1,x')|S\bigr] \geq
\sum\limits_{x, x' \in \mathcal{X}(S^{-}_m)}g\bigl[(1,0,x)|S\bigr] \cdot g\bigl[(1,0,x)|S\bigr],
\end{equation*}
which is equivalent to:
\begin{align*}
&\sum\limits_{x, x' \in \mathcal{X}(S^{-}_m)}g(x|S^{-}_m)  g(x'|S^{-}_m) \exp(\theta_{11} + \theta_{22} + 2\theta_{12} + 2\sum_{i \in S^{-}_m}(\theta_{1i} + \theta_{2i})x'_i) \geq \\
&\hspace{1cm}
\sum\limits_{x, x' \in \mathcal{X}(S^{-}_m)}g(x|S^{-}_m)  g(x'|S^{-}_m)\exp(\theta_{11} + \theta_{22} + 2\sum_{i \in S^{-}_m}\theta_{1i}x_i + 2\sum_{i \in S^{-}_m}\theta_{2i}x'_i).
\end{align*}
Since $\theta_{12} > 0$ and $\exp(\theta_{12})$ is a multiplier that is only present on the left-hand side of the above inequality, it is sufficient to show that this inequality holds for $\theta_{12} = 0$. Therefore, we want to prove the following inequality:
\begin{equation}
\label{th3_main_ineq_rewr}
\sum\limits_{x, x' \in \mathcal{X}(S^{-}_m)}g(x|S^{-}_m)  g(x'|S^{-}_m) \exp(2\sum_{i \in S^{-}_m}\theta_{2i} x'_i)\Bigl(\exp(2\sum_{i \in S^{-}_m}\theta_{1i} x'_i) - \exp(2\sum_{i \in S^{-}_m}\theta_{1i} x_i)\Bigr) \geq 0.
\end{equation}
It is easy to see that if $\theta_{2i} = 0$ $\forall i \in S^{-}_m$, then inequality~(\ref{th3_main_ineq_rewr}) is satisfied as it turns into an equality. Therefore, one can prove inequality~(\ref{th3_main_ineq_rewr}) by showing that its left-hand side is increasing in $\theta_{2i}$ for any $i \in S^{-}_m$. Without loss of generality, let us fix $i = m$ and show that the left-hand side of inequality~(\ref{th3_main_ineq_rewr}) is increasing in $\theta_{2m}$. We do this by taking a derivative of this function with respect to $\theta_{2m}$ and checking that it is always nonnegative under the given assumptions. In other words, we prove the following inequality:
\begin{equation*}
\sum\limits_{\substack{x, x' \in \mathcal{X}(S^{-}_m) \\ x'_m = 1}}g(x|S^{-}_m)  g(x'|S^{-}_m) \exp(2\sum_{i \in S^{-}_m}\theta_{2i} x'_i)\Bigl(\exp(2\sum_{i \in S^{-}_m}\theta_{1i} x'_i) - \exp(2\sum_{i \in S^{-}_m}\theta_{1i} x_i)\Bigr) \geq 0.
\end{equation*}
Clearly, the above inequality holds if the following inequalities are true:

\begin{equation}
\label{th3_ind_step_1}
\sum\limits_{\substack{x, x' \in \mathcal{X}(S^{-}_m) \\ x_m = 0, \, x'_m = 1}}g(x|S^{-}_m)  g(x'|S^{-}_m) \exp(2\sum_{i \in S^{-}_m}\theta_{2i} x'_i)\Bigl(\exp(2\sum_{i \in S^{-}_m}\theta_{1i} x'_i) - \exp(2\sum_{i \in S^{-}_m}\theta_{1i} x_i)\Bigr) \geq 0,
\end{equation}
\begin{equation}
\label{th3_ind_step_2}
\sum\limits_{\substack{x, x' \in \mathcal{X}(S^{-}_m) \\ x_m = 1, \, x'_m = 1}}g(x|S^{-}_m)  g(x'|S^{-}_m) \exp(2\sum_{i \in S^{-}_m}\theta_{2i} x'_i)\Bigl(\exp(2\sum_{i \in S^{-}_m}\theta_{1i} x'_i) - \exp(2\sum_{i \in S^{-}_m}\theta_{1i} x_i)\Bigr) \geq 0.
\end{equation}
Let $S^{-}_{m-1} = \{3, \dots, m-1\}$. Then, inequality~(\ref{th3_ind_step_1}) can be rewritten as:
\begin{align*}
\sum\limits_{x, x' \in \mathcal{X}(S^{-}_{m-1}) }&g(x|S^{-}_{m-1})  g(x', |S^{-}_{m-1}) \exp(2\sum_{i \in S^{-}_{m-1}}\theta_{mi} x'_i + \theta_{mm}) \cdot \\
& \exp(2\sum_{i \in S^{-}_{m-1}}\theta_{2i} x'_i + \theta_{2m})\Bigl(\exp(2\sum_{i \in S^{-}_{m-1}}\theta_{1i} x'_i + \theta_{1m}) - \exp(2\sum_{i \in S^{-}_{m-1}}\theta_{1i} x_i)\Bigr) \geq 0.
\end{align*}
It is sufficient to show that the above inequality holds for $\theta_{1m} = 0$, in which case it takes the following form: 
\begin{equation}
\begin{aligned}
\label{th3_ind_step_1_rewr}
\sum\limits_{x, x' \in \mathcal{X}(S^{-}_{m-1}) }&g(x|S^{-}_{m-1})  g(x'|S^{-}_{m-1}) \exp(2\sum_{i \in S^{-}_{m-1}}(\theta_{mi} + \theta_{2i}) x'_i) \cdot \\
& \Bigl(\exp(2\sum_{i \in S^{-}_{m-1}}\theta_{1i} x'_i) - \exp(2\sum_{i \in S^{-}_{m-1}}\theta_{1i} x_i)\Bigr) \geq 0.
\end{aligned}
\end{equation}
Inequality~(\ref{th3_ind_step_1_rewr}) is equivalent to inequality~(\ref{th3_main_ineq_rewr}) where assortment $S^{-}_m$ is replaced with assortment $S^{-}_{m-1}$ and parameters $\theta_{2i}$ are replaced with parameters $\theta_{2i} + \theta_{mi}$ $\forall i \in S^{-}_{m-1}$. Therefore, inequality~(\ref{th3_ind_step_1_rewr}) holds by induction hypothesis.

Similarly, inequality~(\ref{th3_ind_step_2}) is equivalent to the following one:
\begin{equation}
\begin{aligned}
\label{th3_ind_step_2_rewr}
\sum\limits_{x, x' \in \mathcal{X}(S^{-}_{m-1}) }&g(x|S^{-}_{m-1})  g(x', |S^{-}_{m-1}) \exp(2\sum_{i \in S^{-}_{m-1}}\theta_{mi} (x_i + x'_i)) \cdot \\
& \exp(2\sum_{i \in S^{-}_{m-1}}\theta_{2i} x'_i)\Bigl(\exp(2\sum_{i \in S^{-}_{m-1}}\theta_{1i} x'_i) - \exp(2\sum_{i \in S^{-}_{m-1}}\theta_{1i} x_i)\Bigr) \geq 0.
\end{aligned}
\end{equation}
Inequality~(\ref{th3_ind_step_2_rewr}) holds true since it is equivalent to inequality~(\ref{th3_main_ineq_rewr}) where assortment $S^{-}_m$ is replaced with assortment $S^{-}_{m-1}$ and parameters $\theta_{ii}$ are replaced with parameters $\theta_{ii} + \theta_{mi}$ $\forall i \in S^{-}_{m-1}$.
Thus, inequality~(\ref{th3_main_ineq_rewr}) holds for dimension~$m$, which concludes our proof by induction.
\qed

\section{Proof of Proposition~\ref{ising_opt_proper5}}
\label{ising_opt_proper5_proof}
First, from Proposition~\ref{ising_opt_proper3} it follows that we can consider products from the isolated subgraph as a separate product portfolio, i.e., we can assume that $\mathcal{H} = \mathcal{N}$.
Then, our goal is to show that $R(S \cup \{k\}) \geq R(S)$ for any assortment $S \subseteq \mathcal{N}$ such that $k \notin S$. Note that
\begin{align*}
R(S\cup\{k\}) &= \dfrac{\sum\limits_{x \in \mathcal{X}(S\cup\{k\})}\exp\Bigl(\sum\limits_{i \in S\cup\{k\}}\theta_{ii}x_i + \sum\limits_{i \ne j, \ i, j \in S\cup\{k\}}x_i\theta_{ij}x_j\Bigr)\sum\limits_{j \in S\cup\{k\}} r_{j} x_j}{\sum\limits_{x \in \mathcal{X}(S\cup\{k\})}\exp\Bigl(\sum\limits_{i \in S\cup\{k\}}\theta_{ii}x_i + \sum\limits_{i \ne j, \ i, j \in S\cup\{k\}}x_i\theta_{ij}x_j \Bigr)} \\ 
&\geq  \dfrac{\sum\limits_{x \in \mathcal{X}(S\cup\{k\})}\exp\Bigl(\sum\limits_{i \in S\cup\{k\}}\theta_{ii}x_i + \sum\limits_{i \ne j, \ i, j \in S\cup\{k\}}x_i\theta_{ij}x_j\Bigr)\sum\limits_{j \in S} r_{j} x_j}{\sum\limits_{x \in \mathcal{X}(S\cup\{k\})}\exp\Bigl(\sum\limits_{i \in S\cup\{k\}}\theta_{ii}x_i + \sum\limits_{i \ne j, \ i, j \in S\cup\{k\}}x_i\theta_{ij}x_j\Bigr)}
\\ 
&=  \dfrac{\sum\limits_{x \in \mathcal{X}(S)}\bigl(1 + \exp(\theta_{kk} + 2\sum\limits_{i \in S} x_i\theta_{ki})\bigr)\exp\Bigl(\sum_{i \in S}\theta_{ii}x_i + \sum\limits_{i \ne j, \ i, j \in S}x_i\theta_{ij}x_j\Bigr)\sum\limits_{j \in S} r_{j} x_j}{\sum\limits_{x \in \mathcal{X}(S)}\bigl(1 + \exp(\theta_{kk} + 2\sum\limits_{i \in S} x_i\theta_{ki})\bigr)\exp\Bigl(\sum_{i \in S}\theta_{ii}x_i + \sum\limits_{i \ne j, \ i, j \in S}x_i\theta_{ij}x_j \Bigr)}
\end{align*}
and
\begin{equation*}
R(S) =  \dfrac{\sum\limits_{x \in \mathcal{X}(S)}\exp\Bigl(\sum_{i \in S}\theta_{ii}x_i + \sum\limits_{i \ne j, \ i, j \in S}x_i\theta_{ij}x_j\Bigr)\sum\limits_{j \in S} r_{j} x_j}{\sum\limits_{x \in \mathcal{X}(S)}\exp\Bigl(\sum_{i \in S}\theta_{ii}x_i + \sum\limits_{i \ne j, \ i, j \in S}x_i\theta_{ij}x_j \Bigr)}.
\end{equation*}
Let us denote $\exp\Bigl(\sum_{i \in S}\theta_{ii}x_i + \sum\limits_{i \ne j, \ i, j \in S}x_i\theta_{ij}x_j\Bigr)$ by $a(x)$, $\sum\limits_{j \in S} r_{j} x_j$ by $b(x)$, and $\bigl(1 + \exp(\theta_{kk} + 2\sum\limits_{i \in S} x_i\theta_{ki})\bigr)$ by $c(x)$. For readability, let us slightly abuse the notation and write just the summation over $x$ instead of the summation over $x \in \mathcal{X}(S)$. Then, our goal is to prove that the following inequality holds for any assortment $S \subseteq \mathcal{N}$:
\begin{equation*}
    \dfrac{\sum_{x} a(x)b(x)c(x)}{\sum_{x} a(x)c(x)} \geq \dfrac{\sum_{x} a(x)b(x)}{\sum_{x} a(x)},
\end{equation*}
or 
\begin{equation*}
\sum_{x} a(x)b(x)c(x) \sum_{x} a(x) \geq \sum_{x} a(x)b(x) \sum_{x} a(x) c(x).
\end{equation*}
The above inequality can be rewritten as:
\begin{equation*}
\sum_{x, x'} a(x)b(x)c(x) a(x') \geq \sum_{x, x'} a(x)b(x) a(x') c(x'),
\end{equation*}
which is equivalent to 
\begin{equation*}
\sum_{x, x'} a(x)b(x) a(x') (c(x) - c(x')) \geq 0.
\end{equation*}
Let us multiply the latter inequality by 2 and split the sum on the left-hand side into pairs in the following way:
\begin{equation*}
\sum_{x, x'} \Bigl(a(x)b(x) a(x') (c(x) - c(x')) + a(x')b(x') a(x) (c(x') - c(x))\Bigr) \geq 0,
\end{equation*}
which can be rewritten as
\begin{equation}
\label{pr4_last_ineq}
\sum_{x, x'} a(x) a(x') (b(x) - b(x')) (c(x) - c(x')) \geq 0.
\end{equation}
Finally, suppose that $b(x) > b(x')$, i.e., $\sum\limits_{j \in S} r_{j} x_j > \sum\limits_{j \in S} r_{j} x_j'$. Then, $\exp(\theta_{kk} + 2\sum\limits_{j \in S} \alpha r_{j} x_j) > \exp(\theta_{kk} + 2\sum\limits_{j \in S} \alpha r_{j} x_j')$, meaning that $c(x) > c(x')$. Similarly, if $b(x) < b(x')$, then $c(x) < c(x')$ as well. Therefore, inequality~(\ref{pr4_last_ineq}) holds true, which concludes the proof.
\qed

\section{Sensitivity of Benchmark Performance to Basket Size}
\label{sens_bench_perfo_basket}

\begin{figure}[b!]
	\centering
	\includegraphics[clip,width=0.55\linewidth]{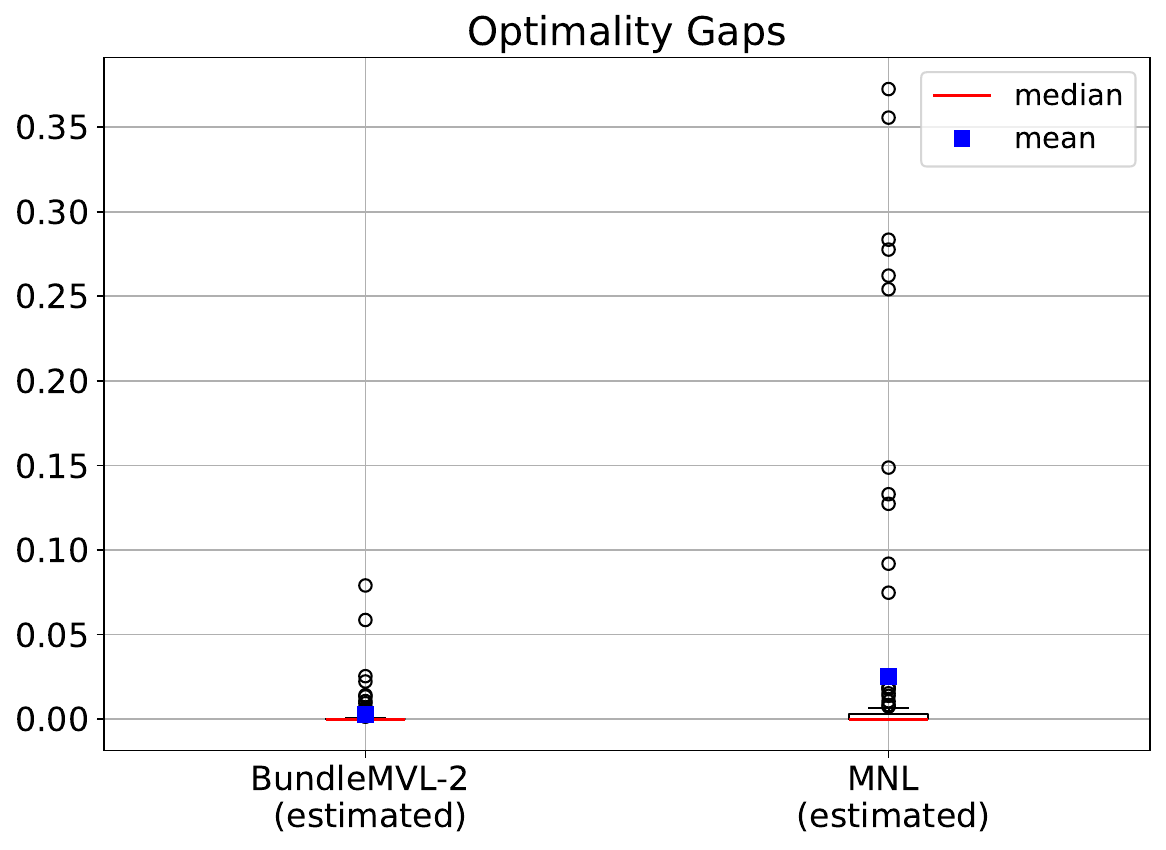}
	\caption{Optimality gap of BundleMVL-2 and MNL with Ising model as the ground truth in small basket size setting.}	
    \label{fig:res2}
\end{figure}

In this appendix, we illustrate how the performance of the considered alternative choice models -- MNL and BundleMVL-2 -- depends on the basket-size distribution induced by the parameter generation procedure (and which reflects different customer purchasing behavior observed in the real world). First, suppose that we modify the procedure from Section~\ref{2subsec:num_analysis} by sampling the values of $\theta_{ii}$ from the interval $[1,2]$ and setting $p_{neg}=0.99$. Under this specification, the vast majority of generated baskets contain at most two products. In particular, the distribution of basket sizes across all sampled baskets and all generated instances has a mean of 1.73 and a standard deviation of 0.58. In this case, the approach based on estimating BundleMVL-2 and optimizing the assortment under this model yields near-optimal solutions, with an average optimality gap of only 0.28\% (see Figure \ref{fig:res2}).
The MNL-based approach also shows good performance, with an average optimality gap of 2.53\%.

\begin{figure}[t!]
	\centering
	\includegraphics[clip,width=0.55\linewidth]{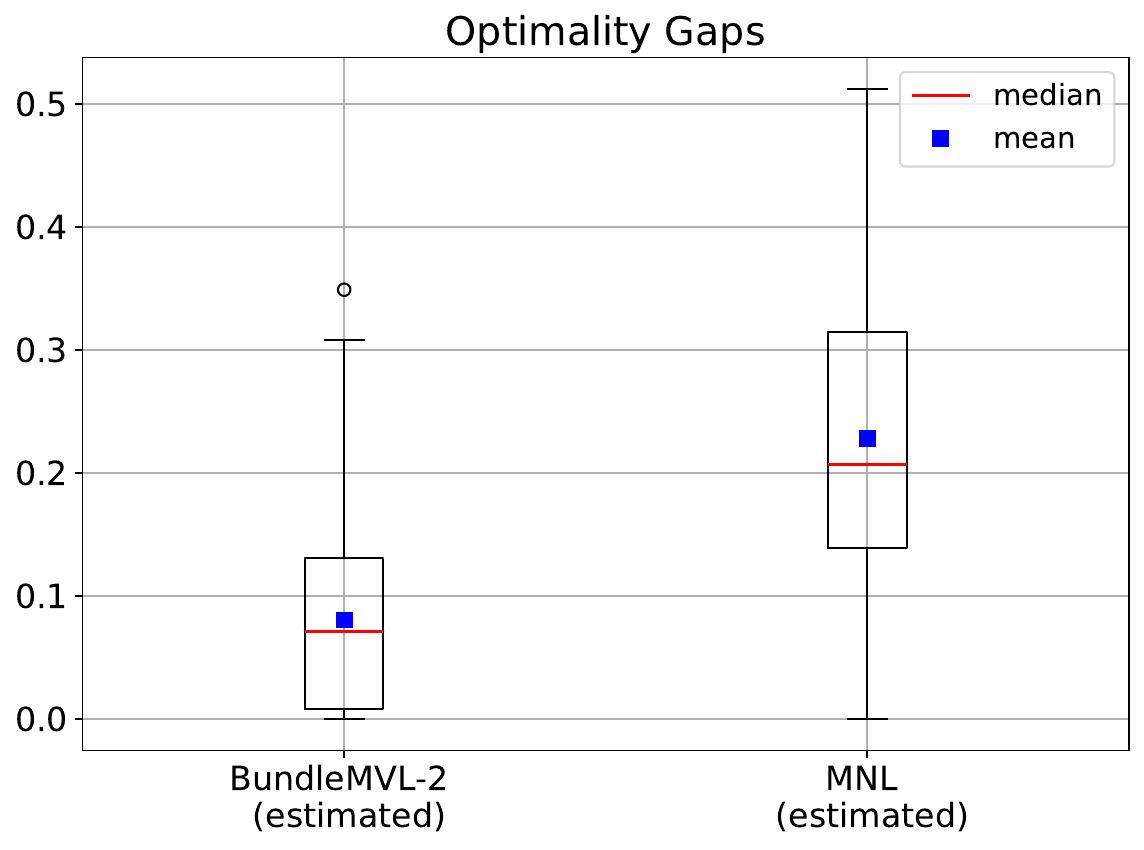}
	\caption{Optimality gap of BundleMVL-2 and MNL with Ising model as the ground truth in large basket size setting.}	\label{fig:res3}
\end{figure}

Next, suppose that we again sample $\theta_{ii}$ from the interval $[1,2]$, but set $p_{\mathrm{neg}}=0.5$. This leads to considerably larger baskets, with the distribution of basket sizes across all sampled baskets and all generated instances having a mean of 6.93 and a standard deviation of 1.17. In this case, the optimality gaps become substantially larger: the average optimality gap equals 8.05\% for BundleMVL-2 and 22.79\% for MNL, as shown by Figure \ref{fig:res3}. 
These additional experiments confirm that benchmark models tailored to small baskets can be highly effective when the underlying choice behavior leads predominantly to small baskets, but their performance deteriorates considerably when customers tend to purchase larger baskets.

\singlespacing
\bibliographystyle{apalike}
\bibliography{literature_ising_ao.bib}

\end{document}